\documentclass[11pt, reqno, oneside]{amsart}

\usepackage{amssymb, amsfonts, mathrsfs, verbatim,graphicx,graphics,epsfig,psfrag,mathtools, bm, bbm}

\DeclareMathAlphabet{\mathpzc}{OT1}{pzc}{m}{it}

\usepackage{mathabx}
\usepackage{slashed}
\usepackage{color, xcolor}

\usepackage{microtype}

\usepackage{abstract}
\usepackage[all]{xy}

\usepackage[left= 3.4 cm, right= 3.4 cm, top= 2.8 cm, bottom=2.7 cm, foot= 1.2 cm, marginparwidth=2.7 cm, marginparsep=0.5 cm]{geometry}

\usepackage[inline]{enumitem}

\usepackage{hyperref}
\hypersetup{colorlinks}
\definecolor{darkred}{rgb}{0.5,0,0}
\definecolor{darkgreen}{rgb}{0,0.5,0}
\definecolor{darkblue}{rgb}{0,0,0.5}
\hypersetup{colorlinks, linkcolor=darkblue, filecolor=darkgreen, urlcolor=darkred, citecolor=darkblue}
\makeatletter 
\@addtoreset{equation}{section}
\makeatother  

\numberwithin{equation}{section}

\usepackage{tikz-cd}
\usepackage{tikz}
\usetikzlibrary{arrows}
\usepackage{float}

\usetikzlibrary{decorations.pathreplacing}
\usepackage{etex}

\newtheorem{thm}{Theorem}[section]

\newtheorem{prop}[thm]{Proposition}
\newtheorem{lemma}[thm]{Lemma}

\theoremstyle{definition}
\newtheorem{defn}[thm]{Definition}
\theoremstyle{remark}
\newtheorem{rem}[thm]{Remark}
\newtheorem{hyp}[thm]{Hypothesis}

\newtheorem{example}[thm]{Example}
\newtheorem{notation}[thm]{Notation}
\newtheorem{convention}[thm]{Convention}

\newcounter{notes}
{\end{list}}

\makeatletter
\@addtoreset{thm}{section}
\makeatother

\newcommand\qu{/\kern-.7ex/} 
\renewcommand{\setminus}{\smallsetminus}

\newcommand{\beq}{\begin{equation}}
\newcommand{\eeq}{\end{equation}}
\newcommand{\beqn}{\begin{equation*}}
\newcommand{\eeqn}{\end{equation*}}
\newcommand{\ov}{\overline}
\newcommand{\mb}{\mathbb}
\newcommand{\mc}{\mathcal}
\newcommand{\mf}{\mathfrak}

\newcommand{\Gammait}{{\mathit{\Gamma}}}

\renewcommand{\i}{{\bm i}}

\newcommand{\uds}[1]{\underline{\smash{#1}}{}}

\newcommand{\ev}{{\rm ev}}

\newcommand{\s}{\slashed}

\begin{document}

\author[Xu]{Guangbo Xu}
\address{
Department of Mathematics\\
Texas A{\&}M University\\
College Station, Texas 77843, USA
}
\email{guangboxu@math.tamu.edu}

\date{\today}

\title{Gluing Affine Vortices}

\maketitle

\begin{abstract}
We construct a gluing map for stable affine vortices over the upper half plane with the Lagrangian boundary condition at a rigid, regular, codimension one configuration. This construction plays an important role in establishing the relation between the gauged linear sigma model and the nonlinear sigma model in the presence of Lagrangian branes. 


\end{abstract}

\setcounter{tocdepth}{2}
\tableofcontents

\section{Introduction}

The vortex equation is a first-order elliptic equation whose solutions, called {\it vortices}, appear in many different areas of mathematics and physics. In physics, vortices first showed up in the Ginzburg--Landau theory in superconductivity (see \cite{Ginzburg_Landau_1950} \cite{Taubes_1990} \cite{Jaffe_Taubes}), and reappeared in many other more abstract physical theories such as the gauged linear sigma model (GLSM) \cite{Witten_LGCY}. In mathematics, especially in symplectic geometry, vortices play the role as equivariant generalizations of pseudoholomorphic curves. Enumerations of vortices lead to the definition of the gauged Gromov--Witten invariants (see \cite{Mundet_thesis, Cieliebak_Gaio_Salamon_2000, Mundet_2003, Cieliebak_Gaio_Mundet_Salamon_2002}) and Floer-type theories (see \cite{Frauenfelder_thesis} \cite{Xu_VHF}). The generalization of the vortex equation, called the gauged Witten equation, plays a fundamental role in the author's project with Tian on a mathematical theory of the GLSM (see \cite{Tian_Xu, Tian_Xu_2, Tian_Xu_3, Tian_Xu_2017, Tian_Xu_geometric}). 

An important argument in the vortex theory is the adiabatic limit. In the adiabatic limit, vortices converge to pseudoholomorphic curves. This phenomenon leads several important results and conjectures. Using the adiabatic limit, Gaio--Salamon \cite{Gaio_Salamon_2005} relates certain gauged Gromov--Witten invariants and the usual Gromov--Witten invariants. Moreover, they conjecture that a more general relation between these two types of invariants should lead to a quantization of the Kirwan map in a similar way as the Gromov--Witten invariants quantize the multiplicative structure of the cohomology ring. In gauge theory the adiabatic limit argument for the instanton equation leads to Dostoglou--Salamon's proof of the $SO(3)$ Atiayh--Floer conjecture for mapping cylinders \cite{Dostoglou_Salamon}. 

In the adiabatic limit the convergence of vortices towards holomorphic curves generally should be ``corrected'' by taking into account the contribution of certain bubbles. These bubbles are called {\it affine vortices} or pointlike instantons.\footnote{Pointlike instantons usually refer to the bubble in the GLSM when there is a nonzero superpotential. They satisfy a more general equation called the gauged Witten equation.} Their contributions, often nontrivial, have already been pointed out by Witten \cite{Witten_LGCY} in the context of the GLSM. Similar to counting pseudoholomorphic curves, to define the contributions of affine vortices, one needs to compactify their moduli spaces and construct manifold-like structures (such as Kuranishi structures) over the moduli space. The compactness results for affine vortices (see \cite{Ziltener_thesis, Ziltener_book} and \cite{Wang_Xu}) and the stratification of the domain moduli (see \cite{Mau_Woodward_2010}) indicate that the involved gluing construction should differ from the gluing of holomorphic curves. 

In this paper we construct the gluing map for affine vortices over the upper half plane satisfying a Lagrangian boundary condition. Instead of working in the most general situation, we restrict to the case when the singular configuration lies in a codimension one stratum of the compactified moduli space and when the singular configuration is rigid. A practical reason for having this restriction is that this is what one needs in the related work \cite{Woodward_Xu} where the authors construct an open version of the quantum Kirwan map. The construction in more general situations will be considered elsewhere. For example, in the forthcoming work \cite{Limit2} a general gluing construction will be provided for pointlike instantons (with nontrivial superpotentials but without a boundary condition) over the complex plane, giving local Kuranishi model on the relevant moduli space.

Our main theorem is the following, whose precise version is restated as Theorem \ref{mainthm}.

\begin{thm}
Let $(V, \omega, \mu)$ be a Hamiltonian $K$-manifold. Assume that the symplectic quotient $X: = \mu^{-1}(0)/K$ is a free quotient. Let $L_V \subset V$ be a $K$-invariant embedded Lagrangian submanifold which is contained in $\mu^{-1}(0)$. 

Given nonnegative integers $l, l^+$ with $l + l^+ \geq 1$, let $\ov{\mc M}_{l, l^+}(V, L_V)$ be the moduli space of gauge equivalence classes of stable affine vortices over ${\mb H}$ with $l$ boundary marked points and $l^+$ interior marked points (with respect to a given family of domain-dependent almost complex structure), equipped with a natural topology. Let $\Gamma$ be a simple combinatorial type (see Definition \ref{simple}) which labels a codimension one stratum ${\mc M}_\Gamma (V, L_V) \subset \ov{\mc M}_{l, l^+}(V, L_V)$. 

Given a regular and rigid point ${\bm p} \in {\mc M}_\Gamma(V, L_V)$ (see Definition \ref{defn32}), there exists an open neighborhood of ${\bm p}$ in $\ov{\mc M}_{l, l^+}(V, L_V)$ which is homeomorphic to the interval $[0, \epsilon)$.
\end{thm}

Our proof is conceptually straightforward and follows a standard protocol used in gauge theory and symplectic geometry. However the complicated behavior of vortices and some subtle features make the argument rather involved. On the technical level, one can view our construction as a nontrivial generalization of the case of Gaio--Salamon \cite{Gaio_Salamon_2005} in two different ways. First, the singular configuration we would like to glue consists of two types of components, affine vortices or holomorphic curves, while in the case of \cite{Gaio_Salamon_2005} they only need to treat limiting configurations which are smooth holomorphic curves. Second, while Gaio--Salamon only treated the vortex equation over compact domains, the noncompact domain of the vortex equation considered here requires particular choices of weighted Sobolev norms and a slightly atypical local model theory (see \cite{Venugopalan_Xu}). With these two types of complications combined, to ensure the gluing protocol can be operated, we need to make extra effort to arrange various ingredients more carefully.

\subsection{Extensions and applications}

As we have explained, the immediate motivation for studying the gluing of affine vortices is from the project of the author with Woodward \cite{Woodward_Xu}, aiming at defining an ``open quantum Kirwan map.'' The open quantum Kirwan map is a morphism of $A_\infty$ algebras from the quasimap Fukaya $A_\infty$ algebra of $L_V$ to a bulk-deformed Fukaya $A_\infty$ algebra of the Lagrangian in the quotient $X$. 

Moreover, using the technique and the analytical setting of this paper, one can construct the gluing map for affine vortices over ${\mb C}$, and the gluing map with respect to the adiabatic limit. This would be an important step towards the resolution of Salamon's quantum Kirwan map conjecture in the symplectic setting, initiated in \cite{Ziltener_thesis} \cite{Ziltener_book} (proved by Woodward \cite{Woodward_15} in the algebraic case). If one allows a nontrivial superpotential in the general setting of the GLSM, where affine vortices are often referred to as pointlike instantons, then the corresponding gluing construction is necessary to establish a relation between the GLSM correlation functions and Gromov--Witten invariants. This is part of the on-going project of the author with G. Tian \cite{Limit2}.

One specific feature of the affine vortex equation is that the equation has only translation invariance but not conformal invariance. In symplectic geometry and gauge theory there are other types of objects which have the same symmetry types. The figure-eight bubble, appeared in the strip shrinking limits of pseudoholomorphic quilts (see \cite{Wehrheim_Woodward_2012} \cite{Wehrheim_Woodward_2015} \cite{Bottman_Wehrheim_2018}), is such an example. There are also examples in gauge theory, such as the anti-self-dual equation or the monopole equation over ${\mb C} \times \Sigma$ (see \cite{Wehrheim_2006} and \cite{Donghao_2018}) and over the product of the real line with a noncompact three-manifold with cylindrical ends (see \cite{Xu_2019}). We hope that the technique of this paper can be used in the gluing construction for other translation invariant equations.

\subsection{Acknowledgments}

I would like to thank Chris Woodward for many stimulating discussions and his generosity in sharing ideas. I would also like to thank Gang Tian and Kenji Fukaya for their support and encouragements, and thank Sushmita Venugopalan for helpful discussions.

\section{Affine Vortices and Their Moduli}

In this section we review the basic facts about affine vortices and their moduli spaces. We also set up notations for trees, domain-dependent almost complex structures, and local models of the domain moduli.
	
\subsection{Preliminaries on affine vortices}

We first recall basic knowledge of affine vortices and fix some notations. Let $K$ be a compact Lie group with its Lie algebra ${\mf k}$ and its complexification $G$. Let $(V, \omega_V)$ be a symplectic manifold with a smooth $K$-action. For each $\eta \in {\mf k}$, let ${\mc X}_\eta\in \Gammait (TV)$ be the infinitesimal action. Our convention is that the map $\eta \mapsto {\mc X}_\eta$ is an anti-homomorphism of Lie algebras from ${\mf k}$ to $\Gammait 	(TV)$. Assume that the $K$-action is Hamiltonian and has a moment map
\beqn
\mu: V \to {\mf k}^*.
\eeqn
This means that $\mu$ is $K$-equivariant and satisfies 
\beqn
\langle d\mu( Z ), \eta \rangle = \omega_V ( {\mc X}_\eta,Z),\ \forall \eta \in {\mf k},\ Z \in TV.
\eeqn

We make the following basic assumption.

\begin{hyp}\label{hyp21}
$0 \in {\mf k}$ is a regular value of $\mu$ and the symplectic quotient $X:= \mu^{-1}(0)/K$ is a free quotient.
\end{hyp}

Now we recall the notion of gauged maps and vortices. Let $\Sigma$ be a Riemann surface with possibly nonempty boundary. A {\it gauged map} from $\Sigma$ to $V$ is a triple ${\bm v} = (P, A, u)$ where $P \to \Sigma$ is a smooth $K$-bundle, $A \in {\mc A}(P)$ is a connection, and $u \in \Gammait (P(V))$ is a smooth section. Here $P(V) =  (P \times V)/K$ is the associated fibre bundle. To a gauged map ${\bm v} = (P, A, u)$ one has the following associated objects: 
\begin{itemize}

\item the {\it curvature} $F_A \in \Omega^2(\Sigma, {\rm ad}P)$,

\item the {\it moment potential} $\mu(u) \in \Omega^0( \Sigma, ({\rm ad} P)^*)$, and

\item the {\it covariant derivative} $d_A u \in \Omega^1(\Sigma, u^* P(TV))$.
\end{itemize}
Suppose $\Sigma$ is equipped with an area form $\nu_\Sigma$, which, together with the complex structure on $\Sigma$, determines a metric on $\Sigma$. Suppose $V$ is equipped with a $K$-invariant Riemannian metric and the Lie algebra ${\mf k}$ is equipped with an invariant inner product. We define the {\it energy} of the gauged map ${\bm v} = (P, A, u)$ to be
\beqn
E({\bm v}) = E(P, A, u):= \frac{1}{2} \left( \| d_A u \|_{L^2(\Sigma)}^2 + \| F_A \|_{L^2(\Sigma)}^2 + \| \mu(u) \|_{L^2(\Sigma)}^2 \right). 
\eeqn

The vortex equation is the equation of energy minimizers. Suppose the Riemannian metric on $V$ used to define the energy is determined by $\omega_V$ and a $K$-invariant $\omega_V$-tamed almost complex structure $J$. Then gauged maps minimizing the energy (locally) solve the {\it symplectic vortex equation}
\begin{align}\label{vortex}
&\ \ov\partial_A u = 0,\ &\ * F_A + \mu(u) = 0.
\end{align}
Here $\ov\partial_A u:= (d_A u)^{0,1}$ is the $(0,1)$-part of the covariant derivative with respect to the domain complex structure on $\Sigma$ and the target almost complex structure $J$ on $V$; $*$ is the Hodge star operator on $\Sigma$ determined by the metric on $\Sigma$. The second equation also uses the metric on the Lie algebra to identify ${\rm ad} P \cong  ({\rm ad} P)^*$. We call a solution to \eqref{vortex} a {\it symplectic vortex} (or simply a {\it vortex}) on $\Sigma$. 

We will impose a Lagrangian boundary condition for the vortex equation. Let 
\beqn
L \subset X
\eeqn
be a Lagrangian submanifold, meaning that $\omega_X|_L \equiv 0$, where $\omega_X$ is the symplectic form of the symplectic reduction $X$ induced from $\omega_V$. Then the preimage 
\beqn
L_V \subset \mu^{-1}(0)
\eeqn
under the map $\mu^{-1}(0) \to \mu^{-1}(0)/K \subset V$ is a $K$-invariant Lagrangian submanifold of $V$. When $\partial \Sigma \neq \emptyset$, we always impose the following condition for the vortex equation \eqref{vortex}
\beq\label{boundary}
u(\partial \Sigma) \subset P(L_V) \subset P(V).
\eeq

The vortex equation and the boundary condition are invariant under gauge transformations. A {\it gauge transformation} on a $K$-bundle $P \to \Sigma$ is a section $g\in \Gammait (P(K))$ where $P(K) \to \Sigma$ is the bundle associated to the adjoint action of $K$ on itself. Viewing them as automorphisms of principal bundles, gauge transformations can pull-back connections $A \in {\mc A}(P)$ and sections $u \in \Gammait (P(V))$. We use the notation such that the action by gauge transformations is a left action: for a gauged map ${\bm v} = (P, A, u)$, denote
\beqn
g\cdot {\bm v} = (P, A \circ g^{-1},\ u \circ g^{-1}).
\eeqn

In this paper we only consider gauged maps defined over rather simple domains. These domains are the complex plane ${\mb C}$, the upper half plane ${\mb H}$, or their open subsets. Over such a domain there is a standard flat metric and standard complex coordinate $s + {\bm i} t$. Let ${\mb A}$ denote either ${\mb C}$ or ${\mb H}$. Then a $K$-bundle $P \to {\mb A}$ is always trivial; with respect to a trivialization, a gauged map from ${\mb A}$ to $V$ can be identified with a triple
\beqn
{\bm v} = (u, \phi, \psi): {\mb A} \to V \times {\mf k} \times {\mf k}.
\eeqn
Here we identify the connection $A$ with $d + \phi ds + \psi dt$ with respect to the trivialization. In such a gauge, the vortex equation \eqref{vortex} is equivalent to 
\begin{align}\label{affine}
&\ \partial_s u + {\mc X}_\phi (u) + J (\partial_t u + {\mc X}_\psi (u)) = 0,\ &\ \partial_s \psi - \partial_t \phi + [\phi, \psi] + \mu(u) = 0
\end{align}
and the boundary condition \eqref{boundary} is equivalent to 
\beq\label{boundary2}
u(\partial {\mb A}) \subset L_V.
\eeq

\begin{defn}
An {\it affine vortex} over ${\mb C}$ (resp. ${\mb H}$) is a finite energy solution ${\bm v} = (u, \phi, \psi)$ to \eqref{affine} subject to the boundary condition \eqref{boundary2}. Two affine vortices over ${\mb C}$ resp. ${\mb H}$ are {\it isomorphic} if after a complex resp. real translation they are gauge equivalent.
\end{defn}

One can define a natural topology on the set of isomorphism classes of affine vortices. One can see from the two examples below that this topology cannot be compact. 

\begin{example}
Taubes \cite{Taubes_vortex} classifies affine vortices for the simplest target space. Given a complex polynomial $f(z)$, Taubes shows that there exists a unique solution $h: {\mb C} \to {\mb R}$ to the Kazdan--Warner equation 
\beqn
- \frac{\Delta h}{2\pi} + \frac{1}{2} \left( e^{2h} |f(z)|^2 - 1 \right) = 0
\eeqn
with an appropriate asymptotic condition on $h$. It follows that ${\bm v}_f: = (e^h f, -\partial_t h, \partial_s h)$ is an affine vortex with target $V = {\mb C}$ acted by $K = U(1)$ with moment map 
\beqn
\mu(z) = - \frac{{\bm i}}{2} (|z|^2 - 1).
\eeqn

One can see how these affine vortices degenerate by looking at the zeroes of the polynomials. Consider a sequence of monic polynomials of degree $d = 2$
\beqn
f_k (z) = (z - z_k) (z-z_k').
\eeqn
Suppose as $k \to \infty$, $|z_k - z_k'| \to \infty$. Then as $k \to \infty$, one can see that the sequence of affine vortices ${\bm v}_{f_k}$ split into two components. One can viewed such a degeneration as a complex analogue of breaking trajectories in Morse theory.
\end{example}

\begin{example}
In general, besides splitting into several components, affine vortices can also converge to nontrivial holomorphic maps into the symplectic quotient. One can see this type of limiting behavior in the higher-rank generalization of abelian vortices. For an $N$-tuple of polynomials
\beqn
\vec{f} = (f^1, \ldots, f^N),\ \max_\alpha \{ {\rm deg} f^\alpha\} = d\geq 1,
\eeqn
by the main theorem of \cite{Guangbo_vortex} (see also \cite{VW_affine}), there exists a unique solution $h: {\mb C} \to {\mb R}$ to the Kazdan--Warner equation 
\beqn
- \frac{\Delta h}{2\pi} + \frac{1}{2} \left( e^{2h} \sum_{\alpha=1}^N | f^\alpha(z)|^2 - 1 \right) = 0.
\eeqn
So such an $N$-tuple of polynomials corresponds to an affine vortex with target $V = {\mb C}^N$ acted by the  gauge group $K = U(1)$.

A sequence of affine vortices could converge to a holomorphic sphere in the quotient $X = \mb{CP}^{N-1}$. Consider the simplest case when $N = 2$, $d= 1$. Consider a sequence of pairs of polynomials
\begin{align*}
&\ f_k^1 (z) = z- k,\ &\  f_k^2(z) = z,\ \ \  k = 1, 2, \cdots.
\end{align*}
For the corresponding sequence of affine vortices ${\bm v}_k = (u_k, \phi_k, \psi_k)$, one can argue that in the $k \to \infty$ limit, if one reparametrize the domain using factor $k$, then $\mu(u_k)$ converges to zero uniformly and the induced sequence of maps $\bar u_k: {\mb C} \to \mb{CP}^1$ converges to the holomorphic map $z \mapsto [z- 1, z]$ into $\mb{CP}^1$.
\end{example}

We will describe compactifications of moduli spaces of affine vortices in the next subsection. A description relies on the following property of affine vortices. 

\begin{prop}\cite{Gaio_Salamon_2005}\cite{Ziltener_Decay}\cite{Wang_Xu}
Given an affine vortex ${\bm v} = (u, \phi, \psi)$ over ${\mb C}$ resp. ${\mb H}$, there is a $K$-orbit $Kx$ in $\mu^{-1}(0)$ resp. in $L_V$ such that
\beqn
\lim_{z \to \infty} K u(z) = Kx \in V/K.
\eeqn
\end{prop}

The above proposition allows on to define the {\it evaluation at infinity} of affine vortices. Given an affine vortex ${\bm v} = (u, \phi, \psi)$ over ${\mb C}$ resp. ${\mb H}$, its limit at the infinity is denoted by
\beq\label{evaluation}
{\rm ev}_\infty({\bm v}) \in X\ \ {\rm resp.}\ \ {\rm ev}_\infty({\bm v}) \in L.
\eeq

\subsection{Scaled trees and domain moduli}

In this section we provide the detailed combinatorial ingredients needed to compactify the moduli space of affine vortices. We first fix the notations about trees. A tree $\Gamma$ has a set of vertices $V_\Gamma$ and a set of edges $E_\Gamma$. We allow both {\it finite edges}, which connect two vertices, and {\it semi-infinite edges}, which are only attached to one vertex. Let 
\beqn
E_\Gamma = E_\Gamma^- \cup E_\Gamma^\infty
\eeqn
be the decomposition into finite and semi-infinite edges. We always assume that the tree has a distinguished vertex $v_0 \in V_\Gamma$ called the {\it root}.\footnote{In some related works, we have the convention that a distinguished semi-infinite edge called the {\it output} is always attached to the root. However in this paper we do not use the output.} The root induces a partial order among vertices: we write $v_\alpha \leq v_\beta$ if $v_\alpha$ is closer to the root. We denote $v_\beta \succ v_\alpha$ if $v_\alpha \leq v_\beta$ and these two vertices are adjacent. One can view a tree as a 1-complex in which the semi-infinite edges are open cells. A {\it ribbon tree} is a tree $\Gamma$ together with an isotopy class of embeddings $\Gamma \hookrightarrow {\mb R}^2$. 

\begin{defn}[Based tree]
A {\it based tree} is a rooted tree $\Gamma$ together with a subtree $\uds\Gamma$ containing the root, a ribbon tree structure on the base, and a bijective labeling  
\beqn
\{ 1, \ldots, l^+\} \cong E_\Gamma^\infty \setminus E_{\uds \Gamma}^\infty
\eeqn
of semi-infinite edges not in the base.
\end{defn}
\noindent Notice that semi-infinite edges in the base (representing boundary marked points) are canonically ordered (counterclockwise) and labelled by integers $1, \ldots, l$ via the ribbon tree structure of the base. A based tree can be used to model the combinatorial type of a stable holomorphic disk such that vertices in the base represent disk components and vertices not in the base represent sphere components.

\begin{notation} Let $\Gamma$ be a rooted tree. 
\begin{enumerate}
\item For any edge $e \in E_\Gamma$, if $e$ is a finite edge, then let $v_{\alpha(e)} \in V_\Gamma$ denote the vertex on one the side of $e$ which is closer to the root; if $e$ is a semi-infinite edge, then let $v_{\alpha(e)} \in V_\Gamma$ denote the vertex to which $e$ is attached. 

\item Let $V_{\Gamma^{\rm sup}} \subset V_\Gamma$ denote the set of vertices not in the base. 

\end{enumerate}
\end{notation}

To model the combinatorial types of stable affine vortices one needs an extra structure on the tree. A {\it scaling} on a rooted tree $\Gamma$ is a map 
\beqn
{\mf s}: V_\Gamma \to \{0, 1, \infty\}
\eeqn
satisfying the following conditions.
\begin{enumerate}
\item ${\mf s}$ is order-reversing. 

\item For any path $v_1\succ \cdots \succ v_k$ in $\Gamma$, if ${\mf s}(v_1) \leq 1$, ${\mf s}(v_k) \geq 1$, then there is exactly one vertex $v_l$ in this path with ${\mf s}(v_l) = 1$. 
\end{enumerate}
A tree with a scaling is called a {\it scaled tree}. We often abbreviate $(\Gamma, {\mf s})$ by $\Gamma$. Denote
\beqn
V_\Gamma^0 = {\mf s}^{-1}(0),\ \hspace{0.2cm}\ V_\Gamma^1 = {\mf s}^{-1}(1),\ \hspace{0.2cm} \ V_\Gamma^\infty = {\mf s}^{-1}(\infty).
\eeqn
Vertices in these subsets will represent holomorphic curves in $V$, affine vortices, and holomorphic curves in the quotient $X$ respectively.

\begin{figure}[ht]
\centering

\includegraphics[scale=1]{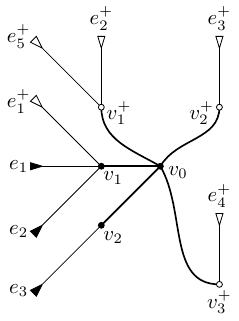}

\caption{A picture of a stable scaled tree. Its base has three vertices $v_0$, $v_1$, and $v_2$. The vertices $v_1^+$, $v_2^+$, and $v_3^+$ are not in the base. We declare that the root $v_0$ has scale $\infty$ and other vertices have scale $1$. The three black incoming arrows are semi-infinite edges in the base and the other white incoming arrows are semi-infinite edges not in the base. The ribbon tree structure of the base is induced from the way we draw the picture. This tree is a simple tree in the sense of Definition \ref{simple}. This tree can model a stable affine vortex with one holomorphic disk component, two affine vortex components over ${\mb H}$, three affine vortex components over ${\mb C}$, three boundary marked points, and five interior marked points.}\label{figure2}
\end{figure}

\begin{defn}
Let $(\Gamma, {\mf s})$ be a scaled tree. A vertex $v_\alpha \in V_\Gamma$ is called {\it stable} if the following conditions are satisfied.
\begin{enumerate}
\item If $v_\alpha \in V_{\Gamma^{\rm sup}}^0 \cup V_{\Gamma^{\rm sup}}^\infty$, then 
\beqn
\# \{ e\in E_\Gamma\ |\ \alpha(e) = \alpha \} \geq 2.
\eeqn

\item If $v_\alpha \in V_{\uds\Gamma}^0 \cup V_{\uds\Gamma}^\infty$, then 
\beqn
2 \# \{ e \in E_\Gamma \setminus E_{\uds\Gamma} |\ \alpha(e) = \alpha\} +  \# \{ e \in E_{\uds\Gamma} \ |\ \alpha(e) = \alpha \} \geq 2.
\eeqn

\item If $v_\alpha \in V_{\Gamma}^1$, then 
\beqn
\# \{ e \in E_\Gamma\ |\ \alpha(e) = \alpha \} \geq 1.
\eeqn
\end{enumerate}
Otherwise, $v_\alpha$ is called {\it unstable}. The scaled tree is called {\it stable} if all vertices are stable. 

See Figure \ref{figure2} for an illustration of a stable scaled tree. 
\end{defn}

\begin{convention}
In this paper we impose the following conditions on scaled trees. We require that all semi-infinite edges are attached to vertices in $V_\Gamma^0 \cup V_\Gamma^1$. 
\end{convention}

\begin{defn}\label{domain} Let $(\Gamma, {\mf s})$ be a scaled tree. A {\it marked scaled curve} of type $(\Gamma, {\mf s})$ is a collection 
\beqn
{\mc C} = \left( (\Sigma_\alpha, \sigma_\alpha, f_\alpha)_{v_\alpha \in V_\Gamma},\ (z_e)_{e \in E_\Gamma} \right)
\eeqn
where 
\begin{itemize}

\item For each $v_\alpha \in V_{\uds\Gamma}$ resp. $V_{\Gamma^{\rm sup}}$, $\Sigma_\alpha$ is a bordered resp. closed Riemann surface, $\sigma_\alpha \in \Omega^2(\Sigma_\alpha) \cup \{\infty\}$, and $f_\alpha: \Sigma_\alpha \cong {\mb H}$ resp. $f_\alpha: \Sigma_\alpha \cong {\mb C}$ is a biholomorphism such that ${\mf s}(v_\alpha) f_\alpha^* ds dt = \sigma_\alpha$.

\item For each $e \in E_\Gamma$, $z_e \in \Sigma_{\alpha(e)}$. 
\end{itemize}
They satisfy the following conditions.
\begin{enumerate}

\item If $e \in E_{\uds\Gamma}$, then $z_e \in \partial \Sigma_\alpha$; otherwise $z_e \in {\rm Int} \Sigma_\alpha$.

\item For each vertex $v_\alpha$, the collection of points defining the set
\beqn
Z_\alpha:= \{ z_e\ |\ \alpha(e) = \alpha \}
\eeqn
are distinct. 
\end{enumerate}
The marked scaled curve is {\it stable} if the tree $(\Gamma, {\mf s})$ is stable. 
\end{defn}

\begin{defn}\label{defn211}
Let ${\mc C} = ((\Sigma_\alpha, \sigma_\alpha, f_\alpha), (z_e))$ and ${\mc C}' = ( (\Sigma_\alpha', \sigma_\alpha', f_\alpha'), (z_e'))$ be two marked scaled curve of type $\Gamma$. An {\it isomorphism} from ${\mc C}$ to ${\mc C}'$ is a collection of biholomorphism 
\beqn
(\varphi_\alpha: (\ov\Sigma_\alpha, \infty) \to (\ov\Sigma_\alpha, \infty))
\eeqn
satisfying the following conditions
\begin{enumerate}

\item For each $v_\alpha \in V_\Gamma$, one has $f_\alpha = f_\alpha' \circ \varphi_\alpha$. 

\item For each $e \in E_\Gamma$, one has $\varphi_{\alpha(e)} (z_e) = z_e'$. 
\end{enumerate}
\end{defn}

We say a marked scaled curve is smooth if the underlying tree $\Gamma$ has a single vertex $v_0$ with ${\mf s}(v_0) = 1$. In this case we can represent a smooth marked scaled curve as
\beqn
\big( ( {\mb H}, ds \wedge dt, {\rm Id}), {\bm z} = (z_e)_{e \in E_\Gamma^\infty} \big)
\eeqn
where $(z_e)_{e\in E_\Gamma^\infty}$ is a configuration of distinct points in the upper half plane. Since semi-infinite edges are labeled by integers, the point configuration is also denoted by 
\beqn
{\bm z} = (z_1, \ldots, z_l; z_1^+, \ldots, z_{l^+}^+)
\eeqn
where $l = \# E_{\uds\Gamma}^\infty$ and $l^+ = \#(E_\Gamma^\infty \setminus E_{\uds\Gamma}^\infty)$. By Definition \ref{defn211}, two smooth marked scaled curves with $l$ boundary markings and $l^+$ interior markings are isomorphic if the corresponding point configurations differ by a real translation in the upper half plane. 

It is easy to see that the automorphism group of a marked scaled curve is finite if and only if it is stable. In this case the automorphism group is actually trivial (as we are in genus zero). For each stable scaled tree $\Gamma$, let ${\mc M}_\Gamma$ be the set of isomorphism classes of marked scaled curves of type $\Gamma$. When $\Gamma$ is stable has a single vertex $v_0$ with ${\mf s}(v_0)=1$ (as described above), then the moduli space ${\mc M}_\Gamma$ can be identified with the space of point configurations modulo translation, and hence is a manifold. For a general stable $\Gamma$, the moduli space ${\mc M}_\Gamma$ has a topology induced from the topology of the space of configurations of points on different components.

Denote by ${\rm Tree}(l, l^+)$ the set of isomorphism classes of scaled trees with $l$ boundary semi-infinite edges and $l^+$ interior semi-infinite edges. By abuse of language, we regard each element of ${\rm Tree}(l, l^+)$ as a tree rather than an isomorphism class. Let ${\rm Tree}^s(l, l^+) \subset {\rm Tree}(l, l^+)$ be the subset of stable trees. Then define 
\beqn
\ov{\mc M}_{l, l^+} = \bigsqcup_{ \Gamma \in {\rm Tree}^s(l, l^+)} {\mc M}_\Gamma. 
\eeqn
A natural topology can be defined on $\ov{\mc M}_{l, l^+}$ stratified by ${\mc M}_\Gamma$. By the work of \cite{Mau_Woodward_2010} the moduli space $\ov{\mc M}_{l, l^+}$ has the topology of a CW complex. We do not explore the detailed structure of this moduli space in general, though the structure near codimension one strata will be  discussed later.

\subsection{Moduli spaces of stable affine vortices}

Now we describe the Gromov--Uhlenbeck compactification of the moduli space of affine vortices. For affine vortices over ${\mb C}$, such a compactification was firstly given by Ziltener \cite{Ziltener_thesis, Ziltener_book}. 

We first define the notion of stable vortices with respect to a fixed almost complex structure $J$. Later we will allow domain-dependent almost complex structures. Let $I$ be the induced almost complex structure on the symplectic quotient $X$. 

\begin{defn}\label{stablevortex} Let $\Gamma$ be a scaled tree. A stable affine vortex of type $\Gamma$ is a collection 
\beqn
{\mc V} = ({\mc C}, ({\bm v}_\alpha)_{v_\alpha\in V_\Gamma} ) = \Big( (\Sigma_\alpha, \sigma_\alpha, f_\alpha)_{v_\alpha\in V_\Gamma}, (z_e)_{e \in E_\Gamma},  ({\bm v}_\alpha)_{v_\alpha \in V_\Gamma}  \Big).
\eeqn
Here
\begin{itemize}

\item ${\mc C} = ((\Sigma_\alpha, \sigma_\alpha, f_\alpha)_{v_\alpha \in V_\Gamma}, (z_e)_{e\in E_\Gamma} )$ is a marked scaled curve of type $\Gamma$. 

\item For each $v_\alpha \in V_\Gamma^0$, ${\bm v}_\alpha$ is a $K$-orbit of holomorphic maps $u_\alpha: \Sigma_\alpha \to V$ satisfying the boundary condition $u_\alpha(\partial \Sigma_\alpha) \subset L_V$.

\item For each $v_\alpha \in V_{\Gamma}^1$, ${\bm v}_\alpha$ is a gauge equivalence class of solutions to \eqref{vortex} for volume form $\sigma_\alpha$ with boundary condition \eqref{boundary} (${\bm v}_\alpha$ is identified via $f_\alpha$ with a gauge equivalence class of affine vortices over ${\mb C}$ or ${\mb H}$). 

\item For each $v_\alpha \in V_{\Gamma}^\infty$, ${\bm v}_\alpha$ is an $I$-holomorphic map $\bar u_\alpha: \Sigma_\alpha \to X$ satisfying the boundary condition $\bar u_\alpha(\partial \Sigma_\alpha) \subset L$.

\end{itemize}
They satisfy the following conditions. 
\begin{enumerate}

\item {\bf (Matching condition)} For each finite edge $e \in E_\Gamma^-$ connecting $v_\beta$ and $v_{\alpha(e)}$, the evaluation of ${\bm v}_\beta$ at infinity \footnote{When $v_\beta \in V_\Gamma^0 \cup V_\Gamma^\infty$, the evaluation at infinity exists by Gromov's removable singularity theorem; when $v_\beta \in V_\Gamma^1$ the evaluation at infinity is defined in \eqref{evaluation}.} is equal to the evaluation of ${\bm v}_{\alpha(e)}$ at $z_e$. 

\item {\bf (Stability condition)} For each unstable vertex $v_\alpha$, the energy of ${\bm v}_\alpha$ is positive. 
\end{enumerate}
\end{defn}

\begin{convention}
We often use the same type of notations (i.e. ${\bm v}_\alpha$) denote different types of components. By abuse of notation, we also use ${\bm v}_\alpha = (u_\alpha, \phi_\alpha, \psi_\alpha)$ to denote an actual gauged map, not its gauge equivalence class.  
\end{convention}

\begin{defn}
Let ${\mc V}$, ${\mc V}'$ be stable affine vortices of type $\Gamma$ where ${\mc V}$ is as described in Definition \ref{stablevortex} and ${\mc V}' = ({\mc C}', ({\bm v}_\alpha')_{v_\alpha \in V_\Gamma})$. An {\it isomorphism} from ${\mc V}$ to ${\mc V}'$ is an isomorphism of marked scaled curves ${\bm \varphi} = (\varphi_\alpha)_{v_\alpha \in V_\Gamma}$ from ${\mc C}$ to ${\mc C}'$ such that
\begin{enumerate}
\item For each $v_\alpha \in V_\Gamma^0$, ${\bm v}_\alpha = {\bm v}_\alpha' \circ \varphi_\alpha$ as $K$-orbits of maps from $\Sigma_\alpha$ to $V$;

\item For each $v_\alpha \in V_\Gamma^1$, ${\bm v}_\alpha = {\bm v}_\alpha' \circ \varphi_\alpha$ as gauge equivalence classes of gauged maps;

\item For each $v_\alpha \in V_\Gamma^\infty$, ${\bm v}_\alpha = {\bm v}_\alpha' \circ \varphi_\alpha$ as maps from $\Sigma_\alpha$ to $X$. 
\end{enumerate}
\end{defn}

Now we define the notion of sequential convergence of affine vortices towards stable affine vortices. For affine vortices over ${\mb C}$, this is defined by Ziltener \cite{Ziltener_book, Ziltener_thesis}. So we only describe the case for affine vortices over ${\mb H}$. We only state the definition for the case when elements of the sequence have smooth domains but it is standard to extend the following definition to the most general situation.

\begin{defn}\label{convergence}
Let 
\beqn
{\mc V}_k = ( {\mc C}_k, {\bm v}_k) = ( (\Sigma_k, \sigma_k, f_k), (z_{k, e}), (u_k, \phi_k, \psi_k))
\eeqn
be a sequence of stable affine vortices with smooth domains $\Sigma_k \cong {\mb H}$, $l$ boundary marked points and $l^+$ interior marked points, and let 
\beqn
{\mc V} = ({\mc C}, ({\bm v}_\alpha)_{v_\alpha \in V_\Gamma})
\eeqn
be a stable affine vortex of type $\Gamma\in {\rm Tree}(l, l^+)$ as described in Definition \ref{stablevortex}. We say that ${\mc V}_k$ converges (modulo gauge transformation and translation) to ${\mc V}$ if there exist a collection of sequences of M\"obius transformations 
\beqn
\varphi_{k,\alpha}: ({\mb C}, \infty) \to ({\mb C}, \infty)
\eeqn
satisfying the following conditions.

\begin{enumerate}

\item For each $v_\alpha \in V_{\uds\Gamma}$, $\varphi_{k,\alpha}$ is real, namely $\varphi_{k,\alpha}$ maps ${\mb H}$ to ${\mb H}$. 

\item For each $v_\alpha \in V_\Gamma^1$, $\varphi_{k,\alpha}$ is a translation, namely $\varphi_{k,\alpha}(z) = z + t_{k,\alpha}$ for some $t_{k,\alpha} \in {\mb C}$. 

\item For each semi-infinite edge $e \in E_\Gamma^\infty$, there holds
\beqn
\lim_{k\to \infty} f_{\alpha(e)}^{-1}\varphi_{k,\alpha}^{-1}( f_k(z_{k,e})) = z_e \in \Sigma_{\alpha(e)}.
\eeqn
For each finite edge $e\in E_\Gamma^-$ connecting $v_\alpha$ and $v_\beta$ (where $v_\beta$ is closer to the root), there holds
\beqn
f_\beta^{-1} \circ \varphi_{k, \beta}^{-1} \circ \varphi_{k, \alpha} \circ f_\alpha \xrightarrow[k\to \infty]{C^\infty_{\rm loc}(\Sigma_\alpha)} z_e \in \Sigma_\beta. \footnote{When $\Gamma$ is stable, the first three items of this definition constitute the definition of the notion of sequential convergence of the domain moduli.}
\eeqn

\item For each $v_\alpha \in V_\Gamma^0 \cup V_\Gamma^1$, the sequence ${\bm v}_k \circ (f_k^{-1} \circ \varphi_{k,\alpha} \circ f_\alpha)$ converges modulo gauge transformations to ${\bm v}_\alpha$.

\item For each $v_\alpha \in V_\Gamma^\infty$, there holds 
\beqn
\mu(u_k \circ (f_k^{-1} \circ \varphi_{k,\alpha} \circ f_\alpha))  \xrightarrow[k\to \infty]{C^0_{\rm loc}(\Sigma_\alpha)} 0.
\eeqn
Moreover, there is a smooth map $u_\alpha: \Sigma_\alpha \to \mu^{-1}(0)$ whose projection to $X$ is the holomorphic map $\bar u_\alpha: \Sigma_\alpha \to X$ such that 
\beqn
u_k \circ \varphi_{k,\alpha} \xrightarrow[k\to \infty]{C^0_{\rm loc}(\Sigma_\alpha)} u_\alpha
\eeqn

\item There is no energy lost, namely, 
\beqn
\lim_{k \to \infty} E({\bm v}_k) = E({\mc C}):= \sum_{v_\alpha \in V_\Gamma} E({\bm v}_\alpha).
\eeqn

\end{enumerate}
\end{defn}

\subsection{Families of almost complex structures}

We would like to allow almost complex structures to be domain-dependent and depend on the domain moduli. This is the perturbation schemed used in \cite{Woodward_Xu}. We want the almost complex structures to be close to a {\it base almost complex structure} $J_V$. It is a $K$-invariant, $\omega_V$-compatible almost complex structure on $V$. We also require that $L_V$ is totally real with respect to $J_V$. Then $\omega_V$ and $J_V$ determine a Riemannian metric 
\beqn
g_V(\cdot, \cdot):= \omega_V( \cdot, J_V \cdot).
\eeqn

We would like the domain dependent almost complex structure to respect certain splittings of the tangent bundle. Let 
\beqn
(TV)_K \subset TV
\eeqn
be the set of infinitesimal $K$-actions. Define 
\beqn
(TV)_G:= (TV)_K \oplus J_V (TV)_K \subset TV
\eeqn
which is formally the set of ``infinitesimal $G$-actions. By Hypothesis \ref{hyp21}, $(TV)_G$ is a subbundle of $TV$ near $\mu^{-1}(0)$. It is easy to to see that the orthogonal complement 
\beqn
((TV)_G)_{g_V}^\bot \subset TV
\eeqn
with respect to the metric $g_V$ is tangent to $\mu^{-1}(0)$. Moreover, if $\pi_X: \mu^{-1}(0) \to X$ is the projection, then there is a $K$-equivariant isomorphism 
\beq\label{eqn26}
((TV)_G)_{g_V}^\bot|_{\mu^{-1}(0)} \cong \pi_X^* TX. 
\eeq

Now we describe the type of domain-dependent almost complex structures we will use. Fix a pair of integers $(l, l^+)$. There is a {\it universal curve}
\beqn
\ov{\mc U}_{l, l^+} \to \ov{\mc M}_{l, l^+}
\eeqn
whose total space is $\ov{\mc U}_{l, l^+} \cong \ov{\mc M}_{l, l^+ + 1}$ and the projection $\ov{\mc U}_{l,l^+} \to \ov{\mc M}_{l, l^+}$ is forgetting the last interior marking. The fibre of the universal curve at a point $p\in \ov{\mc M}_{l, l^+}$ can be canonically identified with a stable marked scaled curve (up to isomorphism) representing the point $p$. The universal curve has the structure of smooth manifold with certain polyhedral corners (see the main result of \cite{Mau_Woodward_2010} for the case $l= 0$ or $l^+ = 0$) which allows one to have a notion of smoothness. There is a closed subset $\ov{\mc U}_{l, l^+}^{\rm node} \subset \ov{\mc U}_{l, l^+}$ corresponding to nodal points. There is also a closed subset $\ov{\mc U}_{l, l^+}^{\rm boundary} \subset \ov{\mc U}_{l, l^+}$ corresponding to boundary points. More generally, for each stable scaled tree $\Gamma \in {\rm Tree}^s (l, l^+)$, denote by
\beqn
{\mc U}_\Gamma \subset \ov{\mc U}_{l, l^+}
\eeqn
the preimage of ${\mc M}_\Gamma \subset \ov{\mc M}_{l, l^+}$ under the projection $\ov{\mc U}_{l, l^+} \to \ov{\mc M}_{l, l^+}$. There is a subset 
\beqn
{\mc U}_\Gamma^\infty \subset {\mc U}_\Gamma
\eeqn
corresponding to points on vertices in $V_\Gamma^\infty$. 

\begin{defn}\label{defn216}
A family of domain-independent almost complex structures is a smooth map $J: \ov{\mc U}_{l, l^+} \to {\mc J}_{\rm tame}(V)$ where ${\mc J}_{\rm tame}(V)$ is the space of $K$-invariant $\omega_V$-tamed almost complex structures on $V$ satisfying the following conditions. 
\begin{enumerate}

\item $J = J_V$ in a neighborhood of $\ov{\mc U}{}_{l, l^+}^{\rm node} \cup \ov{\mc U}{}_{l, l^+}^{\rm boundary}$. 

\item For each stable tree $\Gamma \in {\rm Tree}^s(l, l^+)$, the restriction of $J$ to ${\mc U}{}_\Gamma^\infty \times \mu^{-1}(0)$ preserves the subbundle $((TV)_G)_{g_V}^\bot|_{\mu^{-1}(0)}$. 
\end{enumerate}

\end{defn}

From now on we fix a family of domain-dependent almost complex structures. For each (not necessarily stable) scaled tree $\Gamma\in {\rm Tree}(l, l^+)$ and a marked scaled curve ${\mc C}$ of type $\Gamma$, the family induces a domain-dependent almost complex structure $J_\alpha: \Sigma_\alpha \to {\mc J}_{\rm tame}(V)$ on each component $\Sigma_\alpha \subset {\mc C}$; if $v_\alpha \in V_\Gamma^\infty$ then via the isomorphism \eqref{eqn26} $J_\alpha$ induces a map
\beqn
I_\alpha: \Sigma_\alpha \to {\mc J}_{\rm tame}(X).
\eeqn
Denote by 
\beqn
{\mc M}_\Gamma(V, L_V)
\eeqn
the moduli space of isomorphism classes of stable affine instantons of type $\Gamma$ with respect to the family of almost complex structures. Moreover, define 
\beqn
\ov{\mc M}_{l, l^+}(V, L_V):= \bigsqcup_{\Gamma} {\mc M}_\Gamma(V, L_V)
\eeqn
where the union is taken over all (not necessarily stable) scaled trees $\Gamma\in {\rm Tree}(l, l^+)$. One can define a notion of sequential convergence by generalizing Definition \ref{convergence}. The notion of sequential convergence then induces a topology on $\ov{\mc M}_{l, l^+}(V, L_V)$. Moreover, for each real number $E> 0$, let 
\beqn
\ov{\mc M}_{l, l^+}^{\leq E}(V, L_V)\subset \ov{\mc M}_{l, l^+}(V, L_V)
\eeqn
be the subset of stable affine vortices whose energy (computed using a fixed Riemannian metric on $V$) is at most $E$. 

We recall the compactness theorem for the moduli space of stable affine vortices. The case of affine vortices over ${\mb C}$ was proved by Ziltener \cite{Ziltener_book, Ziltener_thesis}. The case of affine vortices over ${\mb H}$, though not explicitly stated and proved, can be produced from \cite{Wang_Xu}.\footnote{In both \cite{Ziltener_book, Ziltener_thesis} and \cite{Wang_Xu} the compactness theorem is proved for a fixed domain-independent almost complex structure. However the proof extends without difficulty to the current case.} 

\begin{thm}[Compactness theorem]\label{compactness} Fix $l, l^+\geq 0$ with $l+ l^+ \geq 1$ and a family of domain-dependent almost complex structures $J: \ov{\mc U}_{l, l^+} \to {\mc J}_{\rm tame}(V)$. For any $E>0$, the subset $\ov{\mc M}_{l, l^+}^{\leq E} (V, L_V)$ is sequentially compact. 
\end{thm}

\subsection{Codimension one degeneration}\label{subsection25}

The goal of this paper is to provide a local model for the compactified moduli space near a point in a codimension one stratum. In this subsection, we first describe the local model of the domain moduli for such a stratum.

\begin{defn}\label{simple}
A scaled tree $\Gamma$ is called {\it simple} if $V_\Gamma^\infty = \{v_0\}$ and $V_\Gamma^1 = V_\Gamma \setminus \{v_0\}$. 
\end{defn}
\noindent From now on we fix a simple stable scaled tree $\Gamma$ with $l$ boundary semi-infinite edges, $l^+$ interior semi-infinite edges, with 
\begin{align*}
&\ \# V_{\uds\Gamma}^1 = m,\ &\ \# V_{\Gamma^{\rm sup}}^1 = m^+.
\end{align*}
From now on throughout this paper, we order and label the finite edges by integers as 
\beqn
e_1, \ldots, e_m;\ e_1^+, \ldots, e_{m^+}^+.
\eeqn
Then the corresponding vertices $v_\alpha \in V_\Gamma^1$ are denoted by $v_1, \ldots, v_m; v_1^+, \ldots, v_{m^+}^+$. We also label them in a different way as
\beqn
v_i,\ i = 1, \ldots, m+m^+,\ {\rm where}\ v_{m + i}:= v_{i}^+.
\eeqn
For each $e \in E_\Gamma^\infty$, let $i(e)$ be the integer labeling the vertex $v_{i(e)}$ containing $e$. 

We would like to see the local structure of $\ov{\mc M}_{l, l^+}$ near a point $p \in {\mc M}_\Gamma$. First, it is easy to see that ${\mc M}_\Gamma$ is a smooth manifold. Indeed, for $v_i \in V_\Gamma^1$, let $\Gamma_i$ be the scaled tree with only one vertex $v_i$ (with scale $1$) whose semi-infinite edges are all semi-infinite edges in $\Gamma$ attached to  $v_1$. Then there is a canonical homeomorphism
\beqn
{\mc M}_\Gamma  \cong {\mc M}_{m + 1, m^+}^{\rm disk} \times \prod_{i=1}^{m+m^+} {\mc M}_{\Gamma_i}.
\eeqn
Here ${\mc M}_{m + 1, m^+}^{\rm disk}$ is the moduli space of marked disks with $m+1$ boundary markings and $m^+$ interior markings. A point $p \in {\mc M}_\Gamma$ can be written as $(p_0, \ldots, p_{m+m^+})$. Denote 
\begin{align*}
&\ W_i:= T_{p_i} {\mc M}_{\Gamma_i},\ i=1, \ldots, m+m^+;\ &\ W_0:= T_{p_0} {\mc M}_{m + 1, m^+}^{\rm disk} \cong {\mb R}^{m + 2m^+ -2}.
\end{align*}
We fix the following data.

\begin{enumerate}

\item Metrics on $W_i$ for $i=0, \ldots, m+m^+$. For $\tau >0$ let $W_i^\tau \subset W_i$ be the radius $\tau$ ball centered at the origin.

\item For $i=1, \ldots, m + m^+$, a configuration of points ${\bm z}_i = (z_e)_{i(e) = i}$ in ${\mb C}$ or ${\mb H}$ representing the point $p_i \in {\mc M}_{\Gamma_i}$.

\item A configuration of points ${\bm z}_0 = (z_1, \ldots, z_m; z_1^+, \ldots, z_{m^+}^+)$ in ${\mb H}$ such that the marked disk $({\mb H} \cup \{\infty\}; z_1, \ldots, z_m, \infty; z_1^+, \ldots, z_{m^+}^+)$ represents the point $p_0 \in {\mc M}_{m + 1, m^+}^{\rm disk}$.

\item A sufficiently small $\tau >0$ and for each $v_i \in V_\Gamma^1$, a family of configurations 
\beqn
{\bm z}_i ({\bm w}_i) = (z_e({\bm w}_i))_{i(e) = i}
\eeqn
parametrized smoothly by ${\bm w}_i \in W_i^\tau$ such that the map 
\beqn
{\bm w}_i \mapsto [ ( {\mb A}_i, ds \wedge dt, {\rm Id}), {\bm z}_i ({\bm w}_i)]
\eeqn
(where $[\cdot]$ means the isomorphism class of marked scaled curves) is a homeomorphism onto a neighborhood of $p_i$ in ${\mc M}_{\Gamma_i}$.

\item A smooth family of diffeomorphisms 
\beq\label{eqn27}
f_{{\bm w}_0} = s_{{\bm w}_0} + {\bm i} t_{{\bm w}_0}: {\mb H} \to {\mb H}
\eeq
such that the support of $df_{{\bm w}_0} -  {\rm Id}$ is contained in a common compact subset which is disjoint from $z_1, \ldots, z_m; z_1^+, \ldots, z_{m^+}^+$, and such that the map 
\beqn
W_0^\tau \ni {\bm w}_0 \mapsto p_0({\bm w}_0) = [ \Sigma_{{\bm w}_0}; z_1, \ldots, z_m, \infty; z_1^+, \ldots, z_{m^+}^+]
\eeqn
is a homeomorphism onto a neighborhood of $p_0$ in ${\mc M}_{m + 1, m^+}^{\rm disk}$. Here $\Sigma_{{\bm w}_0}$ is the bordered Riemann surface whose underlying surface is ${\mb D} \cong {\mb H} \cup \{\infty\}$ and whose complex structure is induced from the coordinate $s_{{\bm w}_0} + {\bm i} t_{{\bm w}_0}$. Equivalently, The diffeomorphism $z_{{\bm w}_0}$ gives a family of configurations
\beqn
{\bm z}_0({\bm w}_0):= (z_i({\bm w}_0))_{i=1}^{m+m^+} = (z_{{\bm w}_0}(z_i))_{i=1}^{m+m^+}
\eeqn
so that the marked disk (with the standard complex structure)
\beqn
({\mb H} \cup \{\infty\}; z_1({\bm w}_0), \ldots, z_m({\bm w}_0), \infty; z_1^+({\bm w}_0), \ldots, z_{m^+}^+({\bm w}_0))
\eeqn
represents the same point $p_0({\bm w}_0) \in {\mc M}_{m+1, m^+}^{\rm disk}$.
\end{enumerate}

Denote
\begin{align*}
&\ W_{\rm def}:= \bigoplus_{v_\alpha \in V_\Gamma} W_\alpha,\ &\ W_{\rm def}^\tau:= \prod_{v_\alpha \in V_\Gamma} W_\alpha^\tau.
\end{align*}
Elements of $W_{\rm def}$, called {\it deformation parameters}, are denoted by ${\bm w} = ({\bm w}_\alpha)_{v_\alpha \in V_\Gamma}$. The data we just fixed induce a family of marked stable scaled curves of type $\Gamma$, denoted by 
\beqn
{\mc C}_{0, {\bm w}},\ {\bm w} \in W_{\rm def}^\tau.
\eeqn

We would like to turn on the gluing parameter and obtain a local model of the moduli $\ov{\mc M}_{l, l^+}$. Denote 
\beqn
z_{i,\epsilon}({\bm w}_0):= \frac{z_i({\bm w}_0)}{\epsilon},\ i = 1, \ldots, m + m^+
\eeqn
and for each $e \in E_\Gamma^\infty$, define
\beqn
z_{e,\epsilon}({\bm w}):= z_e({\bm w}_{i(e)}) + z_{i(e), \epsilon}({\bm w}_0).
\eeqn
This gives a family of configurations of points in ${\mb H}$. For $\tau>0$ small enough, define
\beq\label{domainglue}
\Phi_p: [0, \tau)\times W_{\rm def}^\tau \to \ov{\mc M}_{l, l^+},  (\epsilon, {\bm w}) \mapsto \left\{ \begin{array}{cc} {[} {\mc C}_{0, {\bm w}} {]},\ &\ \epsilon = 0,\\
                                                         {[} {\mc C}_{\epsilon, {\bm w}} {]}:={[} ({\mb H}, ds\wedge dt, {\rm Id}); (z_{e, \epsilon}({\bm w}))_{e \in E_\Gamma^\infty} {]},\ &\ \epsilon \neq 0. \end{array} \right.
\eeq
Then one has
\begin{lemma}
For $\tau$ small enough, the map $\Phi_p$ is a homeomorphism onto an open neighborhood of $p$ in $\ov{\mc M}_{l, l^+}$. 
\end{lemma}

Using the domain local model given above the family of domain-dependent almost complex structures can be regarded as depending on the gluing parameter and the deformation parameter. More precisely, for any stable marked scaled curve ${\mc C}$ representing a point in the image of $\Phi_p$, we can identify ${\mc C}$ with a certain ${\mc C}_{\epsilon, {\bm w}}$. Then $J: \ov{\mc U}_{l, l^+} \to {\mc J}_{\rm tame}(V)$ induces a domain-dependent almost complex structure
\beqn
J_{\epsilon, {\bm w}}: {\mc C} \cong {\mc C}_{\epsilon, {\bm w}} \to {\mc J}_{\rm tame}(V).
\eeqn

For the purpose of constructing the gluing map, we slightly modify the local model given above. The idea is to fix the positions of the nodal points while varying the complex structure on the disk component away from the nodal points. Denote 
\beq\label{epsilonrescaling}
s_\epsilon: {\mb H} \to {\mb H},\ z\mapsto \epsilon z
\eeq
and call it the {\it $\epsilon$-rescaling}. Then let $\Sigma_{\epsilon, {\bm w}}$ be bordered Riemann surface whose underlying surface is ${\mb H}$ and whose complex structure is induced from the coordinate
\beqn
f_{\epsilon, {\bm w}} = s_{\epsilon, {\bm w}} + {\bm i} t_{\epsilon, {\bm w}}:= s_\epsilon^{-1} \circ f_{{\bm w}_0} \circ s_\epsilon.
\eeqn
Here $f_{{\bm w}_0}: {\mb H} \to {\mb H}$ is the family of diffeomorphisms \eqref{eqn27}. Furthermore, define 
\beqn
\sigma_{\epsilon, {\bm w}}':= ds_{\epsilon, {\bm w}} \wedge dt_{\epsilon, {\bm w}}
\eeqn
which is a volume form on $\Sigma_{\epsilon, {\bm w}}$. Define
\beqn
{\mc C}_{\epsilon, {\bm w}}': =\big( (\Sigma_{\epsilon, {\bm w}},  \sigma_{\epsilon, {\bm w}}', f_{\epsilon, {\bm w}}) ; {\bm z}_{\epsilon, {\bm w}}' \big)
\eeqn
where ${\bm z}_{\epsilon, {\bm w}}'$ is the configuration 
\beqn
\Big(z_{e, {\bm w}}' = z_e({\bm w}_{i(e)}) + \frac{z_i}{\epsilon}\Big)_{e\in E_\Gamma^\infty}.
\eeqn
Then ${\mc C}_{\epsilon, {\bm w}}'$ also parametrizes a neighborhood of $p$ in $\ov{\mc M}_{l, l^+}$ and gives the same map as $\Phi_p$ (see \ref{domainglue}) and the family of biholomorphisms
\beqn
f_{\epsilon, {\bm w}}: {\mc C}_{\epsilon, {\bm w}}' \cong {\mc C}_{\epsilon, {\bm w}}
\eeqn
intertwines the volume form $\sigma_{\epsilon, {\bm w}}'$ resp. marked points ${\bm z}_{\epsilon, {\bm w}}'$ on the former with the standard volume form $ds \wedge dt$ resp. marked points ${\bm z}_{\epsilon, {\bm w}}$ on the latter. 

Using this new family the vortex equation \eqref{affine} takes a slightly different form. The family of almost complex structures $J$ induces a domain-dependent almost complex structure $J_{\epsilon, {\bm w}}: \Sigma_{\epsilon, {\bm w}} \to {\mc J}_{\rm tame}(V)$. Then \eqref{affine} is equivalent to the following 
\beq\label{eqn210}
\begin{split}
\partial_{s_{\epsilon, {\bm w}}} u + {\mc X}_\phi(u) + J_{\epsilon, {\bm w}}( \partial_{t_{\epsilon, {\bm w}}} u + {\mc X}_\psi(u)) = & 0,\\
\partial_{s_{\epsilon, {\bm w}}} \psi - \partial_{t_{\epsilon, {\bm w}}} \phi + [\phi, \psi] + \mu(u) = & 0.
\end{split}
\eeq

\section{The Main Theorem}\label{section3}

In this section we state the main theorem of this paper, under a precise version of the transversality assumption.

\subsection{Local model of affine vortices}

We recall the main result of \cite{Venugopalan_Xu} which provides the local model of affine vortices.

Since the domain of affine vortices is noncompact, Fredholm theory requires certain Sobolev completions of the domains and codomains of the linearized operator. The particular choices of these norms reflects the complicated asymptotic behavior of affine vortices. Choose a smooth weight function $\rho_{\mb A}: {\mb A} \to [1, +\infty)$ such that outside the unit disk, $\rho_{\mb A} (z) = |z|^2$. Let $\delta \in {\mb R}$ be a real parameter. For $f \in L{}^p_{\rm loc}({\mb A})$ define
\beqn
\| f\|_{L^{p, \delta}}:= \left[ \int_{\mb A} |f(z)|^p (\rho_{\mb A} (z))^{\frac{p\delta}{2}} ds dt \right]^{\frac{1}{p}}.
\eeqn
We can generalize this norm to sections of certain Euclidean bundle $E$ over certain subset of ${\mb A}$. On the other hand, we will use to different norms to complete the space of functions with local $W^{1,p}$-regularity. Define 
\beq\label{hnorm}
\| f \|_{L^{1,p,\delta}_H}:= \| f \|_{L^\infty} + \| \nabla f \|_{L^{p,\delta}}
\eeq
and 
\beq\label{gnorm}
\| f \|_{L^{1,p,\delta}_G}:= \| f \|_{L^{p,\delta}} + \| \nabla f \|_{L^{p,\delta}}.
\eeq

The completions of the domains and codomains of the linearized operators are modelled on the above Sobolev spaces. We first introduce a few new notations. Let ${\bm v} = (u, \phi, \psi)$ be an affine vortex. Denote the connection form by 
\beqn
a = \phi ds + \psi dt \in \Omega^1( {\mb A}, {\mf k}). 
\eeqn
Introduce a covariant derivative of sections $\xi$ of $u^* TV$ resp. ${\mf k}$-valued functions $f$ on ${\mb A}$ by
\beqn
\begin{split}
\nabla_s^a \xi = \nabla_s \xi + \nabla_\xi {\mc X}_\phi\ & {\rm resp.}\ \nabla_s^a f = \partial_s f + [\phi, f],\\
\nabla_t^a \xi = \nabla_t \xi + \nabla_\xi {\mc X}_\psi\ & {\rm resp.}\ \nabla_t^a f = \partial_t f + [\psi, f].
\end{split}
\eeqn
Then for ${\bm \xi} = (\xi, \eta, \zeta) \in \Gammait( u^* TV \oplus {\mf k} \oplus {\mf k})$, define the norm 
\beq\label{eqn33}
\| {\bm \xi} \|_{L^{1, p, \delta}_m}:= \| \xi \|_{L^\infty} + \| \nabla^{a} {\bm \xi} \|_{L^{p, \delta}} + \| d\mu(u)(\xi) \|_{L^{p, \delta}} + \| d \mu(u)(J_V \xi) \|_{L^{p, \delta}} + \| \eta \|_{L^{p, \delta}} + \| \zeta \|_{L^{p,\delta}}.
\eeq
Here $m$ stands for ``mixed.'' It was proved in \cite{Ziltener_book} for ${\mb A} = {\mb C}$ and in \cite{Venugopalan_Xu} for ${\mb A} = {\mb H}$ that the above norm is complete. We define 
\beqn
\hat{\mc B}_{\bm v}:=  \Big\{ {\bm \xi} = (\xi, \eta, \zeta) \in W^{1,p}_{\rm loc}( {\mb A}, u^* TV \oplus {\mf k} \oplus {\mf k} )\ |\ \| {\bm \xi} \|_{L^{1,p,\delta}_m} < \infty,\ \xi|_{\partial {\mb H}} \subset TL_V,\  \zeta|_{\partial {\mb H}} = 0 \Big\}.
\eeqn
The ``Coulomb slice'' is the subspace
\beq\label{eqn34}
{\mc B}_{\bm v}:= \Big\{ {\bm \xi} = (\xi, \eta, \zeta) \in \hat{\mc B}_{\bm v}\ |\ \partial_s \eta + [\phi, \eta] + \partial_t \zeta + [\psi, \zeta] + d\mu(u)(J_V \xi) = 0  \Big\}.
\eeq

The main results of \cite{Venugopalan_Xu} can be summarized as follows. 

\begin{thm}[Local model] \label{localmodel}
Fix $p > 2$ and $\delta \in (1 -\frac{2}{p}, 1)$. Let ${\mb A}$ be either ${\mb C}$ or ${\mb H}$. Let $(J_z)$ be a smooth family of $K$-invariant almost complex structures on $V$ parametrized by $z \in {\mb A}$ such that $J_z = J_V$ for $|z|$ large. There is a Banach manifold ${\mc B}$, a Banach vector bundle ${\mc E} \to {\mc B}$, and a Fredholm section $F: {\mc B} \to {\mc E}$ satisfying the following conditions. 
\begin{enumerate}
\item Each element of ${\mc B}$ is a gauge equivalence class of gauged maps from ${\mb A}$ to $V$ of regularity $W^{1,p}_{\rm loc}$ having limits at infinity as $K$-orbits. Moreover, the evaluation map at infinity $\ev_\infty:  {\mc B} \to X$ is a smooth map.

\item The zero locus of $F$ consists of gauge equivalence classes of affine vortices over ${\mb A}$ and the natural map from $F^{-1}(0)$ to the moduli space of gauge equivalence classes of solutions to the vortex equation over ${\mb A}$ (for the almost complex structure $J_z$) is a homeomorphism between the Banach manifold topology on the former and the $C_{\rm loc}^\infty$ topology on the latter. 

\item For each point ${\bm p} \in F^{-1}(0)$, choose a smooth representative ${\bm v} = (u, \phi, \psi)$, the tangent space of ${\mc B}$ at ${\bm p}$ is isomorphic to the space ${\mc B}_{\bm v}$ (see \eqref{eqn34}), the fibre of ${\mc E}$ at ${\bm p}$ is isomorphic to ${\mc E}_{\bm v} \cong L^{p,\delta}({\mb A}, u^* TV \oplus {\mf k})$, and the linearization of $F$ at ${\bm p}$ is 
\beq
D_{\bm p} F({\bm \xi}) = \left[ \begin{array}{c} \nabla_s^a \xi + J_z \nabla_t^a \xi + (\nabla_\xi J_z) ( \partial_t u + {\mc X}_\psi) + {\mc X}_\eta + J_z {\mc X}_\zeta\\
                                             \nabla_s^a \zeta - \nabla_t^a \eta + d\mu (\xi) \end{array} \right].
\eeq

\item The linearized operator $D_{\bm p} F: {\mc B}_{\bm v} \to {\mc E}_{\bm v}$ is a Fredholm operator.\footnote{The case when ${\mb A} = {\mb C}$ is proved by Ziltener, see \cite[(1.27)]{Ziltener_book}.}
\end{enumerate}
\end{thm}

\subsection{Regular configurations}

The existence of the gluing map relies on a certain transversality assumption. In this subsection we specify the transversality conditions. Let $\Gamma$ be a stable simple domain type (see Definition \ref{simple}) with $l$ boundary semi-infinite edges and $l^+$ interior semi-infinite edges. Let $J: \ov{\mc M}_{l, l^+} \to {\mc J}_{\rm tame}(V)$ be a smooth family of domain-dependent almost complex structures. Let ${\mc V} = ({\mc C}, ({\bm v}_\alpha)_{v_\alpha \in V_\Gamma})$ be a stable affine vortex of type $\Gamma$ representing a point 
\beqn
{\bm p} = [ {\mc V}] \in {\mc M}_\Gamma(V, L_V) \subset \ov{\mc M}_{l, l^+}(V,L_V).
\eeqn

We order vertices underlying affine vortices components in the same way as we label vertices of the tree in Subsection \ref{subsection25}. Namely, let the components which are affine vortices over ${\mb H}$ be ${\bm v}_1, \ldots, {\bm v}_m$ and the components which are affine vortices over ${\mb C}$ be ${\bm v}_1^+, \ldots, {\bm v}_{m^+}^+$. We also introduce
\beqn
{\bm v}_{m+j}:= {\bm v}_j^+,\ j=1, \ldots, m^+
\eeqn
so ${\bm v}_1, \ldots, {\bm v}_{m+m^+}$ are all the affine vortex components of ${\mc V}$. Up to isomorphism, we can assume that the domain of ${\bm v}_i$ be ${\mb A}_i$ which is either ${\mb H}$ or ${\mb C}$ with the standard structures. 

For each vertex $v_i \in V_\Gamma^1$, let ${\mc B}_i$ be the Banach manifold parametrizing gauge equivalence classes of gauged maps defined over ${\mb A}_i$ and let ${\mc E}_i \to {\mc B}_i$ be the corresponding Banach vector bundle for a chosen $p \in (2, 4)$ and $\delta = \delta_p = 2- \frac{4}{p}$ (see Theorem \ref{localmodel}). On the other hand, for the holomorphic disk component, there is a Banach manifold
\beqn
{\mc B}_0 = \big\{ \bar u_0' \in W^{1,p}({\mb D}^2, X)\ |\ \bar u_0' (\partial {\mb D}^2) \subset L \big\}.
\eeqn
There is a Banach vector bundle ${\mc E}_0 \to {\mc B}_0$ whose fibre over a map $\bar u_0'$ is
\beqn
{\mc E}_0 |_{\bar u_0'} = L^p({\mb D}^2, \Lambda^{0,1} \otimes (\bar u_0')^* TX).
\eeqn

Now we describe the local model for the singular configuration. Define
\beqn
{\mc B}_\Gamma^\#:= {\mc B}_0 \times \prod_{i=1}^{m+m^+} {\mc B}_i.
\eeqn
This is direct product of spaces of component-wise configurations whose values at nodal points may not match. Define
\beq\label{eqn36}
{\mc B}_\Gamma:= \Big\{ (\bar u_0', {\bm v}_1', \ldots, {\bm v}_{m+m^+}') \in {\mc M}_\Gamma \times {\mc B}_\Gamma^\# \ |\ {\rm the\ ``matching\ condition"} \Big\}.
\eeq
Here the matching condition means the following: 
\begin{itemize}
\item {\bf (Matching condition)} For each vertex $v_i \in V_\Gamma^1$, one has 
\beqn
{\bm v}_i' (\infty) = \bar u_0' ( z_i)
\eeqn
as points in $X$ (if $v_\alpha \in V_{\Gamma^{\rm sup}}$) or $L$ (if $v_\alpha \in V_{\uds\Gamma}$). 
\end{itemize}

There is also the direct sum of Banach vector bundles. Denote
\beqn
{\mc E}_\Gamma^\#:= {\mc E}_0 \oplus \bigoplus_{i=1}^{m + m^+} {\mc E}_i\to {\mc B}_\Gamma^\#.
\eeqn
By pulling back and restriction one has the Banach vector bundle 
\beqn
{\mc E}_\Gamma \to {\mc B}_\Gamma.
\eeqn
Then using the family of domain-dependent almost complex structure $J$ specified previously, using Theorem \ref{localmodel} and the standard facts about local models of holomorphic disks, the vortex equation and the holomorphic curve equation provide a Fredholm section 
\beqn
{\mc F}_\Gamma: W_{\rm def} \times {\mc B}_\Gamma \to {\mc E}_\Gamma.
\eeqn
By Theorem \ref{localmodel}, there is a natural homeomorphism 
\beqn
{\mc F}_\Gamma^{-1}(0) \cong {\mc M}_\Gamma(V, L_V).
\eeqn
The linearization of ${\mc F}_\Gamma$ at point ${\bm p} \in {\mc F}_\Gamma^{-1}(0)$ (which is independent of how one locally trivializes the bundle ${\mc E}_\Gamma$) is a linear Fredholm operator 
\beq\label{linear}
D_{\bm p} {\mc F}_\Gamma : W_{\rm def}\oplus T_{\bm p} {\mc B}_\Gamma \to {\mc E}_\Gamma|_{{\bm p}}.
\eeq

\begin{defn}\label{defn32}
A point ${\bm p} \in {\mc F}_\Gamma^{-1}(0)$ is called {\it regular} if the linear operator $D_{\bm p} {\mc F}_\Gamma$ in \eqref{linear} is surjective; the point ${\bm p}$ is called {\it rigid} if the index of $D_{\bm p} {\mc F}_\Gamma$ is zero.
\end{defn}

Now we can state our main theorem.

\begin{thm}\label{mainthm}
Suppose $\Gamma$ is a simple stable scaled tree with $l$ boundary semi-infinite edges and $l^+$ interior semi-infinite edges. Given a regular and rigid point ${\bm p} \in {\mc M}_\Gamma( V, L_V)$, there is an open neighborhood ${\mc U} \subset \ov{\mc M}_{l, l^+} (V, L_V)$ of ${\bm p}$ which is homeomorphic to an interval $[0, \tau)$ where $0$ corresponds to ${\bm p}$.  
\end{thm}

\begin{rem}
There are two reasons why we do not consider the case of having unstable components. If there are unstable components one just adds extra marked points to stabilize, which does not make an essential difference. On the other hand, in the case of using stabilizing divisors and domain dependent almost complex structures (which is the approach of \cite{Woodward_Xu}), there cannot be unstable components in a stable object. 
\end{rem}

\begin{rem}
The construction of this paper extends without difficulty to the situation when $\Gamma$ is a stable simple scaled tree with empty base. In that case the gluing parameter is a complex number and all equations have no boundary condition.
\end{rem}

\section{Preparation}\label{section4}

The proof of Theorem \ref{mainthm} is straightforward and resembles the gluing construction of holomorphic curves. However, because the singular configuration we start with mixes two different types of equations, and because the analysis of affine vortices is much more complicated than holomorphic curves, to prepare for the straightforward argument, one needs to adjust the standard gluing apparatus in many aspects. In this section we make various technical preparations. In Subsection \ref{subsection41} we specify a Riemannian metric allowing one to construct nearby gauged maps and a connection needed to linearize the vortex equation at a general gauged map. In Subsection \ref{subsection42} we recall the ``augmented linearization'' of the vortex equation which is convenient to use in the case when the equation has gauge symmetry. In Subsection \ref{subsection43} we fix the gauge for affine vortex components such that the affine vortices have good asymptotic behavior. In Subsection \ref{subsection44} we rewrite the linear theory of the holomorphic disk component, putting it into a setting closer to the affine vortex equation. In Subsection \ref{subsection45} and Subsection \ref{subsection46}, we define an $\epsilon$-dependent family of singular configurations, define an $\epsilon$-dependent bounded operator $\hat{\mc Q}_{\Gamma, \epsilon}$ (playing the role as a right inverse of the linearization at the singular configuration), and show that the operator norm of $\hat{\mc Q}_{\Gamma, \epsilon}$ has a uniform upper bound.

\subsection{Metric and connection}\label{subsection41}

We will need a good metric on $V$ and will use its exponential map to identify small infinitesimal deformations with nearby gauged maps. 

\begin{lemma}(cf. \cite[Lemma 49]{Wang_Xu}) \label{metric}
There exists a $K$-invariant Riemannian metric $h_V$ on $V$ satisfying the following conditions. 
\begin{enumerate}

\item The base almost complex structure $J_V$ is $h_V$-isometric along $L_V$.\footnote{We do not require that $J_V$ is isometric everywhere.}

\item $J_V (TL_V)$ is orthogonal to $TL_V$ with respect to $h_V$.

\item $L_V$ and $\mu^{-1}(0)$ are totally geodesic with respect to $h_V$.

\item $T\mu^{-1}(0)$ is orthogonal to the distribution spanned by $J_V {\mc X}_a$ for all $a \in {\mf k}$. 

\item The subbundle $((TV)_G)_{g_V}^\bot|_{\mu^{-1}(0)} \subset TV|_{\mu^{-1}(0)}$ (which is the orthogonal complement of $(TV)_G$ with respect to $g_V$) is also the orthogonal complement of $(TV)_G$ with respect to $h_V$.

\end{enumerate}
\end{lemma}

\begin{proof}
Recall that one has the $K$-equivariant decomposition
\beq\label{eqn41}
T\mu^{-1}(0) = ((TV)_G)_{g_V}^\bot|_{\mu^{-1}(0)} \oplus (TV)_K|_{\mu^{-1}(0)}.
\eeq
Moreover one has a $K$-equivariant isomorphisms
\begin{align}\label{eqn42}
&\ ((TV)_G)_{g_V}^{\bot}|_{\mu^{-1}(0)} \cong \pi_X^* TX,\ &\ (TV)_K|_{\mu^{-1}(0)} \cong X \times {\mf k}.
\end{align}
Here 
\beqn
\pi_X: \mu^{-1}(0) \to X  = \mu^{-1}(0)/K
\eeqn
is the projection. By \cite[Lemma A.3]{Frauenfelder_disk}, there is a Riemannian metric $h_X$ on the symplectic quotient $X$ such that the quotient almost complex structure $J_X$ is isometric, $J_X(TL)$ is orthogonal to $TL$, and $L \subset X$ is totally geodesic. On the other hand, choose a $K$-invariant metric $h_{\mf k}$ on the Lie algebra ${\mf k}$. Then using the decomposition \eqref{eqn41} and the isomorphisms \eqref{eqn42}, one takes the direct sum $\pi_X^* h_X \oplus h_{\mf k}$ which is a Riemannian metric $h_{\mu^{-1}(0)}$ on $\mu^{-1}(0)$. In the proof of \cite[Lemma 49]{Wang_Xu} it was proved that $L_V$ is totally geodesic inside $\mu^{-1}(0)$. 

Next we extend the metric $h_{\mu^{-1}(0)}$ to a metric in $V$. For a small $r >0$, let ${\mf k}_r \subset {\mf k}$ be the radius $r$ ball centered at the origin with respect to the metric $h_{\mf k}$. Define a map 
\beqn
\iota: \mu^{-1}(0) \times {\mf k}_r \to V,\ \ (x, \eta) \mapsto \exp_x (J_V {\mc X}_\eta)
\eeqn
where $\exp$ is the exponential map of the metric $g_V$. Then the map $\iota$ is a $K$-equivariant open embedding ($K$ acts diagonally on $\mu^{-1}(0) \times {\mf k}_r$). Furthermore, the product metric $h_{\mu^{-1}(0)} \times h_{\mf k}$ is $K$-invariant. We define
\beqn
h_V:= \iota_* ( h_{\mu^{-1}(0)} \times h_{\mf k})
\eeqn
over a neighborhood of $\mu^{-1}(0)$ and extend it to a $K$-invariant metric on $V$. The required properties (a)---(e) are straightforward to verify. 
\end{proof}

From now on we fix a metric $h_V$ satisfying properties listed in Lemma \ref{metric}. We use this good metric to decompose the tangent bundle. Choose a $K$-invariant neighborhood $U \subset V$ of $\mu^{-1}(0)$ such that the $K$-action on $U$ is free. Then the set $(TV)_G$ is a subbundle of $TV|_U$. Let $(TV)_H \subset TV|_U$ be the orthogonal complement of $(TV)_G$ with respect to $h_V$. Then one has the orthogonal decomposition
\beq\label{splitting}
TV|_U = (TV)_G \oplus (TV)_H. 
\eeq

\begin{lemma}\label{connection}
There exists a $K$-invariant connection $\nabla$ on $TV$ which preserves the orthogonal decomposition \eqref{splitting} and whose restriction to $L_V$ preserves $TL_V$. 
\end{lemma}

\begin{proof}
Let $\tilde\nabla$ be the Levi--Civita connection of $h_V$. With respect to the splitting \eqref{splitting}, one can write $\tilde \nabla$ in the block matrix form as 
\beqn
\tilde \nabla = \left[ \begin{array}{cc}  P_H \circ \tilde \nabla \circ P_H &  E \\
                                           E^*           &  P_G \circ \tilde \nabla \circ P_G    \end{array} \right]
\eeqn
were $P_H$ resp. $P_G$ is the orthogonal projections onto $(TV)_H$ resp. $(TV)_G$. Then one can construct a $K$-invariant connection $\nabla$ on $TV$ which preserves the metric $h_V$ and which coincides with the diagonal part of $\tilde \nabla$ in the neighborhood $U$ of $\mu^{-1}(0)$. We claim that $\nabla$ satisfies the required property, namely,
\beqn
\nabla_{V_1} V_2 \in \Gammait (TL_V),\ \forall V_1, V_2 \in TL_V.
\eeqn
Indeed, this follows from the fact that $\tilde \nabla_{V_1} V_2 \in TL_V$ and the two orthogonal projections $P_H$ and $P_G$ preserve the subbundle $TL_V$.
\end{proof}

Via the isomorphism $\pi_X^* TX \cong (TV)_H|_{\mu^{-1}(0)}$, the metric $h_V$ induces a Riemannian metric $h_X$ on the symplectic quotient $X$ and the connection $\nabla$ induces a connection $\bar \nabla$ on $TX$ which preserves $h_X$. 

\subsection{The augmented linearization}\label{subsection42}

It is often convenient to use an augmented linearized operator which incorporate the gauge fixing condition instead of the original linearized operator. We first introduce a few new notations. Let ${\bm v} = (u, \phi, \psi)$ be a gauged map from ${\mb A}$ to $V$ (not necessarily an affine vortex), using the preferred connection $\nabla$ chosen from Lemma \ref{connection} one can write down the linearization at ${\bm v}$ of the affine vortex equation with respect to a domain-dependent almost complex structure $J$. The {\it augmented linearization} at ${\bm v}$ is the operator 
\beqn
\hat D_{\bm v}: \hat{\mc B}_{\bm v} \to \hat {\mc E}_{\bm v}:= L^{p,\delta}({\mb A}, u^* TV \oplus {\mf k} \oplus {\mf k})
\eeqn
which reads
\beq\label{augmented}
\hat D_{\bm v}(\xi, \eta, \zeta) = \left[ \begin{array}{c} \nabla_s^a \xi + (\nabla_\xi J) (\partial_t u + {\mc X}_\psi) + J \nabla_t^a \xi + {\mc X}_\eta + J {\mc X}_\zeta \\
                                                           \nabla_s^a \eta + \nabla_t^a \zeta + d\mu(u) (J_V\xi)\\
                                                           \nabla_s^a \zeta - \nabla_t^a \eta + d\mu(u) (\xi) \end{array} \right].
\eeq
(As a convention one put the gauge-fixing condition in the second row.)

\subsection{Representatives of affine vortices}\label{subsection43}

To perform the gluing construction we would like the affine vortices to have good asymptotic behaviors. Let 
\beqn
B_R:= B_R(0) \subset {\mb A}
\eeqn
be the radius $R$ (half) disk centered at the origin of ${\mb A}$.

\begin{lemma}\label{lemma43}
Let $p \in (2,4)$ and $\delta \in (\delta_p, 1) = (2- \frac{4}{p}, 1)$. Let ${\bm v} = (u, a)$ be an affine vortex over ${\mb A}$. Then there exist $R>0$, a gauge transformation $g: {\mb A} \setminus B_R \to K$, an element $\lambda \in {\mf k}$, a point $x \in \mu^{-1}(0)$ satisfying the following conditions. We write 
\beqn
g\cdot {\bm v} = ( \check u, \check a). 
\eeqn
\begin{enumerate}
\item The map $\check u: {\mb A} \setminus B_R \to V$ converges to $x$ at the infinity. Hence if $R$ is sufficiently large, one can write $\check u(z) = \exp_x \check \vartheta (z)$ for $\check \vartheta \in \Gammait( {\mb A} \setminus B_R, T_x V)$.

\item There holds
\beq\label{eqn45}
\| \check \vartheta{}^G \|_{L_G^{1,p, \delta} ({\mb A} \setminus B_R) } + \| \check \vartheta{}^H \|_{L_H^{1,p, \delta}({\mb A} \setminus B_R) }  + \| \check a \|_{L_G^{1,p, \delta}({\mb A} \setminus B_R)}  < \infty.
\eeq
Here $\check \vartheta{}^G$ and $\check \vartheta{}^H$ are the components of $\check \vartheta$ with respect to the decomposition \eqref{splitting} and $L_G^{1,p, \delta}$ and $L_H^{1,p, \delta}$ are the weighted Sobolev norms \eqref{hnorm} and \eqref{gnorm}. 
\end{enumerate}
\end{lemma}

\begin{proof}
By \cite[Lemma 6.1]{Venugopalan_Xu} (for the case ${\mb A}= {\mb C}$) and the straightforward extension to the case ${\mb A} = {\mb H}$ (see \cite[Appendix]{Venugopalan_Xu}), for any $\delta \in (1- \frac{2}{p}, 1) = ( \frac{\delta_p}{2}, 1)$, one can find $R$, $g$, $\lambda$, and $x$ satisfying the first item of this lemma and 
\beq\label{eqn46}
\| \check a \|_{W^{1,p,\delta}({\mb A} \setminus B_R)} + \|  \nabla \check \vartheta \|_{W^{1,p,\delta}({\mb A} \setminus B_R)} < \infty.
\eeq
Here $\nabla$ is the standard derivative of $T_x V$-valued functions. We would like to apply a further gauge transformation. Indeed, there exists a unique $s: {\mb A} \setminus B_R \to {\mf k}$ satisfying
\beq\label{eqn47}
e^{s(z)}  \check\vartheta (z)  \in  ((TV)_K)^\bot.
\eeq
By \eqref{eqn46} there holds
\beq\label{eqn48}
\| \nabla s \|_{W^{1, p,\delta}({\mb A} \setminus B_R)}< \infty.
\eeq
If we replace $g(z)$ by $e^{s(z)} g(z)$, then by \eqref{eqn46} and \eqref{eqn48} there still holds
\beq\label{eqn49}
\| \check a \|_{W^{1,p,\delta}({\mb A} \setminus B_R)} + \| \nabla \check \vartheta \|_{L^{p,\delta}({\mb A} \setminus B_R)}  +  \| \check \vartheta \|_{L^\infty({\mb A} \setminus B_R)} < \infty.
\eeq
Moreover, \eqref{eqn47} implies that for certain $C, C'>0$ (depending only on the metric of $V$)
\beq\label{eqn410}
\| \check \vartheta{}^G \|_{L^{p,\delta}({\mb A}\setminus B_R)} \leq C \| d\mu(x) \cdot \check \vartheta \|_{L^{p,\delta}({\mb A} \setminus B_R)} \leq C' \| \mu( \check{u})  \|_{L^{p,\delta}({\mb A} \setminus B_R)} < \infty.
\eeq
Here the last inequality follows from energy decay property of affine vortices (see \cite[Theorem 1.3]{Ziltener_Decay} and \cite[Proposition A.4]{Venugopalan_Xu}). Then \eqref{eqn45} follows from  \eqref{eqn49} and \eqref{eqn410}. 
\end{proof}

The above lemma allows us to fix a good gauge for the affine vortex components. We fix $\delta_0$ satisfying 
\beq\label{delta}
\max \left\{ 2- \frac{4}{p}, 1- \frac{1}{p} \right\} < \delta_0 < 1.
\eeq
In the context of the main theorem (Theorem \ref{mainthm}), for each of the affine vortex component of the point $ {\bm p} \in {\mc M}_\Gamma(V, L_V)$ labelled by a vertex $v_i \in V_\Gamma^1$, we fix a representative ${\bm v}_i = (u_i, \phi_i, \psi_i)$ which is a smooth affine vortex satisfying properties listed in Lemma \ref{lemma43} for $\delta = \delta_0$. In this gauge the convergence of $\check u_i (z)$ to the limit 
\beq\label{eqn412}
x_i:= \check u_i (\infty)
\eeq
is exponentially fast in the cylindrical coordinates. Indeed, there is a constant $c(p,\delta)>0$ such that
\beq\label{hardy}
\| f - f(\infty) \|_{L^{p,\delta-1}} \leq c(p,\delta) \| \nabla f \|_{L^{p,\delta}},\ \forall f \in L_H^{1,p,\delta}({\mb A}). \footnote{This estimate is called a Hardy-type inequality. See \cite[Proposition 91]{Ziltener_book}.}
\eeq
Therefore
\beqn
\| \check \vartheta_i \|_{L^{p,\delta_0 -1}} \leq c(p,\delta_0) \| \nabla \check \vartheta_i \|_{L^{p,\delta_0}}.
\eeqn
If we identify ${\mb A}_i \setminus B_R$ with $[\log R, +\infty) \times S^1$ (or $[\log R, +\infty) \times [0, \pi])$ if ${\mb A}_i = {\mb H}$) via the exponential map $(s, t)\mapsto e^{s + {\bm i} t}$, then it follows that $\check \vartheta_i$ is of class $W^{1,p, \delta_0 - 1 + \frac{2}{p}}$ with respect to the cylindrical metric. By the usual Sobolev embedding $W^{1,p}\hookrightarrow C^0$ in two dimensions, one sees that 
\beq\label{exponential}
|\check \vartheta_i (z)| \leq C |z|^{- (\delta_0 - \frac{\delta_p}{2})}.
\eeq

\subsection{Linearization of the marked disk}\label{subsection44}

The purpose of this subsection is to reformulate the linearization of the holomorphic disk. We would like to use the Euclidean metric on the upper half plane instead of the original metric on the disk. We also allow the disk to be deformed out of the symplectic quotient.

\subsubsection{Lifting the holomorphic disk}

We would like to lift the disk to a gauged map. Recall that one has a holomorphic disk $\bar u_0: ({\mb D}^2, \partial {\mb D}^2) \to (X, L)$ with interior marked points $z_1^+, \ldots, z_{m^+}^+ \in {\rm Int} {\mb D}^2$ and boundary marked points $z_1, \ldots, z_m$. Since the domain ${\mb D}^2$ is contractible, there exists a smooth map from ${\mb D}^2$ to $\mu^{-1}(0)$ whose composition with the projection $\mu^{-1}(0) \to X$ agrees with $\bar u_0$. Moreover, we choose the lift such that the matching condition
\beqn
\check u_i(\infty) = x_i = u_0(z_i),\ i = 1, \ldots, m + m^+
\eeqn
holds for all nodes (cf. \eqref{eqn412}). Then there exists a unique 1-form $a_0 \in \Omega^1({\mb D}^2, {\mf k})$ such that for all tangent vector $W \in T{\mb D}^2$ one has
\beq\label{eqn415}
du_0(W) + {\mc X}_{a_0(W)} \in (TV)_H.
\eeq
In this way the lift $u_0$ determines a gauged map from ${\mb D}^2$ to $V$, denoted by
\beqn
{\bm v}_0:= ( u_0, a_0).
\eeqn
We regard ${\mb D}$ as the completion ${\mb H} \cup \{\infty\}$ where ${\mb H}$ has the standard coordinate $z = s + {\bm i} t$. Using the standard coordinates $s$ and $t$ in ${\mb H}$ we write $\alpha_0 = \phi_0 ds + \psi_0 dt$. We also use ${\bm v}_0$ to denote the triple $(u_0, \phi_0, \psi_0)$, which satisfies the equation
\beqn
\begin{split}
\partial_s u_0 + {\mc X}_{\phi_0} + J (\partial_t u_0 + {\mc X}_{\psi_0})= &\ 0,\\
\mu(u_0 ) \equiv &\ 0,\\
u_0 (\partial {\mb H}) \subset &\  L_V.
\end{split}
\eeqn

We would like to give a linear operator which can be identified with the linearization of the holomorphic disk. Since $a_0$ is determined by $u_0$ by \eqref{eqn415}, the above system is equivalent to
\beq\label{eqn416}
\begin{split}
P_H ( \partial_s u_0 + J \partial_t u_0 ) = &\ 0,\\
\mu( u_0 ) \equiv &\ 0,\\
u_0 (\partial {\mb H}) \subset &\ L_V.
\end{split}
\eeq
The space of infinitesimal deformations of $u_0$ modulo gauge transformations is the space of sections of the bundle $u_0^* (TV)_H$. Then using the connection $\nabla$ specified by Lemma \ref{connection} one can write down the linearization of \eqref{eqn416}, which reads
\beqn
D_0^H(\xi_0) = P_H \big( \nabla_s \xi_0 +  (\nabla_{\xi_0} J) (\partial_t u_0)  + J \nabla_t \xi_0  \big).
\eeqn
Since $\nabla$ and $J$ commute with $P_H$, it is also equal to 
\beq\label{eqn417}
D_0^H (\xi_0) = \nabla_s \xi_0 + (\nabla_{\xi_0} J) (\partial_t u_0 + {\mc X}_{\psi_0}) + J \nabla_t \xi_0. 
\eeq

\subsubsection{Weighted Sobolev norms}

We would like to complete the space of infinitesimal deformations using a particular weighted Sobolev norm.  Choose a pair of half disks  
\begin{align*}
&\ U_0' = B_r^+(0) \subset {\mb H},\ &\ U_0	 = B_{r+1}^+(0)
\end{align*}
such that 
\begin{align*}
&\ B_1^+(z_i) \subset U_0'\ \forall i = 1,\ldots, m; &\ B_1 (z_j^+) \subset U_0'\ \forall j = 1, \ldots, m^+.
\end{align*}
Choose a smooth function $\rho_0: {\mb H} \to [1, +\infty)$ such that 
\beqn
\rho_0 (z) = \left\{ \begin{array}{cl} 1,\ &\ z \in U_0' \\
                                            |z|^2,\ &\ z \in {\mb H} \setminus U_0. 
                                            \end{array} \right.
\eeqn
Define weighted Sobolev norms $\| f \|_{W^{k,p,\delta}}$ for functions on ${\mb H}$ using the weight function $\rho_0$. Then define norms similar to \eqref{hnorm} and \eqref{gnorm}
\begin{align*}
&\ \| f \|_{L_G^{1,p,\delta}}:= \| f \|_{W^{1,p,\delta}},\ &\ \| f \|_{L_H^{1,p,\delta}}:= \| f \|_{L^\infty} + \| \nabla f \|_{L^{p,\delta}}.
\end{align*}
Using the connection $\nabla$ on $(TV)_H$ (resp. $\bar \nabla$ on $TX$), one can complete the space of infinitesimal deformations of $u_0$ (resp. $\bar u_0 = \pi_X \circ u_0$) to 
\beqn
L_H^{1,p,\delta}({\mb H}, u_0^* (TV)_H)_L \ \ {\rm resp.}\ \ L_H^{1,p,\delta}({\mb H}, \bar u_0^* TX)_L
\eeqn
where the subscript ${}_L$ indicates the Lagrangian boundary condition. 

Now we relate the linear operator $D_0^H$ with the linearization of the holomorphic disk $\bar u_0$. It is well-known that using the induced connection $\bar \nabla$ on $X$, the linearization of $\bar u_0$ is a Cauchy--Riemann operator 
\beqn
D_{\bar u_0} : W^{1,p}({\mb D}, \bar u_0^* TX)_L \to L^p({\mb D}, \Lambda^{0,1}\otimes \bar u_0^* TX)
\eeqn
which reads
\beq\label{eqn418}
D_{\bar u_0} (\bar \xi_0 ) = \left(  \nabla_s \bar \xi_0 + (\nabla_{\bar \xi_0 } I)( \partial_t \bar u_0) + I (\nabla_t \bar \xi_0) \right) \otimes d \bar z.
\eeq
Here $I$ is the domain-dependent almost complex structure on $X$ induced from the family $J$. To compare the two linear operators $D_0^H$ and $D_{\bar u_0}$, define a map between their domains
\begin{align*}
&\ \rho_1: W^{1,p}({\mb D}, \bar u_0^* TX)_L \to L_H^{1,p,\delta_p}({\mb H}, u_0^* (TV)_H)_L,\ &\ \rho_1^H (\bar \xi_0):= \bar \xi_0|_{{\mb D} \setminus \{\infty\}}
\end{align*}
and codomains
\begin{align*}
&\ \rho_0^H: L^p({\mb D}, \Lambda^{0,1} \otimes \bar u_0^* TX) \to L^{p,\delta_p}({\mb H}, u_0^*(TV)_H),\ &\ \rho_0^H( \nu \otimes d \bar z):= \nu|_{ {\mb D} \setminus \{\infty\}}.
\end{align*}

\begin{lemma}\label{lemma44}
The following facts hold.
\begin{enumerate}

\item The maps $\rho_1^H$ and $\rho_0^H$ are equivalences of Banach spaces.

\item There is a commutative diagram 
\beq\label{eqn419}
\xymatrix{     W^{1,p}({\mb D}^2, \bar u_0^* TX)_L  \ar[r]^-{\rho_1^H} \ar[d]_{D_{\bar u_0}} &  L_H^{1, p, \delta_p} ({\mb H}, u_0^* TX)_L \ar[d]^{D_0^H} \\
               L^p({\mb D}^2,  \bar u_0^* TX)   \ar[r]^-{\rho_0^H} &   L^{p, \delta_p} ({\mb H}, u_0^* TX).}
\eeq
\end{enumerate}
\end{lemma}

\begin{proof}
Let $w$ be the complex coordinate on ${\mb D}\subset {\mb C}$. We know that as $|z|\to \infty$, $|\frac{dz}{dw}|$ is comparable to $|z|^{-2}$ which is also comparable to $\rho_0^{-1}$. Then for $f \in L^p({\mb D})$, one has
\beqn
\| f \|_{L^p({\mb D})} \sim \int_{\mb H} |f|^p ds dt ( 1 + |z|^2)^{-2} \sim \int_{\mb H} |f|^p |\frac{dz}{dw}|^p (\rho_0(z))^{\frac{p \delta_p}{2}} ds dt  = \| \frac{dz}{dw} f \|_{L^{p,\delta_p}({\mb H})}.
\eeqn
Here $\sim$ means the two sides are comparable to each other. This calculation shows that $\rho_0^H$ is an equivalence of Banach spaces. Turn to the case of $\rho_1^H$. Let $\nabla$ be the differentiation on ${\mb H}$ with respect to the Euclidean coordinates and let $\nabla^\sim$ be the differentiation on ${\mb D}$. Then, using the same calculation and the Sobolev embedding, one can show that
\begin{multline*}
\| f \|_{L^\infty({\mb H})} + \| \nabla f \|_{L^{p,\delta_p}({\mb H})} \leq C \big( \| f \|_{L^\infty({\mb D})} +  \| \nabla^\sim f \|_{L^p({\mb D})} \big) \\
\leq C \big( \| f \|_{W^{1,p} ({\mb D})} + \| \nabla^\sim f \|_{L^p({\mb D})} \big) \leq C \big( \| f \|_{L^p({\mb D})} + \| \nabla^\sim f \|_{L^p({\mb D})} \big)\\
\leq C \big( \| f \|_{L^\infty({\mb D})} + \| \nabla^\sim f \|_{L^p({\mb D})} \big)  \leq C \big( \| f \|_{L^\infty({\mb H})} + \| \nabla f \|_{L^{p,\delta_p}({\mb H})} \big).
\end{multline*}
Here we followed the convention that the value of $C$ can vary from line to line. 

Lastly, the assertion that the diagram \eqref{eqn419} commutes follows from the comparison between the two expressions \eqref{eqn417} and \eqref{eqn418} and the fact that the connection $\bar \nabla$ resp. the almost complex structure $\bar I$ is induced from the restriction of the $K$-invariant connection $\nabla$ resp. $K$-invariant almost complex structure $J$ to the distribution $(TV)_H$. 
\end{proof}

\subsubsection{Deformation parameters}

The variation of the deformation parameter also deforms the Cauchy--Riemann equation. By the way we set the domain local model, the deformation parameter ${\bm w}$ does not only vary the target almost complex structure $J_{\bm w}$, but also the domain complex structure, given by the complex coordinate $f_{{\bm w}_0} = s_{{\bm w}_0} + {\bm i} t_{{\bm w}_0}$ (see \eqref{eqn27}). Therefore the linearization of the Cauchy--Riemann operator with respect to ${\bm w}$ is a finite rank operator
\beqn
\Phi_{\bar u_0}({\bm w})\in \Gammait({\mb H}, \Lambda^{0,1}\otimes \bar u_0^* T\bar X)
\eeqn
which is linear in $d \bar u_0$ and which is supported in a compact set of ${\mb H}$ disjoint from the marked points (i.e., the support of $df_{{\bm w}_0} - {\rm Id}$). Using the map $\rho_0^H$ one can identify it with a finite rank operator
\beq\label{eqn420}
\Phi_{u_0}^H({\bm w}) = \rho_0^H( \Phi_{\bar u_0}) \in \Gammait({\mb H}, u_0^* (TV)_H).
\eeq

\subsection{Rescaling the disk}\label{subsection45}

After the gluing, the disk component with a rescaling (and a small correction from the implicit function theorem) will become part of an ${\mb H}$-vortex. Using the $\epsilon$-rescaling \eqref{epsilonrescaling} we define a gauged map ${\bm v}_{0, \epsilon}:= (u_{0, \epsilon}, a_{0, \epsilon}) = (u_{0, \epsilon}, \phi_{0, \epsilon}, \psi_{0, \epsilon})$ on $\Sigma_0 \cong {\mb H}$ by 
\begin{align*}
&\ u_{0, \epsilon}: = s_\epsilon^* u_0,\ &\ a_{0, \epsilon} = \phi_{0,\epsilon} ds + \psi_{0,\epsilon} dt: = s_\epsilon^* a_0 = \epsilon s_\epsilon^* \phi_0 ds + \epsilon s_\epsilon^* \psi_0 dt.
\end{align*}
We regard ${\bm v}_{0, \epsilon}$ as an approximate affine vortex. 

We would like to consider the linear theory of the affine vortex equation at the gauged map ${\bm v}_{0,\epsilon}$ just defined. First we need certain $\epsilon$-dependent Sobolev norms. 

\begin{defn}\label{norms} Let $\epsilon>0$ be sufficiently small. 
\begin{enumerate}

\item For $f \in L^p_{\rm loc}({\mb H})$, define
\beqn
\| f \|_{L_\epsilon^{p,\delta_p}}:= \left[ \int_{\mb H}   | f(z)|^p \left( \frac{\rho_0 (\epsilon z)}{\epsilon} \right)^{p-2} ds dt \right]^{\frac{1}{p}}.
\eeqn

\item For $\xi^H \in W^{1,p}_{\rm loc}({\mb H}, u_{0,\epsilon}^* (TV)_H)_{L_V}$, define
\beqn
\| \xi^H \|_{L_{H, \epsilon}^{1,p,\delta_p}}:= \| \xi^H \|_{L^\infty} + \| \nabla^{a_{0,\epsilon}} \xi^H \|_{L_\epsilon^{p,\delta_p}}
\eeqn
and 
\beqn
{\mc B}_{0,\epsilon}^H:= \Big\{ \xi^H \in W^{1,p}_{\rm loc}({\mb H}, u_{0,\epsilon}^* (TV)_H)_{L_V}\ |\ 
\| \xi^H \|_{L_{H, \epsilon}^{1,p,\delta_p}} < \infty \Big\};
\eeqn
for ${\bm \xi}^G \in W^{1,p}_{\rm loc}({\mb H}, u_{0,\epsilon}^* (TV)_G \oplus {\mf k} \oplus {\mf k})_{L_V}$, define
\beqn
\| {\bm \xi}^G \|_{L_{G, \epsilon}^{1,p,\delta_p}}:= \| {\bm \xi}^G \|_{L_\epsilon^{p,\delta_p}} + \| \nabla^{a_{0,\epsilon}} {\bm \xi}^G \|_{L_\epsilon^{p,\delta_p}}.
\eeqn
and 
\beqn
{\mc B}_{0,\epsilon}^G:= \Big\{ {\bm \xi}^G \in W^{1,p}_{\rm loc}({\mb H}, u_{0,\epsilon}^* (TV)_G \oplus {\mf k} \oplus {\mf k})_{L_V}\ |\ 
\| {\bm \xi}^G \|_{L_{G, \epsilon}^{1,p,\delta_p}} < \infty \Big\}.
\eeqn
Define
\beqn
\hat{\mc B}_{0, \epsilon}:= {\mc B}_{0, \epsilon}^H \oplus {\mc B}_{0, \epsilon}^G
\eeqn
with the direct sum norm, which is equal to 
\beqn
\| {\bm \xi}\|_{L_\epsilon^{1,p, \delta_p}}:=  \|  \xi^H \|_{L^\infty} + \| {\bm \xi}^G \|_{L_\epsilon^{p,\delta_p}} + \| \nabla^{a_{0,\epsilon}} {\bm \xi} \|_{L_\epsilon^{p,\delta_p}}.
\eeqn
\end{enumerate}
On the other hand, define 
\beqn
\hat{\mc E}_{0,\epsilon}:= L_\epsilon^{p,\delta_p}( {\mb H}, u_{0,\epsilon}^* TV \oplus {\mf k} \oplus {\mf k}).
\eeqn
\end{defn}

We would like to invert the augmented linearization of the affine vortex equation at the rescaled gauged map ${\bm v}_{0,\epsilon}$. We first make a few calculations. Let $J: {\mb H}\to {\mc J}_{\rm tame}(V)$ be the restriction of the domain-dependent almost complex structure $J: \ov{\mc U}_{l, l^+} \to {\mc J}_{\rm tame}(V)$ to the disk component of the domain of ${\mc V}$ (see properties of such almost complex structures stated in Definition \ref{defn216}). By the expression of the augmented linearization (see \eqref{augmented}), the first coordinate of the augmented linearization reads
\begin{multline*}
\nabla_s^{a_{0,\epsilon}} \xi + (\nabla_\xi J ) (\partial_t u_{0, \epsilon} + {\mc X}_{\psi_{0, \epsilon}}) + J \nabla_t^{a_{0, \epsilon}} \xi + {\mc X}_\eta + J {\mc X}_\zeta \\
= \left[\begin{array}{c}  \nabla_s^{a_{0, \epsilon}} \xi^H + ( \nabla_\xi J ) (\partial_t u_{0, \epsilon} + {\mc X}_{\psi_{0,\epsilon}} ) + J \nabla_t^{a_{0, \epsilon}} \xi^H\\
\nabla_s^{a_{0, \epsilon}} \xi^G + J \nabla_t^{a_{0, \epsilon}} \xi^G + {\mc X}_\eta + J{\mc X}_\zeta \end{array} \right].
\end{multline*}
Notice that since $\nabla$ and $J$ respect the splitting $(TV)_G \oplus (TV)_H$, with respect to this splitting one can write
\beq\label{eqn421}
\hat D_{0,\epsilon} = \left[ \begin{array}{cc} D_{0, \epsilon}^H & E_{0,\epsilon} \\  0 & D_{0, \epsilon}^G \end{array} \right]
\eeq
where 
\beqn
D_{0,\epsilon}^H (\xi^H) =  \nabla_s^{a_{0, \epsilon}} \xi^H + ( \nabla_\xi J ) (\partial_t u_{0, \epsilon} + {\mc X}_{\psi_{0,\epsilon}} ) + J \nabla_t^{a_{0, \epsilon}} \xi^H,
\eeqn
\beq\label{eqn422}
E_{0,\epsilon}({\bm \xi}^G) = E_{0,\epsilon}(\xi^G, \eta, \zeta) = (\nabla_{\xi^G} J) (\partial_t u_{0,\epsilon} + {\mc X}_{\psi_{0,\epsilon}}),
\eeq
and 
\beq\label{eqn423}
D_{0,\epsilon}^G({\bm \xi}^G) = D_{0,\epsilon}^G (\xi^G, \eta, \zeta) = \left[ \begin{array}{l}  \nabla_s^{a_{0,\epsilon}} \xi^G + J \nabla_t^{a_{0,\epsilon}} \xi^G + {\mc X}_\eta + J {\mc X}_\zeta \\
\nabla_s^{a_{0,\epsilon}} \eta + \nabla_t^{a_{0,\epsilon}} \zeta + d\mu(u_{0,\epsilon}) (J_V \xi^G) \\
\nabla_s^{a_{0,\epsilon}} \zeta - \nabla_t^{a_{0,\epsilon}} \eta + d\mu(u_{0,\epsilon}) (\xi^G) \end{array} \right].
\eeq

The deformation parameter ${\bm w}$, which deforms both the domain and the target almost complex structure, gives rise to another zero-th order linear operator $\Phi_{{\bm v}_{0,\epsilon}}({\bm w})$ coming from \eqref{eqn420}, which is linear in $\partial_s u_{0,\epsilon} + {\mc X}_{\phi_{0,\epsilon}}$ and $\partial_t u_{0,\epsilon} + {\mc X}_{\psi_{0,\epsilon}}$. Because both $\partial_s u_{0,\epsilon} + {\mc X}_{\phi_{0,\epsilon}}$ and $\partial_t u_{0,\epsilon} + {\mc X}_{\psi_{0,\epsilon}}$ are in the distribution $(TV)_H$ and because of the property of the domain-dependent almost complex structure (see item (b) of Definition \ref{defn216}), one has 
\beqn
\Phi_{{\bm v}_{0,\epsilon}}({\bm w})  = \epsilon s_\epsilon^* \Phi_{u_0}^H({\bm w}) \in (TV)_H.
\eeqn
Therefore, the total augmented linearization (including the effect of deformation parameters) of the affine vortex equation at ${\bm v}_{0,\epsilon}$ reads
\begin{align*}
&\ \hat{\mc D}_{0, \epsilon}: W_{\rm def} \oplus \hat{\mc B}_{0,\epsilon} \to \hat{\mc E}_{0,\epsilon},\ &\  \hat{\mc D}_{0,\epsilon}({\bm w}, {\bm \xi}) = \Phi_{{\bm v}_{0,\epsilon}}({\bm w}) + \hat D_{0,\epsilon}({\bm \xi}).
\end{align*}
We can see that from the above calculation, the ``horizontal'' part of the linearization satisfies 
\beq\label{eqn424}
D_{0,\epsilon}^H(s_\epsilon^* \xi^H) + \Phi_{{\bm v}_{0,\epsilon}}({\bm w}) = \epsilon s_\epsilon^* \big( D_0^H(\xi^H) + \Phi_{u_0}^H({\bm w}) \big),\ \forall {\bm w} \in W_{\rm def},\ \xi^H \in \Gammait(u_0^* (TV)_H).
\eeq
The linear theory of the holomorphic disk then can be translated to the horizontal part. On the other hand, the group action part of the augmented linearization is treated in the following lemma. 

\begin{lemma}\label{lemma46}
There exist $\epsilon_0 > 0 $ and $C >0$ such that for all $\epsilon \in (0, \epsilon_0)$, the operator $D_{0,\epsilon}^G: {\mc B}_{0,\epsilon}^G \to {\mc E}_{0, \epsilon}^G$ has a bounded inverse $D_{0,\epsilon}^G$. Moreover, with respect to the norms defined in Definition \ref{norms}, one has
\beqn
\| Q_{0,\epsilon}^G \| \leq C,\ \forall \epsilon \in (0, \epsilon_0).
\eeqn
\end{lemma}
	
\begin{proof}
We define an $\epsilon$-dependent norm on sections of $(u_0^* TV)_G \oplus {\mf k} \oplus {\mf k}$ by
\beq\label{eqn425}
\| {\bm \xi}^G \|_{1,p;\epsilon}^{\rm aux}:= \| {\bm \xi}^G \|_{L^{p,\delta_p}} + \epsilon \| \nabla^{a_0} {\bm \xi}^G \|_{L^{p,\delta_p}}.
\eeq
By straightforward calculation one has
\begin{align}\label{eqn426}
&\ \| s_\epsilon^* {\bm \xi}^G \|_{L_{G,\epsilon}^{1,p,\delta_p}} = \epsilon^{-1} \| {\bm \xi}^G \|_{1,p;\epsilon}^{\rm aux},\ &\ \| s_\epsilon^* {\bm \nu}^G \|_{L_\epsilon^{p,\delta_p}} = \epsilon^{-1} \| {\bm \nu}^G \|_{L^{p,\delta_p}}.
\end{align}
Hence it suffices to consider the conjugated operator
\beqn
( s_\epsilon^*)^{-1} \circ D_{0,\epsilon}^G \circ s_\epsilon^*: W^{1,p,\delta_p}({\mb H}, u_0^* (TV)_G \oplus {\mf k} \oplus {\mf k})_L \to L^{p,\delta_p}({\mb H}, u_0^* (TV)_G \oplus {\mf k} \oplus {\mf k}).
\eeqn
The formula \eqref{eqn423} shows that after the above conjugation, one can write 
\beqn
( s_\epsilon^*)^{-1} \circ D_{0,\epsilon}^G \circ s_\epsilon^* = \left[ \begin{array}{cc} \epsilon D_0^G & L_0 \\  L_0^* & \epsilon \s D_0^G \end{array} \right]
\eeqn
where 
\begin{align*}
&\ D_0^G (\xi^G) = \nabla_s^{a_0} \xi^G + J \nabla_t^{a_0} \xi^G,\ &\ 
\s D_0^G (\eta, \zeta) = ( \nabla_s^{a_0} \eta + \nabla_t^{a_0} \zeta, \nabla_s^{a_0}\zeta - \nabla_t^{a_0} \eta),
\end{align*}
\begin{align*}
&\ L_0 (\eta, \zeta) = {\mc X}_\eta + J{\mc X}_\zeta,\ &\ L_0^* (\xi^G) = (d\mu(u_0)(J_V \xi^G), d\mu(u_0)(\xi^G)).
\end{align*}

We would to rewrite this operator as a differential operator on matrix-valued functions. Define a trivialization $\tau_G: {\mb H} \times  ({\mf k} \otimes {\mb C}) \cong u_0^* (TV)_G$ by 
\beqn
\tau_G ( f' + {\bm i} f'') = {\mc X}_{f'} + J {\mc X}_{f''}.
\eeqn
Then one has 
\begin{multline*}
\tau_G^{-1} \circ D_0^G \circ \tau_G (f' + {\bm i} f'') \\
= \tau_G^{-1}( \nabla_s ({\mc X}_{f'} + J {\mc X}_{f''}) + \nabla_{{\mc X}_{f'} + J {\mc X}_{f''}} {\mc X}_{\phi_0} + J(\nabla_t ({\mc X}_{f'} + J {\mc X}_{f''}) + \nabla_{{\mc X}_{f'} + J {\mc X}_{f''}} {\mc X}_{\psi_0} ) \\
= 2 \partial_{\bar z} (f' + {\bm i} f'') + B(z)(f' + {\bm i} f'')
\end{multline*}
where $B(z)$ is a zero-th order operator with $B(z) \to 0$ as $z\to \infty$. Moreover, 
\beqn
\tau_G^{-1} \circ L_0( \eta + {\bm i} \zeta) = \eta + {\bm i} \zeta
\eeqn
and 
\beqn
L_0^* \circ \tau_G (f' + {\bm i} f'') = A(z) (f' + {\bm i} f'')
\eeqn
where $A(z)$ is a family of Hermitian matrices whose eigenvalues are bounded from below through ${\mb H}$. Then with respect to the trivialization $\tau_G$, one can write $(s_\epsilon^*)^{-1}\circ D_{0,\epsilon}^G \circ s_\epsilon^*$  as 
\beqn
\left[ \begin{array}{cc} 2\epsilon \partial_{\bar z} & {\rm Id}_{{\mf k}\otimes {\mb C}} \\ A(z) & 2\epsilon \partial_z \end{array}\right] + \epsilon B(z)
\eeqn
Then this lemma follows from Lemma \ref{model} below.\end{proof}

\begin{lemma}\label{model}
Given $l \geq 1$. Let 
\beqn
A_1, A_2: {\mb H} \to {\mb C}^{l\times l}
\eeqn
be a pair of matrix-valued functions satisfying the following conditions.
\begin{enumerate}
\item $A_1(z)$, $A_2(z)$ are positive Hermitian matrices for all $z$, converge as $z \to \infty$, and commute for all $z$.

\item $|\nabla A_1 (z)|$ and $|\nabla A_2 (z)|$ are bounded.

\item The smallest eigenvalues of $A_1 (z)$ and $A_2 (z)$ are bounded away from zero.
\end{enumerate}
Let $B: {\mb H} \to {\rm End}_{\mb R}({\mb C}^{2l})$ be a matrix-valued continuous function which converges to zero as $z \to \infty$. For each $\epsilon>0$, define 
\beqn
D^\epsilon:  = \left[ \begin{array}{cc} 2\epsilon \partial_{\bar z} & A_1 \\ A_2 & 2 \epsilon \partial_z  \end{array}\right] + \epsilon B. 
\eeqn
Then there exists $\epsilon_0>0$ and $C >0$ such that for all $\epsilon \in (0, \epsilon_0)$, $D_\epsilon$ has a bounded inverse $Q^\epsilon$ satisfying
\beqn
\|Q^\epsilon( h ) \|_{L^{p,\delta_p}} + \epsilon \| \nabla Q^\epsilon(h) \|_{L^{p,\delta_p}} \leq C \| h \|_{L^{p,\delta_p}}.
\eeqn
\end{lemma}

\begin{proof}
Since $B$ is uniformly bounded, it suffices to consider the situation when $B = 0$. We first prove for the situation when the weighted Sobolev norms are replaced by the ordinary Sobolev norms. By a conjugation trick (see \cite[p.59]{Ziltener_book}), $D^\epsilon$ is conjugated to 
\beqn
D_A^\epsilon = \left[ \begin{array}{cc} 2 \epsilon \partial_{\bar z} & \sqrt{ A_1 } \sqrt{A_2} \\ \sqrt{A_2} \sqrt{A_1} & 2 \epsilon \partial_z \end{array}\right].
\eeqn
Since $A_1$ and $A_2$ commute, we may assume that $A_1 = A_2 = A = A^*$. Using the $\epsilon$-rescaling, define
\beqn
A^\epsilon(z):= (s_\epsilon^* A)(z) = A (\epsilon z).
\eeqn
Then it suffices to show that the operator
\beqn
D_{A^\epsilon}:= s_\epsilon^*  \circ D^\epsilon \circ (s_\epsilon^*)^{-1} = \left[ \begin{array}{cc} 2\partial_{\bar z} & A^\epsilon \\ A^\epsilon & 2 \partial_z \end{array}\right]
\eeqn
has an inverse $Q$ such that 
\beq\label{eqn427}
\| Q (h) \|_{L^p} + \| \nabla Q(h) \|_{L^p} \leq C \| h \|_{L^p}
\eeq
for some $C>0$. To see this, consider the formal adjoint of $D_{A^\epsilon}$ which reads
\beqn
(D_{A^\epsilon})^* = \left[ \begin{array}{cc} -2 \partial_z & A^\epsilon \\ A^\epsilon & - 2 \partial_{\bar z} \end{array}\right].
\eeqn
Then we see that  
\beqn
(D_{A^\epsilon})^* D_{A^\epsilon} \left[ \begin{array}{c} f_1 \\ f_2 \end{array} \right] = \left[\begin{array}{c} -\Delta f_1 + (A^\epsilon)^2 f_1 \\ - \Delta f_2 + (A^\epsilon)^2 f_2  \end{array} \right] + T_{A^\epsilon}
\eeqn
where $T_{A^\epsilon}$ is a zero-th order operator coming from derivatives of $A^\epsilon$. Then 
\beqn
\| T_{A^\epsilon}\|_{L^\infty} \leq C \epsilon
\eeqn
for a certain $C>0$. Since the eigenvalues of $A$ are uniformly bounded from below, $(D_{A^\epsilon})^* D_{A^\epsilon}$ is strictly positive definite. Moreover, we know that $D_{A^\epsilon}$ is Fredholm with index zero. Hence $D_{A^\epsilon}$ has an inverse satisfying \eqref{eqn427}. This proves this lemma if we replace the weighted norms $L^{p,\delta_p}$ by the ordinary norm $L^p$.  

To prove the original statement, consider the map
\beqn
m_\delta: f \mapsto \rho_0^{-\delta}  f.
\eeqn
This is an equivalence of norms between the un-weighted and weighted Sobolev spaces. Then consider the conjugated operator between the un-weighted Sobolev spaces	
\beqn
(m_\delta^{-1} \circ D^\epsilon \circ m_\delta) \left[ \begin{array}{c} f_1 \\ f_2 \end{array}\right] = D^\epsilon \left[ \begin{array}{c} f_1 \\ f_2 \end{array}\right] + 2\epsilon  \left[ \begin{array}{c} \rho_0^\delta (\partial_{\bar z} \rho_0^{-\delta}) f_1 \\ \rho_0^\delta (\partial_z \rho_0^{-\delta}) f_2 \end{array} \right].
\eeqn
The last term is bounded by a multiple of $\epsilon$ as $|\rho_0^\delta \nabla \rho_0^{-\delta}|$ is uniformly bounded. Hence when $\epsilon$ is small enough, the assertion of this lemma follows from the un-weighted case. 
\end{proof}

\subsection{The singular configuration}\label{subsection46}

In this subsection we describe the inverse of the linearization map for an $\epsilon$-dependent singular configuration. Recall that one has the Banach space $\hat{\mc B}_{0, \epsilon}$ parametrizing certain infinitesimal deformations of the rescaled gauged map ${\bm v}_{0,\epsilon}$. Using the exponential map of $h_V$ we may regard the space $\hat{\mc B}_{0,\epsilon}$ (more precisely a small neighborhood of the origin) as gauged maps from ${\mb H}$ to $V$ near ${\bm v}_{0,\epsilon}$. Moreover one has the Banach spaces $\hat{\mc B}_{{\bm v}_i}$ of infinitesimal deformations of ${\bm v}_i$, which can also be viewed as a space of genuine gauged maps near ${\bm v}_i$. Denote 
\beqn
\hat{\mc B}_{\Gamma, \epsilon}^\#:= \hat{\mc B}_{0, \epsilon} \oplus \bigoplus_{i=1}^{m+m^+} \hat{\mc B}_{{\bm v}_i}. 
\eeqn
Using the evaluation map at infinities and at nodal points, one has the subspace of configurations satisfying the matching conditions:
\beqn
\hat {\mc B}_{\Gamma, \epsilon}:= \big\{ ({\bm \xi}_0, {\bm \xi}_1, \ldots, {\bm \xi}_{m+m^+} ) \in \hat{\mc B}_{\Gamma, \epsilon}^\# \ |\ {\rm the\ ``matching\ condition"} \big\}
\eeqn
where the matching condition is defined as the linearization of the matching condition in defining the Banach manifold $\hat {\mc B}_\Gamma$. More precisely, the matching condition is 
\beqn
{\bm \xi}_i(\infty) = {\bm \xi}_0 ( z_i ) \in T_{\bar x_i} X,\ i = 1, \ldots, m + m^+
\eeqn
where $\bar x_i \in X$ is the projection image of $x_i \in \mu^{-1}(0)$ (cf. \eqref{eqn412}). Taking the direct sum of the codomains of the augmented linearized operators provides another Banach space
\beqn
\hat {\mc E}_{\Gamma, \epsilon}:= \hat{\mc E}_{0, \epsilon} \oplus \bigoplus_{i=1}^{m+m^+} \hat {\mc E}_{{\bm v}_i}.
\eeqn
The direct sum of augmented linearizations at the singular configuration together with deformations of the almost complex structure defines a linear operator
\beqn
\hat {\mc D}_{\Gamma, \epsilon}: W_{\rm def} \oplus \hat{\mc B}_{\Gamma, \epsilon} \to \hat{\mc E}_{\Gamma, \epsilon}.
\eeqn

The following proposition concludes this section.

\begin{prop}\label{inverse}
There exist $\epsilon(Q) > 0$ and $c(Q) > 0$ such that for all $\epsilon \in (0, \epsilon(Q))$, the operator $\hat {\mc D}_{\Gamma, \epsilon}$ has an inverse $\hat {\mc Q}_{\Gamma, \epsilon}$. Moreover, with respect to the $\epsilon$-dependent norms on $\hat {\mc B}_{\Gamma, \epsilon}$ and $\hat{\mc E}_{\Gamma, \epsilon}$, there holds
\beqn
\| \hat {\mc Q}_{\Gamma, \epsilon} \| \leq c(Q).
\eeqn
\end{prop}

\begin{proof}
Recall the transversality assumption of Theorem \ref{mainthm}. Abbreviate the linearization over the $i$-th component (including the disk component) by 
\beqn
D_i: {\mc B}_i \to {\mc E}_i, i = 0, 1, \ldots, m+m^+
\eeqn
and the extension including deformation parameters by 
\beqn
{\mc D}_i: W_{\rm def} \oplus {\mc B}_i \to {\mc E}_i.
\eeqn
Define ${\mc D}_\Gamma$ to be the restriction of the direct sum 
\beqn
{\mc D}_\Gamma^\#:= \bigoplus_{i=0}^{m+m^+} {\mc D}_i: W_{\rm def} \oplus \bigoplus_{i=0}^{m+m^+} {\mc B}_i \to \bigoplus_{i=0}^{m+m^+} {\mc E}_i
\eeqn
to the finite codimensional subspace ${\mc B}_\Gamma$ defined by the matching condition. The transversality assumption of Theorem \ref{mainthm} implies the existence of a bounded inverse
\beqn
{\mc Q}_\Gamma: \bigoplus_{i=0}^{m+m^+} {\mc E}_i \to {\mc B}_\Gamma \hookrightarrow W_{\rm def} \oplus \bigoplus_{i=0}^{m+m^+} {\mc B}_i.
\eeqn

We first extend ${\mc Q}_\Gamma$ over the affine vortex components to get right inverses to the augmented linearizations. For $i = 1, \ldots, m+m^+$, let ${\mc B}_i^\bot \subset \hat{\mc B}_i$ be the $L^2$-orthogonal complement of the Coulomb slice ${\mc B}_i$. Then the Coulomb gauge-fixing condition defines a bounded invertible operator (called the Coulomb operator)
\begin{align*}
&\ {\mc B}_i^\bot \to {\mc E}_i^\bot:= L^{p,\delta_p}({\mb A}_i, {\mf k}),\ &\ (\xi, \eta, \zeta) \mapsto \nabla_s^{a_i} \eta + \nabla_t^{a_i} \zeta + d\mu(u_i)(J_V \xi).
\end{align*}
Let 
\beqn
Q_i^\bot: {\mc E}_i^\bot \to {\mc B}_i^\bot.
\eeqn
be the inverse of the Coulomb operator. Then define 
\beqn
Q_\Gamma^\bot: \bigoplus_{i=1}^{m+m^+} {\mc E}_i^\bot \to \bigoplus_{i=1}^{m+m^+} {\mc B}_i^\bot. 
\eeqn
Then denote the direct sum of ${\mc Q}_\Gamma$ and $Q_\Gamma^\bot$ by
\beqn
\hat {\mc Q}_\Gamma': {\mc E}_0 \oplus \bigoplus_{i=1}^{m+m^+} \hat{\mc E}_i \to W_{\rm def} \oplus {\mc B}_0 \oplus \bigoplus_{i=1}^{m+m^+} \hat{\mc B}_i.
\eeqn
Its image lies in the subspace defined by the matching condition. 

Using the horizontal maps of the commutative diagram \eqref{eqn419} and the $\epsilon$-rescaling, one can identify
\begin{align*}
&\ s_\epsilon^* \circ \rho_1^H: {\mc B}_0 \cong {\mc B}_{0,\epsilon}^H,\ &\ \epsilon s_\epsilon^* \circ \rho_0^H: {\mc E}_0 \cong {\mc E}_{0,\epsilon}^H.
\end{align*}
By Lemma \ref{lemma44}, these identifications are equivalences of Banach spaces making the norms uniformly comparable. Then after these identifications, one can regard $\hat {\mc Q}_\Gamma'$ as an operator
\beqn
\hat {\mc Q}_{\Gamma, \epsilon}': {\mc E}_{0, \epsilon}^H \oplus \bigoplus_{i=1}^{m+m^+} \hat{\mc E}_i \to W_{\rm def} \oplus {\mc B}_{0, \epsilon}^H \oplus \bigoplus_{i=1}^{m+m^+} \hat{\mc B}_i.
\eeqn
Moreover, over the holomorphic disk component, Lemma \ref{lemma46} shows that there is a uniformly bounded inverse
\beqn
Q_{0, \epsilon}^G: {\mc E}_{0, \epsilon}^G \to {\mc B}_{0,\epsilon}^G.
\eeqn
to the operator $D_{0, \epsilon}^G$. Recall that $\hat{\mc B}_{0, \epsilon} = {\mc B}_{0,\epsilon}^H \oplus {\mc B}_{0, \epsilon}^G$ and $\hat{\mc E}_{0, \epsilon} = {\mc E}_{0, \epsilon}^H \oplus {\mc E}_{0,\epsilon}^G$. Then one can form the direct sum
\beqn
\hat {\mc Q}_{\Gamma, \epsilon}' \oplus Q_{0, \epsilon}^G: \hat{\mc E}_{0, \epsilon} \oplus \bigoplus_{i=1}^{m+m^+} \hat{\mc E}_i \to W_{\rm def} \oplus \hat{\mc B}_{0,\epsilon} \oplus \bigoplus_{i=1}^{m+m^+} \hat{\mc B}_i.
\eeqn
By the construction, its image lies in the finite codimensional subspace defined by the matching condition at nodes. We can see that if the off-diagonal term $E_{0, \epsilon}$ in the expression of $\hat D_{0,\epsilon}$ (see \eqref{eqn421} and \eqref{eqn422}) is zero, then the above direct sum is an inverse of $\hat {\mc D}_{\Gamma, \epsilon}$ (cf. \eqref{eqn424}), which has a uniform upper bound on the norm. Since the off-diagonal term is uniformly bounded, one can correct the direct sum to obtain a true inverse $\hat{\mc Q}_{\Gamma, \epsilon}$ whose norm is uniformly bounded. This finishes the proof.
\end{proof}

\section{The Gluing Construction}\label{section5}

In this section we apply the standard gluing protocol (with modifications) to construct a family of smooth affine vortices near the limiting singular configuration. The basic procedure is as follows. In Subsection \ref{subsection51} we construct a family of approximate solutions parametrized by the gluing parameter $\epsilon$. In Subsection \ref{subsection52} we introduce an $\epsilon$-dependent weighted Sobolev norms along the approximate solutions. Inn Subsection \ref{subsection53} we state the major estimates required to apply the implicit function theorem, which immediately implies the (set-theoretic) construction of the gluing map. Lastly, in Subsection \ref{subsection54}, \ref{subsection55}, and \ref{subsection56} we prove the major estimates stated in Subsection \ref{subsection53}. 

We still follow the convention that the letter $C$ represents a constant which is allowed to vary its value from line to line. 

\subsection{Pregluing}\label{subsection51}

\subsubsection{The cut-off functions}

The glued affine vortex approximately agrees with the singular affine vortex on different regions. We first describe these regions. Recall that the domain of the holomorphic disk component is identified with a copy of the upper half plane ${\mb H}$ with marked points
\begin{align*}
&\ z_1, \ldots, z_m \in \partial {\mb H},\ &\ z_1^+ = z_{m+1}, \ldots, z_{m^+}^+ = z_{m+ m^+} \in {\rm Int}{\mb H}.
\end{align*}
Take a number $b >0$ such that 
\beq\label{eqn51}
100 c(p,\delta_p) c(Q) \leq \log b \leq 1000 c(p,\delta_p ) c(Q)
\eeq
where $c(p, \delta_p)$ is the constant for the Hardy-type inequality \eqref{hardy} for $\delta = \delta_p$ and $c(Q)$ is the upper bound of the operator norm of $\hat {\mc Q}_{\Gamma, \epsilon}$ given by Proposition \ref{inverse}. Define
\beqn
z_{i,\epsilon}:= \frac{z_i}{\epsilon},\ \hspace{2cm}\ \Sigma_{i,\epsilon}:= B\left( z_{i,\epsilon}, \frac{1}{\sqrt\epsilon} \right),\ \hspace{2cm}\ i = 1, \ldots, m+m^+.
\eeqn
(Here $B(z, r) \subset {\mb A}$ means the radius $r$ open (half) disk centered at $z$.) Define their enlargements/shrinkings by
\beqn
\begin{split}
\Sigma_{i, \epsilon}^>:= B\left( z_{i,\epsilon}, \frac{b}{\sqrt\epsilon} \right),\ &\ \Sigma_{i, \epsilon}^\gg:= B\left( z_{i,\epsilon}, \frac{2b}{\sqrt\epsilon} \right),\\
\Sigma_{i,\epsilon}^<:= B\left( z_{i,\epsilon}, \frac{1}{b \sqrt\epsilon} \right),\ &\ \Sigma_{i,\epsilon}^\ll:= B \left( z_{i,\epsilon}, \frac{1}{2b\sqrt{\epsilon}} \right).
\end{split}
\eeqn
Define
\beqn
\Sigma_{0, \epsilon}:= {\mb H} \setminus \bigcup_{i=1}^{m+m^+} \Sigma_{i,\epsilon}
\eeqn
and
\begin{align*}
&\ \Sigma_{0,\epsilon}^>:= {\mb H} \setminus \bigcup_{i=1}^{m+m^+} \Sigma_{i,\epsilon}^<,\ &\ \Sigma_{0,\epsilon}^<:= {\mb H} \setminus \bigcup_{i=1}^{m+m^+} \Sigma_{i,\epsilon}^>,
\end{align*}
\begin{align*}
&\ \Sigma_{0,\epsilon}^\gg:= {\mb H} \setminus \bigcup_{i=1}^{m+m^+} \Sigma_{i,\epsilon}^\ll,\ &\ \Sigma_{0,\epsilon}^\ll:= {\mb H} \setminus \bigcup_{i=1}^{m+m^+} \Sigma_{i,\epsilon}^\gg
\end{align*}

We need certain cut-off functions to interpolate different components. Let $\gamma: {\mb R} \to [0, 1]$ be a smooth cut-off function such that 
\beqn
\gamma |_{(-\infty, -1]} \equiv 1,\hspace{1cm} \gamma |_{[0, +\infty)} = 0,\hspace{1cm} |\nabla \gamma | \leq 2.
\eeqn
Then for a given gluing parameter $\epsilon$ define
\beqn
 \gamma_i^\epsilon (z):= \gamma \left( \frac{\log |z - z_{i, \epsilon}| +  \log b + \log \sqrt\epsilon }{\log 2}\right),\ i = 1, \ldots, m + m^+.
\eeqn
Define $\gamma_0^\epsilon: {\mb H} \to [0, 1]$ by 
\beqn
\gamma_0^\epsilon (z) = \left\{ \begin{array}{cc} \displaystyle \gamma \left( \frac{ - \log |z- z_{i, \epsilon}| + \log b - \log \sqrt\epsilon }{\log 2}\right),&\ z \in \Sigma_{i,\epsilon}^\gg, \\ 1, &\ z \in \Sigma_{0, \epsilon}^\ll. \end{array}\right.
\eeqn
Then for $i = 0, \ldots, m + m^+$, $\gamma_i^\epsilon$ equals to one inside $\Sigma_{i, \epsilon}^\ll$ and  and equals to zero outside $\Sigma_{i,\epsilon}^<$. We see that 
\beq\label{eqn52}
\|\nabla \gamma_i^\epsilon \|_{L^\infty} \leq  \frac{2 b \sqrt\epsilon }{\log 2} \sup |\nabla \gamma | \leq 2b \sqrt\epsilon,\ i = 0, \ldots, m + m^+.
\eeq

\subsubsection{The approximate solutions}

For each $\epsilon$ small enough, we would like to define a gauged map ${\bm v}_\epsilon = (u_\epsilon, \phi_\epsilon, \psi_\epsilon)$ on ${\mb H}$ from the components ${\bm v}_i$ and ${\bm v}_0$. We first need to choose a different gauge of the gauged map representing the holomorphic disk. Let $(r_i, \theta_i)$ be the polar coordinates centered at $z_i$. Let 
\beqn
g_0: {\mb H} \smallsetminus \{ z_1^+, \ldots, z_{m^+}^+ \} \to K
\eeqn
be a gauge transformation satisfying the following conditions
\begin{enumerate}
\item For $i = 1, \ldots, m$, $g_0 ( r_i, \theta_i) = e^{- \lambda_i^+ \theta_i}$ in a small neighborhood of $z_i$. 

\item $g_0$ equals identity outside a compact subset of ${\mb H}$ and 
\beqn
* dg_0 |_{\partial {\mb H}} = 0.
\eeqn 
\end{enumerate}
We will replace the gauged map ${\bm v}_0 = (u_0, a_0)$ by $g_0 \cdot {\bm v}_0 = (g_0 \cdot u_0, g_0 \cdot a_0)$. After the gauge transformation, the connection form has singularities at the interior markings. For convenience, define $\lambda_i = 0$ for $j = 1, \ldots, m$. Then since $u_0$ has value $x_i \in \mu^{-1}(0)$ at $z_i$, one has 
\beqn
\lim_{z \to z_i} e^{\lambda_i \theta_i} u_0 (z) = x_i,\ i=1, \ldots, m + m^+.
\eeqn
Moreover, the gauged map ${\bm v}_0$ satisfies the boundary condition 
\begin{align*}
&\ u_0 (\partial \Sigma) \subset L_V,\ &\ * a_0 |_{\partial \Sigma} = 0.
\end{align*}
One can still use the $\epsilon$-rescaling $s_\epsilon$ to pull back ${\bm v}_{0, \epsilon}$ to a family of gauged maps ${\bm v}_{0, \epsilon} = (u_{0, \epsilon}, a_{0, \epsilon})$ defined over 
\beqn
{\mb H} \setminus \{ z_{1, \epsilon}^+, \ldots, z_{m^+, \epsilon}^+\}.
\eeqn
The original (smooth) gauged map is denoted by $\check{\bm v}_0 = (\check u_0, \check a_0)$ and the corresponding rescaled ones denoted by $\check {\bm v}_{0, \epsilon} = (\check u_{0, \epsilon}, \check a_{0, \epsilon})$.

Now we describe the pregluing construction which gives an approximate solution. Define a gauged map ${\bm v}_\epsilon = (u_\epsilon, a_\epsilon)$ from ${\mb H}$ to $V$ as follows. Recall notations used in Lemma \ref{lemma43}.

\begin{enumerate}

\item The connection part $a_\epsilon$ is defined as follows. Over the neck region $A_{i, \epsilon}:= \Sigma_{i,\epsilon}^\gg \smallsetminus \Sigma_{i,\epsilon}^\ll$ with polar coordinates $(r_i, \theta_i)$. Then define
\beq\label{eqn53}
a_\epsilon = \left\{ \begin{array}{cl}
                             a_{0, \epsilon}(z),\ &\ z \in \Sigma_{0, \epsilon}^\ll;\\
														 e^{- \lambda_i \theta_i} \cdot \Big[ \gamma_0^\epsilon \check a_{0, \epsilon}(r_i, \theta_i) + \gamma_i^\epsilon \check a_i (r_i, \theta_i) \Big]	,\ &\ z \in A_{i,\epsilon};\\
														  a_i(z - z_{i, \epsilon}),\ &\  z\in \Sigma_{i,\epsilon}^\ll.
															\end{array}\right.
															\eeq

\item The map part $u_\epsilon$ is defined as follows. Recall one has $\check{u}_i = e^{\lambda_i \theta_i} u_i$ and 
\beqn
\check{u}_{0, \epsilon}( z_{i, \epsilon} + r_i e^{\i \theta_i}) = \exp_{x_i} \check\vartheta_{0, \epsilon}(r_i, \theta_i),\ \check{u}_i(r_i, \theta_i) = \exp_{x_i} \check \vartheta_i(r_i, \theta_i),\ \frac{1}{2b \sqrt\epsilon} \leq r_i \leq \frac{2b}{\sqrt\epsilon}.
\eeqn
Define
\beq\label{eqn54}
u_\epsilon(z) = \left\{ \begin{array}{cl} u_{0, \epsilon}(z),&\ z \in \Sigma_{0, \epsilon}^\ll;\\
e^{- \lambda_i \theta_i} \exp_{x_i} \big( \gamma_0^\epsilon \check\vartheta_{0, \epsilon}(r_i, \theta_i) + \gamma_i^\epsilon \check\vartheta_i(r_i, \theta_i) \big),\ &\  z \in A_{i,\epsilon};\\
                                           u_i( z - z_{i, \epsilon} ),\ &\ z\in \Sigma_{i,\epsilon}^\ll.
																																									 
\end{array} \right.
\eeq
Notice that since $L_V$ is totally geodesic with respect to $h_V$, the map $u_\epsilon$ satisfies the Lagrangian boundary condition $u_\epsilon(\partial {\mb H}) \subset L_V$. 
\end{enumerate}

Over a neck region the approximate solution agrees with a ``constant'' gauged map. Indeed, over the smaller neck region $A_{i,\epsilon}':= \Sigma_{i, \epsilon}^> \setminus \Sigma_{i,\epsilon}^<$ one has
\beq\label{eqn55}
{\bm v}_\epsilon|_{A_{i,\epsilon}'} = {\bm c}_i:= (e^{- \lambda_i \theta_i} x_i, \lambda_i d \theta_i ),\ i = 1, \ldots m + m^+.
\eeq
We introduce a collection of auxiliary gauged maps ${\bm v}_{i, \epsilon}'$ (defined over ${\mb A}_i$) and ${\bm v}_{0, \epsilon}'$ (defined over $\Sigma_{\epsilon}^*$):
\begin{align}\label{eqn56}
&\ {\bm v}_{i, \epsilon}' = \left\{ \begin{array}{cl} {\bm v}_\epsilon,\ &\ {\rm on}\ \Sigma_{i, \epsilon},\\
                                                    {\bm c}_i,\ &\ {\rm on}\ {\mb A}_i \setminus \Sigma_{i,\epsilon};
\end{array} \right.\ &\ {\bm v}_{0, \epsilon}'= \left\{ \begin{array}{cl} {\bm v}_\epsilon,\ &\ {\rm on}\ \Sigma_{0,\epsilon},\\
                                                                              {\bm c}_i,\ &\ {\rm on}\ \Sigma_{i, \epsilon},\ i=1, \ldots, m + m^+.
\end{array} \right.
\end{align}
Notice that ${\bm v}_{i, \epsilon}'$ is very close to ${\bm v}_i$ and ${\bm v}_{0, \epsilon}'$ is very close to ${\bm v}_{0, \epsilon}$. More precisely, one has the following fact. 

\begin{lemma}\label{lemma51}
If we identify ${\bm v}_{i, \epsilon}'$ with a point ${\bm \xi}_{i, \epsilon}'\in \hat{\mc B}_i$ and identify ${\bm v}_{0, \epsilon}'$ with a point ${\bm \xi}_{0, \epsilon}'\in \hat{\mc B}_{0, \epsilon}$, then 
\beqn
\lim_{\epsilon \to 0} \| {\bm \xi}_{i, \epsilon}' \|_{L_m^{1,p, \delta_p}} = \lim_{\epsilon \to 0} \| {\bm \xi}_{0, \epsilon}'\|_{L_{m, \epsilon}^{1,p,\delta_p}} =0.
\eeqn
\end{lemma}

\begin{proof}
By the definition of ${\bm v}_{i,\epsilon}$, one has 
\begin{align*}
&\ e^{\lambda_i \theta_i} \cdot {\bm v}_i = ( \exp_{x_i} \check \vartheta_i, \check a_i),\ &\ e^{\lambda_i \theta} \cdot {\bm v}_{i,\epsilon} = ( \exp_{x_i} \gamma_i^\epsilon  \check \vartheta_i, \gamma_i^\epsilon \check a_i).
\end{align*}
Denote $\check{\bm \vartheta}_i = (\check\vartheta_i, \check a_i)$. To prove the statement for ${\bm \xi}_{i,\epsilon}'$, by the asymptotic behavior of the affine vortex at infinity, it suffices to bound the norm
\beqn
\| ( 1 - \gamma_i^\epsilon ) \check {\bm \vartheta}_i \|_{L_m^{1,p,\delta_p}}
\eeqn
where the norm is defined using the trivial connection. By Lemma \ref{lemma43}, one has
\beqn
\| \check \vartheta_i^H \|_{L_H^{1,p,\delta_0}({\mb A}_i \setminus B_R)} + \| \check \vartheta_i^G \|_{L_G^{1,p,\delta_0}({\mb A}_i \setminus B_R)} + \| \check{a}_i \|_{L_G^{1,p,\delta_0}({\mb A}_i \setminus B_R)} < \infty.
\eeqn
Notice that $(1-\gamma_i^\epsilon )$ is supported in ${\mb A}_i \setminus B_{\frac{1}{2b\sqrt\epsilon}}$; the derivative of $\gamma_i^\epsilon$ which is supported in the annular region whose area is of the scale $\epsilon^{-1}$, is bounded by a multiple of $\sqrt\epsilon$. Then one has
\begin{multline*}
\| ( 1- \gamma_i^\epsilon) \check {\bm \vartheta}_i \|_{L_m^{1,p, \delta_p}( {\mb A}_i \setminus B_{\frac{1}{2b\sqrt\epsilon}})} \\
\leq \| \check {\bm \vartheta}_i \|_{L_m^{1,p, \delta_p} ( {\mb A}_i \setminus B_{\frac{1}{2b \sqrt\epsilon}})} +  \| \nabla \gamma_i^\epsilon \|_{L^\infty} \| \check {\bm \vartheta}_i \|_{L^\infty} \big( {\rm Area} ({\rm Supp} \nabla \gamma_i^\epsilon ) \big)^{\frac{1}{p}}\\
\leq \| \check {\bm \vartheta}_i \|_{L_m^{1,p, \delta_p}({\mb A}_i \setminus B_{\frac{1}{2b \sqrt\epsilon}})} + C \epsilon^{\frac{1}{2} - \frac{1}{p}} \| \check{\bm \vartheta}_i \|_{L^\infty}.
\end{multline*}
Since $p > 2$ and $\delta_p < \delta_0$, we see that the right hand side above converges to zero as $\epsilon \to 0$. The estimate for ${\bm \xi}_{0, \epsilon}'$ is similar and omitted. 
\end{proof}

\begin{rem}\label{rem52}
There is a simplification one can take if there are constant affine vortex components. If an affine vortex component is constant, the stability condition requires that there are marked points on this component. We merely treat the marked points as combinatorial data and we do not need to glue this constant affine vortex component with the disk component. 
\end{rem}

\subsection{The weighted Sobolev norms}\label{subsection52}

We would like to define a Banach space of infinitesimal deformations of the approximate solutions. First we choose a weight function interpolating between weight functions on different components. Define 
\begin{align}\label{weight}
&\ \rho_\epsilon: {\mb H} \to [1, +\infty),\ &\ \rho_\epsilon(z) = \left\{ \begin{array}{cc}  \rho_i( z-z_{i,\epsilon}),\ &\ z \in \Sigma_{i, \epsilon},\\
                                            \displaystyle  \frac{1}{\epsilon} \rho_0 (\epsilon z),\ &\ z \in \Sigma_{0,\epsilon}.   \end{array} \right.
                                            \end{align}
Notice that this function is continuous and its value on the (half) circles $\partial \Sigma_{i,\epsilon}$ is equal to $\epsilon^{-1}$. Then for a section $f\in L^p_{\rm loc}({\mb H}, E)$ of an Euclidean bundle $E \to {\mb H}$, define
\beqn
\| f \|_{L_\epsilon^{p,\delta}} = \left( \int_{\mb H} | f (z)|^p (\rho_\epsilon(z))^{\frac{p\delta}{2}}ds dt \right)^{\frac{1}{p}}.
\eeqn
Define Banach spaces $\hat{\mc B}_\epsilon$ and $\hat{\mc E}_\epsilon$ as
\beq\label{eqn58}
\hat{\mc E}_\epsilon:= \big\{ {\bm \nu} = (\nu, \varsigma', \varsigma'') \in L^p_{\rm loc}( {\mb H}, u_\epsilon^* TV \oplus {\mf k} \oplus {\mf k})\ |\ \|{\bm \nu}\|_{L_\epsilon^{p,\delta_p}} < \infty \big\};
\eeq
\beq\label{eqn59}
\hat{\mc B}_\epsilon:= \big\{ {\bm \xi} = (\xi, \eta, \zeta) \in W^{1,p}_{\rm loc}({\mb H}, u_\epsilon^* TV \oplus {\mf k} \oplus {\mf k})_L \ |\ \| {\bm \xi} \|_{L_{m,\epsilon}^{1,p,\delta_p}} < \infty\big\},
\eeq
where for a general $\delta \in {\mb R}$, 
\beq\label{eqn510}
 \| {\bm \xi} \|_{L_{m,\epsilon}^{1,p,\delta}}:= \| \eta \|_{L_\epsilon^{p,\delta}} + \| \zeta \|_{L_\epsilon^{p,\delta}} + \| d\mu(u_\epsilon) (\xi) \|_{L_\epsilon^{p,\delta}} + \| d \mu(u_\epsilon) (J_V \xi )\|_{L_\epsilon^{p,\delta}} + \big\| \nabla^{a_\epsilon} {\bm \xi} \big\|_{L_\epsilon^{p,\delta}} + \| \xi \|_{L^\infty}. 
\eeq

We need the following uniform Sobolev estimate. 

\begin{lemma}\label{Sobolev}
There exist $C $ and $\epsilon_2$ such that for all $\epsilon\in (0, \epsilon_2)$ and ${\bm \xi} \in \hat{\mc B}_\epsilon$, one has 
\beqn
\| {\bm \xi} \|_{L^\infty} \leq C \| {\bm \xi} \|_{L_{m, \epsilon}^{1,p,\delta_p}}.
\eeqn
\end{lemma}

\begin{proof}
Given ${\bm \xi} = (\xi, \beta) \in \hat{\mc B}_\epsilon$, by the definition of the norm $\| {\bm \xi}\|_{L_{m,\epsilon}^{1,p,\delta}}$ (see \eqref{eqn510}) one has $\| \xi \|_{L^\infty} \leq \| {\bm \xi}\|_{L_{m, \epsilon}^{1,p,\delta_p}}$. On the other hand, notice that the norm on the components $\eta$ and $\zeta$ is a weighted Sobolev norm with weight function bigger than one. Hence by the usual Sobolev embedding $W^{1,p}\hookrightarrow C^0$ for $p>2$ in two dimensions, for some $C>0$ (which only depends on $p$) one has 
\beqn
\| \eta \|_{L^\infty} + \| \zeta \|_{L^\infty} \leq C \| {\bm \xi}\|_{L_{m,\epsilon}^{1,p,\delta_p}}.
\eeqn
This finishes the proof.
\end{proof}

\subsection{The implicit function theorem}\label{subsection53}

In this subsection we apply the implicit function theorem to find exact solutions close to the approximate ones. First we reformulate the problem using notations introduced above. We identify a small ball of the Banach space $\hat{\mc B}_\epsilon$ with a set of gauged maps near the approximate solution ${\bm v}_\epsilon$. For ${\bm \xi} = (\xi, \eta, \zeta) \in \hat{\mc B}_\epsilon$ with $\| {\bm \xi}\|_{L_{m,\epsilon}^{1,p,\delta_p}}$ sufficiently small, Lemma \ref{Sobolev} implies that $\| \xi\|_{L^\infty}$ is sufficiently small. Then using the exponential map $\exp$ associated to the metric $h_V$, one can identify ${\bm \xi}$ with a nearby gauged map 
\beqn
{\bm v}_\epsilon':= \exp_{{\bm v}_\epsilon} {\bm \xi}:= (u_\epsilon', \phi_\epsilon', \psi_\epsilon'):= (\exp_{u_\epsilon} \xi, \phi_\epsilon + \eta, \psi_\epsilon + \zeta).
\eeqn
By abuse of notations, we still use $\hat{\mc B}_\epsilon$ to denote the set of these nearby gauged maps. For each nearby gauged map ${\bm v}_\epsilon'$, one defines
\beqn
\hat{\mc E}_{{\bm v}_\epsilon'}:= L^{p,\delta_p}( {\mb H}, (u_\epsilon')^*TV \oplus {\mf k} \oplus {\mf k}).
\eeqn
These Sobolev spaces define a Banach vector bundle over the space of these nearby gauged maps, denoted (with abuse of notations) by $\hat{\mc E}_\epsilon$. Using the parallel transport with respect to the connection $\nabla$ (specified by Lemma \ref{connection}) along the shortest pointwise geodesic connection $u_\epsilon$ and $u_\epsilon'$, one trivializes the bundle $\hat{\mc E}_\epsilon$ so that each fibre is identified with the central fibre $\hat{\mc E}_{{\bm v}_\epsilon}$ (which is the space \eqref{eqn58}). The affine vortex equation over the curve ${\mc C}_{\epsilon, {\bm w}}'\cong {\mb H}$ (with respect to the domain dependent almost complex structure $J_{\epsilon, {\bm w}}$ and the volume form $\sigma_{\epsilon, {\bm w}}$) together with the Coulomb gauge-fixing condition relative to the central point ${\bm v}_\epsilon$ then defines a (nonlinear) map
\beqn
\hat{\mc F}_\epsilon: W_{\rm def} \times \hat{\mc B}_\epsilon \to \hat{\mc E}_\epsilon
\eeqn
from (an open set of) the Banach space $\hat{\mc B}_\epsilon$ to the Banach space $\hat{\mc E}_\epsilon$. For fixed $\epsilon$, $\hat{\mc F}_\epsilon$ is a smooth map. We would like to solve the equation 
\beqn
\hat{\mc F}_\epsilon({\bm w}, {\bm v}) = 0
\eeqn
for ${\bm w}$ near ${\bm w}^*$ and ${\bm v}$ near ${\bm v}_\epsilon$, for all sufficiently small $\epsilon$. 

The existence of the exact solution using the implicit function theorem relies on three major ingredients. We state these ingredients below and give their proofs later. First we need to estimate the failure of the approximate solution from being an exact solution. 

\begin{prop}\label{prop54}
There exist $\epsilon_1 > 0$, $C >0$, and $\tau>0$ such that for all $\epsilon \in (0, \epsilon_1)$
\beqn
\| \hat {\mc F}_\epsilon ( {\bm w}^*, {\bm v}_\epsilon) \|_{L_\epsilon^{p,\delta_p}} \leq C \epsilon^\tau.
\eeqn
\end{prop}
It is proved in Subsection \ref{subsection54}.

Next we need to estimate the variation of the linearized operators. Consider a deformation parameter $\uds {\bm w} \in W_{\rm def}$ very close to ${\bm w}^*$ and a small infinitesimal deformation 
\beqn
{\bm \rho} = (\rho, \upsilon) \in \hat {\mc B}_\epsilon ,\ {\rm where}\ \rho \in \Gammait(u_\epsilon^* TV)\ {\rm and} \ \upsilon = \varrho ds + \varsigma dt.
\eeqn
Consider the gauged map 
\beqn
\uds {\bm v}_\epsilon = (\uds u_\epsilon, \uds a_\epsilon):= \exp_{{\bm v}_\epsilon} {\bm \rho} = (\exp_{u_\epsilon}  \rho, a_\epsilon + \upsilon ).
\eeqn
The linearization of $\hat {\mc F}_\epsilon$ at $({\bm w}^*, {\bm v}_\epsilon)$ resp. $( \uds {\bm w}, \uds {\bm v}_\epsilon)$ is a linear operator
\beqn
{\mc D}_\epsilon: W_{\rm def} \oplus \hat{\mc B}_\epsilon \to \hat{\mc E}_\epsilon\ \hspace{2cm}\ {\rm resp.}\ \hspace{2cm}\ \uds{\hat {\mc D}}_\epsilon: W_{\rm def} \oplus \hat{\mc B}_\epsilon \to \hat{\mc E}_\epsilon. 
\eeqn

In Subsection \ref{subsection55} we prove the following proposition.

\begin{prop}\label{prop55}
There exist $\epsilon_2 > 0$ and $C > 0$ such that for all $\epsilon \in (0, \epsilon_2)$ and all ${\bm \rho} \in \hat{\mc B}_\epsilon$ with $\| {\bm \rho} \|_{L_{m, \epsilon}^{1,p,\delta_p}} \leq \epsilon_2$ and $\uds {\bm w} \in W_{\rm def}$ with $| \uds {\bm w} - {\bm w}^* | \leq \epsilon_2$, using the notations above, there holds
\begin{multline}\label{eqn511}
\big\| \uds{\hat{\mc D}}_\epsilon( {\bm w}, {\bm \xi}) - \hat{\mc D}_\epsilon ( {\bm w}, {\bm \xi}) \big\|_{L_\epsilon^{p,\delta_p}} \\
 \leq C \big( |  \uds {\bm w} - {\bm w}^* | +  \| {\bm \rho} \|_{L_{m,\epsilon}^{1,p,\delta_p}} \big) \big( |{\bm w}| + \| {\bm \xi} \|_{L_{m,\epsilon}^{1,p,\delta_p}} \big),\ \forall {\bm \xi} \in \hat{\mc B}_\epsilon,\ {\bm w} \in W_{\rm def}.
 \end{multline}
\end{prop}

Next one needs a right inverse of the linearized operator at the approximate solution. One needs an $\epsilon$-independent upper bound on the norm of this right inverse. 

\begin{prop}\label{prop56}
There exist $\epsilon_3>0$, $C>0$, and, for each $\epsilon \in (0, \epsilon_3)$, a bounded right inverse 
\beqn
\hat{\mc Q}_\epsilon: \hat{\mc E}_\epsilon \to  W_{\rm def} \oplus \hat{\mc B}_\epsilon
\eeqn
to the operator $\hat{\mc D}_\epsilon$ such that $\| \hat{\mc Q}_\epsilon \| \leq C$. 
\end{prop}

The construction of $\hat {\mc Q}_\epsilon$ and the proof of this proposition are given in Subsection \ref{subsection56}. Again there is no need to consider the extension to families of approximate solutions. 

\subsubsection{The gluing map}	

Now we are ready to apply the implicit function theorem. Let us first cite a precise version of it. 

\begin{prop}\cite[Proposition A.3.4]{McDuff_Salamon_2004}\label{prop57}
Let $X$, $Y$ be Banach spaces, $U \subset X$ be an open set and $f: U \to Y$ be a continuously differentiable map. Let $x_\circ \in U$ be such that $df(x_\circ): X \to Y$ is surjective and has a bounded right inverse $Q:Y \to X$. Assume there are constants $a, c>0$ such that 
\beq\label{eqn512}
\| Q \| \leq c;
\eeq
\beq\label{eqn513}
B_a (x_\circ ) \subset U\ {\rm and}\ x \in B_a (x_0) \Longrightarrow \| df(x) - df(x_\circ) \| \leq \frac{1}{2c}.
\eeq
Suppose $x' \in X$ satisfies
\begin{align}\label{eqn514}
&\ \| f(x')\| < \frac{a}{4c},\  &\ \| x' - x_\circ \| < \frac{a}{8},
\end{align}
then there exists a unique $x \in X$ satisfying
\beq\label{eqn515}
f(x) = 0,\ x- x' \in {\rm Im}Q,\ \| x- x_\circ \| \leq a.
\eeq
Moreover there holds
\beqn
\| x- x' \| \leq 2c \| f(x')\|.
\eeqn
\end{prop}

We can apply the implicit function theorem as follows. For each fixed small $\epsilon$, set
\begin{align*}
&\ X:= W_{\rm def} \oplus \hat{\mc B}_\epsilon,\ &\ Y:= \hat{\mc E}_\epsilon,
\end{align*}
set $f: X \to Y$ to be the map $\hat{\mc F}_\epsilon$, and set $x_\circ = x'$ to be the point $({\bm w}^*, {\bm v}_\epsilon)$. Then Proposition \ref{prop56} implies that \eqref{eqn512} holds for a certain $c$ which is independent of sufficiently small $\epsilon$. Further Proposition \ref{prop55} implies \eqref{eqn513} holds for a certain $a$, and Proposition \ref{prop54} implies that, when $\epsilon$ is sufficiently small, \eqref{eqn514} is satisfied. Therefore it follows from Proposition \ref{prop57} that exists a unique $x$ satisfying \eqref{eqn515}, which we denote by 
\beqn
({\bm w}_\epsilon, \tilde {\bm v}_\epsilon) \in W_{\rm def} \times \hat{\mc B}_\epsilon. 
\eeqn
It defines a point in $\ov{\mc M}_{l, l^+}(V, L_V)$ represented by $({\mc C}_{\epsilon, {\bm w}_\epsilon}', \tilde {\bm v}_\epsilon)$.

\begin{defn}[Gluing map]\label{gluingmap}
For $\tau>0$ sufficiently small, define
\begin{align*}
&\ {\it Glue}: [0, \tau) \to \ov{\mc M}_{l, l^+}(V, L_V)\ &\ \epsilon \mapsto \left\{ \begin{array}{cc} {\bm p},\ &\ \epsilon = 0,\\
{[} {\mc C}_{\epsilon, {\bm w}_\epsilon}', \tilde {\bm v}_\epsilon {]},\ &\ \epsilon > 0. \end{array} \right.
\end{align*}
\end{defn}

It is not difficult to show that the gluing map is continuous. Indeed, the Banach space $\hat{\mc B}_\epsilon$ resp. $\hat{\mc E}_\epsilon$ for different positive $\epsilon$'s are isomorphic as topological vector spaces. The $\epsilon$-dependent norms are comparable to each other. Moreover, the family of approximate solutions form a continuous curve in the topological vector space. Hence the path $\epsilon \to \tilde {\bm v}_\epsilon$ can be proved to be continuous with respect to the $W^{1,p}_{\rm loc}$-topology, and hence the $C^\infty_{\rm loc}$-topology. Moreover, it is straightforward to see that as $\epsilon \to 0$, the path $({\mc C}_{\epsilon, {\bm w}_\epsilon}', \tilde {\bm v}_\epsilon)$ converges to the singular configuration we start with. Hence the gluing map is continuous. 

It remains to show that the gluing map is a local homeomorphism to finish the proof of the main theorem. This is done in Section \ref{section6}.

\subsection{Proof of Proposition \ref{prop54}}\label{subsection54}

\begin{proof}[Proof of Proposition \ref{prop54}]
We denote the three components of $\hat {\mc F}_\epsilon$ by $\hat {\mc F}_1$, $\hat {\mc F}_2$ and $\hat {\mc F}_3$ respectively, where only $\hat {\mc F}_1$ depends on the perturbation term. Since the Coulomb gauge fixing condition is with respect to the approximate solution ${\bm v}_\epsilon$ itself, $\hat {\mc F}_2 ( {\bm w}^*, {\bm v}_\epsilon ) = 0$ automatically. The rest of the proof is divided into the following two major steps. 

\vspace{0.2cm}

\noindent {\bf Step 1. Estimate $\hat {\mc F}_1( 0, {\bm v}_\epsilon)$.} Abbreviate $J_{\epsilon, {\bm w}^*}$ by $J_\epsilon$. We estimate in different regions of the domain as follows. 

\vspace{0.1cm}

\noindent {\it Step 1(a).} Inside each $\Sigma_{i,\epsilon}^\ll$, ${\bm v}_\epsilon$ agrees with ${\bm v}_i$ after a translation. Hence
\beqn
\hat {\mc F}_1( {\bm w}^*, {\bm v}_\epsilon) = \partial_s u_i + {\mc X}_{\phi_i} + J_\epsilon( \partial_t u_i + {\mc X}_{\psi_i}) = (J_\epsilon - J_0) (\partial_t u_i + {\mc X}_{\psi_i}).
\eeqn
Here we used the fact that ${\bm v}_i$ is an affine vortex with respect to $J_0$. Then by the smooth dependence of the almost complex structure on the parameter $\epsilon$ and the fact that $J_\epsilon - J_0$ is supported over a bounded subset of ${\mb A}_i$, for some $C>0$ independent of $\epsilon$, one has
\beqn
\big\| \hat {\mc F}_1({\bm w}^*, {\bm v}_\epsilon) \big\|_{L_\epsilon^{p,\delta_p}(\Sigma_{i,\epsilon}^\ll)} \leq \| J_\epsilon - J_0 \|_{L^\infty} \| \partial_t u_i + {\mc X}_{\psi_i} \|_{L^{p,\delta_p}}  \leq C \epsilon.
\eeqn

\vspace{0.1cm}

\noindent {\it Step 1(b).} Over $\Sigma_{0,\epsilon}^\ll$, ${\bm v}_\epsilon$ agrees with the rescaled ${\bm v}_{0, \epsilon}$. Hence one has
\beqn
\begin{split}
\hat{\mc F}_1( {\bm w}^*, {\bm v}_\epsilon) = &\ \partial_s u_{0, \epsilon} + {\mc X}_{\phi_{0, \epsilon}} + J_\epsilon( \partial_t u_{0, \epsilon} + {\mc X}_{\psi_{0, \epsilon}})  \\
= &\ (J_\epsilon - J_0) (\partial_t u_{0, \epsilon} + {\mc X}_{\psi_{0, \epsilon}} )\\
= &\ \epsilon s_\epsilon^* ( J_\epsilon - J_0) (P_H( \partial_t u_0)).
\end{split}
\eeqn
Then by the smooth dependence of $J_\epsilon$ on $\epsilon$ and Lemma \ref{lemma44}, one has
\beqn
\| \hat{\mc F}_1( {\bm w}^*, {\bm v}_\epsilon ) \|_{L_\epsilon^{p,\delta_p}(\Sigma_{0,\epsilon}^\ll)} \leq \epsilon \| J_\epsilon - J_0 \|_{L^\infty} \| s_\epsilon^* (P_H(\partial_t u_0)) \|_{L_\epsilon^{p,\delta_p}(\Sigma_{0,\epsilon}^\ll)} \leq C \epsilon \| P_H( \partial_t u_0) \|_{L^{p, \delta_p}} \leq C \epsilon.
\eeqn

\vspace{0.1cm}

\noindent {\it Step 1(c).} For the estimate over the inner part of the neck region $A_{i,\epsilon}^-:= \Sigma_{i,\epsilon} \setminus \Sigma_{i,\epsilon}^\ll$, we first prove an estimate regarding the affine vortex component ${\bm v}_i$. By the asymptotic decay of the energy density (see \cite[Theorem 1.3]{Ziltener_Decay} and \cite[Proposition A.4]{Venugopalan_Xu}), for all (small) $\alpha>0$ there exists $C(\alpha)>0$ such that 
\beqn
| \partial_s \exp_{x_i} \check\vartheta_i + {\mc X}_{\check \phi_i}(\exp_{x_i} \check\vartheta_i) | \leq C(\alpha) |z-z_{i,\epsilon}|^{-2 +\alpha}.
\eeqn
Fix $\alpha < \frac{2}{p}$. Then there is a constant $C>0$ such that for all $r> 0$, 
\beq\label{www}
\left\| |z|^{-2+\alpha} \right\|_{L^{p, \delta_p}({\mb A}_i \setminus B_r)} =  \left( \int_{{\mb A}_i \setminus B_r} |z|^{-(2-\alpha)p} |z|^{2p-4} ds dt \right)^{\frac{1}{p}} \leq C r^{\alpha p -2} < \infty.
\eeq
Therefore, one has
\beq\label{xxx}
\| \partial_s \exp_{x_i} \check\vartheta_i + {\mc X}_{\check \phi_i}( \exp_{x_i} \check\vartheta_i) \|_{L_\epsilon^{p,\delta_p}(A_{i,\epsilon}^-)} \leq C \left\| |z|^{-2+\alpha} \right\|_{L^{p, \delta_p}({\mb A}_i \setminus B_{\frac{1}{2b\sqrt\epsilon}} )} \leq C (\sqrt\epsilon)^{2- \alpha p}.
\eeq
Moreover, by using the derivatives of the exponential map (see notations in \eqref{eqna1}), one has
\beq\label{yyy}
\partial_s \exp_{x_i} \check\vartheta_i + {\mc X}_{\check \phi_i}(\exp_{x_i} \check\vartheta_i) = E_2 ( \nabla_s \check\vartheta_i + \nabla_{\check\vartheta_i} {\mc X}_{\check \phi_i})+ E_1 {\mc X}_{\check \phi_i}(x_i);
\eeq
The asymptotic behavior of ${\bm v}_i$ (see Lemma \ref{lemma43} and \eqref{eqn45}) implies that 
\begin{multline}\label{zzz}
 \| E_1 {\mc X}_{\check \phi_i}(x_i) \|_{L_\epsilon^{p,\delta_p}(A_{i,\epsilon}^-)} \\
  \leq \Big( \int_{A_{i, \epsilon}^-}   | \check a_i |^p (\rho_\epsilon(z))^{p-2} ds dt \Big)^{\frac{1}{p}} \leq  \Big( \int_{A_{i,\epsilon}^-}  | \check a_i|^p (\rho_i(z))^{p-2} ds dt \Big)^{\frac{1}{p}} \\
 \leq \| \check a_i \|_{L^{p, \delta_p}( {\mb A}_i \setminus B_{\frac{1}{2b\sqrt\epsilon}} )} \leq C (\sqrt\epsilon)^{\delta_0 - \delta_p} \| \check a_i \|_{L^{p,\delta}}.
\end{multline}
It follows from \eqref{xxx}, \eqref{yyy}, and \eqref{zzz} that there is a constant $C>0$ such that 
\beq\label{eqn520}
\| \nabla_s \check\vartheta_i +  \nabla_{\check\vartheta_i} {\mc X}_{\check \phi_i}\|_{L_\epsilon^{p,\delta_p}(A_{i,\epsilon}^-)}   \leq C \epsilon^\tau,\ {\rm where}\ \tau = \min \left\{ 1- \frac{ \alpha p }{2},  \frac{\delta_0 - \delta_p}{2} \right\} > 0.
\eeq

Now we estimate $\hat  {\mc F}_1({\bm w}^*, {\bm v}_\epsilon)$ over the region $A_{i,\epsilon}^-$. By the construction of ${\bm v}_\epsilon$, one has
\begin{align*}
&\ u_\epsilon(z) = e^{- \lambda_i \theta_i} \cdot \exp_{x_i}(\gamma_i^\epsilon \check\vartheta_i (z)),\ &\ a_\epsilon = e^{-\lambda_i \theta_i} \cdot (  \gamma_i^\epsilon \check a_i).
\end{align*}
Therefore, by using the linear maps $E_1$ and $E_2$, one has
\begin{multline*}
e^{\lambda_i \theta_i} \big( \partial_s u_\epsilon + {\mc X}_{\phi_\epsilon}(u_\epsilon) \big) = \partial_s (\exp_{x_i} \gamma_i^\epsilon \check \vartheta_i ) + {\mc X}_{\gamma_i^\epsilon  \check \phi_i}(\exp_{x_i} \gamma_i^\epsilon \check \vartheta_i ) \\
= E_2 \big(  \gamma_i^\epsilon  \nabla_s \check \vartheta_i + (\gamma_i^\epsilon)^2 \nabla_{\check \vartheta_i} {\mc X}_{\check \phi_i} + ( \partial_s \gamma_i^\epsilon ) \check \vartheta_i \big) + E_1( \gamma_i^\epsilon {\mc X}_{\check\phi_i}(x_i)).
\end{multline*}
Hence 
\begin{multline*}
\| \partial_s u_\epsilon + {\mc X}_{\phi_\epsilon}(u_\epsilon) \|_{L_\epsilon^{p,\delta_p}(A_{i,\epsilon}^-)} \\
\leq \| \gamma_i^\epsilon ( \nabla_s \check \vartheta_i + {\mc X}_{\check \vartheta_i} {\mc X}_{\check \phi_i}) \|_{L_\epsilon^{p,\delta_p}(A_{i,\epsilon}^-)} + \| ( (\gamma_i^\epsilon)^2 - \gamma_i^\epsilon ) \nabla_{\check \vartheta_i} {\mc X}_{\check \phi_i} \|_{L_\epsilon^{p,\delta_p}(A_{i,\epsilon}^-)} \\
+ \| E_1(\gamma_i^\epsilon {\mc X}_{\check \phi_i}(x_i) )\|_{ L_\epsilon^{p,\delta_p}(A_{i,\epsilon}^-)} + \| (\partial_s \gamma_i^\epsilon ) \check \vartheta_i \|_{L_\epsilon^{p,\delta_p}(A_{i,\epsilon}^-)}.
\end{multline*}
The first term can be bounded using \eqref{eqn520}, the second and the third term can be bounded by using \eqref{zzz}. For the fourth term, using \eqref{exponential} and \eqref{eqn52}, one has 
\begin{multline*}
\| (\partial_s \gamma_i^\epsilon ) \check \vartheta_i \|_{L_\epsilon^{p,\delta_p}(A_{i,\epsilon}^-)}  \leq C \| \nabla \gamma_i^\epsilon \|_{L^\infty} \| \check \vartheta_i\|_{L^\infty( A_{i,\epsilon}^-)} \|\rho_i(z)^{\frac{\delta_p}{2}} \|_{L^\infty(A_{\epsilon;2}^{i,-})} \left( {\rm Area} (A_{\epsilon;2}^{i,-})\right)^{\frac{1}{p}}\\
\leq C \sqrt\epsilon \cdot (\sqrt\epsilon)^{\delta_0 - \frac{\delta_p}{2}} \cdot \epsilon^{-\frac{\delta_p}{2}} \cdot  \epsilon^{-\frac{1}{p}}\leq C \epsilon^{\frac{1}{2} (\delta_0 - \delta_p)}.
\end{multline*}
Therefore one has
\beqn
\| \partial_s u_\epsilon + {\mc X}_{\phi_\epsilon}(u_\epsilon) \|_{L_\epsilon^{p,\delta_p}(A_{i,\epsilon}^-)} \leq C \epsilon^\tau.
\eeqn
Similarly one has 
\beqn
\| \partial_t u_\epsilon + {\mc X}_{\psi_\epsilon}(u_\epsilon) \|_{L_\epsilon^{p,\delta_p}(A_{i,\epsilon}^-)}  \leq C \epsilon^\tau.
\eeqn
So 
\beqn
\| \hat {\mc F}_1({\bm w}^*, {\bm v}_\epsilon) \|_{L_\epsilon^{p,\delta_p}(A_{i,\epsilon}^-)} = \| \partial_s u_\epsilon + {\mc X}_{\phi_\epsilon}(u_\epsilon) + J_V( \partial_t u_\epsilon + {\mc X}_{\psi_\epsilon}(u_\epsilon)) \|_{L_\epsilon^{p,\delta_p}(A_{i,\epsilon}^-)} \leq C \epsilon^\tau.
\eeqn

\vspace{0.1cm}

\noindent {\it Step 1(d).} In the outer part of the neck region $A_{i,\epsilon}^+:= \Sigma_{i,\epsilon}^\gg \setminus \Sigma_{i,\epsilon}$ we can derive similar estimate by comparing with the rescaled disk ${\bm v}_{0, \epsilon}$. We omit the details of the tedious verification.

\vspace{0.1cm}

\noindent In summary, for a certain $\tau>0$ and a certain $C>0$, for $\epsilon$ sufficiently small, one has
\beq\label{eqn521}
\| \hat{\mc F}_1( {\bm w}^*, {\bm v}_\epsilon) \|_{L_\epsilon^{p,\delta_p}({\mb H})} \leq C \epsilon^\tau.
\eeq

\vspace{0.2cm}

\noindent {\bf Step 2. Estimate $\hat {\mc F}_2( {\bm w}^*, {\bm v}_\epsilon)$.}

\vspace{0.1cm}

\noindent {\it Step 2(a).} Over the region $\Sigma_{i,\epsilon}^\ll$, ${\bm v}_\epsilon$ agrees with ${\bm v}_i$ after a translation. Then 
\beqn
\hat {\mc F}_2({\bm w}^*, {\bm v}_\epsilon)|_{\Sigma_{i,\epsilon}^\ll} = 0.
\eeqn 

\vspace{0.1cm}

\noindent {\it Step 2(b).} Over the interior part of the neck region $A_{i,\epsilon}^-$, one has
\begin{multline*}
e^{\lambda_i \theta_i} \big( \hat {\mc F}_3({\bm w}^*, {\bm v}_\epsilon) \big) = \partial_s \check \psi_\epsilon - \partial_t \check \phi_\epsilon + [\check \phi_\epsilon, \check \psi_\epsilon] + \mu(\check u_\epsilon)\\
= \partial_s (\gamma_i^\epsilon \check \psi_i ) - \partial_t( \gamma_i^\epsilon \check \phi_i) + ( \gamma_i^\epsilon)^2 [\check \phi_i , \check \psi_i] + \mu( \exp_{x_i} \gamma_i^\epsilon \check \vartheta_i)\\
= \gamma_i^\epsilon ( \partial_s \check \psi_i- \partial_t \check \phi_i + [\check \phi_i, \check \psi_i] ) + ((\gamma_i^\epsilon)^2 - \gamma_i^\epsilon ) [\check \phi_i, \check \psi_i] + \mu(\exp_{x_i} \gamma_i^\epsilon \check \vartheta_i).
\end{multline*}
Fix $\alpha < \frac{2}{p}$. Then by the exponential decay of the energy density (see \cite[Theorem 1.3]{Ziltener_Decay} and \cite[Proposition A.4]{Venugopalan_Xu}) there holds
\beqn
| \partial_s \check \psi_i - \partial_t \check \phi_i + [\check \phi_i, \check \psi_i]| + |\mu(u_i)| \leq C |z-z_{i,\epsilon}|^{-2+\alpha}.
\eeqn
On the other hand, since $\mu^{-1}(0)$ is totally geodesic with respect to $h_V$, one has
\beqn
|\mu( \exp_{x_i}\gamma_i^\epsilon \check \vartheta_i)| \leq C |d\mu(x_i) \check\vartheta_i| \leq C |\mu(u_i)| \leq C |z-z_{i,\epsilon}|^{-2+\alpha}.
\eeqn
Using the \eqref{www} one obtains that 
\beq\label{eqn522}
\| \gamma_i^\epsilon ( \partial_s \check \psi_i- \partial_t \check \phi_i + [\check \phi_i, \check \psi_i] ) \|_{L_\epsilon^{p,\delta_p}} + \| \mu(\exp_{x_i} \gamma_i^\epsilon \check \vartheta_i ) \|_{L_\epsilon^{p,\delta_p}} \leq C \epsilon^{\frac{\alpha p}{2} -1}.
\eeq
On the other hand, by \eqref{eqn45} and Sobolev embedding, one has 
\beqn
|\check \phi_i | + |\check \psi_i| \leq C |z|^{-\delta_0} \Longrightarrow  \|\check \phi_i \|_{L^\infty(A_{i,\epsilon}^-)}  + \|\check \psi_i \|_{L^\infty( A_{i, \epsilon}^-)} \leq C \epsilon^{\frac{\delta_0}{2}}. 
\eeqn
Then
\begin{multline*}
\| ((\gamma_i^\epsilon)^2 - \gamma_i^\epsilon ) [\check \phi_i, \check \psi_i] \|_{L_\epsilon^{p,\delta_p}(A_{i,\epsilon}^-)} \\
 \leq C \| \check \phi_i\|_{L^\infty (A_{i,\epsilon}^-)} \| \check \psi_i \|_{L^\infty (A_{i,\epsilon}^-)} \| \rho_i^{\frac{\delta_p}{2}} \|_{L^\infty( A_{i,\epsilon}^-)} \left( {\rm Area} (A_{i,\epsilon}^-) \right)^{\frac{1}{p}}\\
 \leq C \epsilon^{\frac{\delta_0}{2}} \cdot \epsilon^{\frac{\delta_0}{2}} \cdot \epsilon^{-\frac{\delta_p}{2}} \cdot \epsilon^{-\frac{1}{p}} = C \epsilon^{\delta_0 - ( 1- \frac{1}{p})}.
\end{multline*}
Notice that $\delta_0 >  1 - \frac{1}{p}$ (see \eqref{delta}). By combining above with \eqref{eqn522} one can see that for an appropriate $\tau>0$ and $C>0$ one has 
\beqn
\| \hat {\mc F}_2( {\bm w}^*, {\bm v}_\epsilon) \|_{L_\epsilon^{p,\delta_p}(A_{i,\epsilon}^-)}  \leq C \epsilon^\tau.
\eeqn

\vspace{0.1cm}

\noindent {\it Step 2(c).} Over the complement of the union of the above two regions, one has
\begin{multline*}
e^{\lambda_i \theta_i} \hat {\mc F}_2({\bm w}^*, {\bm v}_\epsilon) = \partial_s (\gamma_0^\epsilon \check \psi_{0, \epsilon}) - \partial_t (\gamma_0^\epsilon \check \phi_{0, \epsilon} ) + (\gamma_0^\epsilon)^2 [\check \phi_{0,\epsilon}, \check \psi_{0, \epsilon}] + \mu \big( \exp_{x_i} (\gamma_0^\epsilon \check\vartheta_{0,\epsilon} ) \big)\\
= \gamma_0^\epsilon \big( \partial_s \check \psi_{0, \epsilon} - \partial_t \check \phi_{0, \epsilon} + [\check \phi_{0, \epsilon}, \check \psi_{0, \epsilon} ] \big) + ( (\gamma_0^\epsilon)^2 - \gamma_0^\epsilon ) [\check \phi_{0, \epsilon}, \check \psi_{0, \epsilon} ].
\end{multline*}
For the last equality, we used the fact that $\mu^{-1}(0)$ is totally geodesic with respect to $h_V$ (see Lemma \ref{metric}). Notice that 
\beqn
\partial_s \check \psi_{0, \epsilon} - \partial_t \check \phi_{0, \epsilon} + [\check \phi_{0, \epsilon}, \check \psi_{0, \epsilon}] = \epsilon^2 s_\epsilon^* ( \partial_s \psi_0 - \partial_t \phi_0 + [\phi_0, \psi_0]).
\eeqn
Then one has
\begin{multline}\label{eqn523}
\| \gamma_0^\epsilon  (\partial_s \check \psi_{0, \epsilon} - \partial_t \check \phi_{0, \epsilon} + [\check \phi_{0, \epsilon}, \check \psi_{0, \epsilon}]) \|_{L_\epsilon^{p,\delta_p}(\Sigma_{0,\epsilon})} \\
\leq \| \partial_s \check \psi_{0, \epsilon} - \partial_t \check \phi_{0, \epsilon} + [\check \phi_{0, \epsilon}, \check \psi_{0, \epsilon}] \|_{L_\epsilon^{p,\delta_p}({\mb H})} \\
= \epsilon \| \partial_s \psi_0 - \partial_t \phi_0 + [\phi_0, \psi_0] \|_{L^{p,\delta_p}({\mb H})} \leq C \epsilon. 
\end{multline}
Here the norm in the last line is finite because $\phi_0$, $\psi_0$ are induced from the holomorphic disk. On the other hand, since the holomorphic disk is smooth at the marked point $z_i$, one has 
\beqn
\| \check \phi_{0, \epsilon}\|_{L^\infty(\Sigma_{0,\epsilon})} + \| \check \psi_{0, \epsilon}\|_{L^\infty ( \Sigma_{0,\epsilon})} \leq C \epsilon.
\eeqn
Then one has
\begin{multline*}
\| ((\gamma_0^\epsilon)^2 - \gamma_0^\epsilon) [ \check \phi_{0, \epsilon}, \check \psi_{0, \epsilon} ] \|_{L_\epsilon^{p,\delta_p}(\Sigma_{0,\epsilon})} \\
\leq C \| \check \phi_{0, \epsilon}\|_{L^\infty} \| \check \psi_{0, \epsilon}\|_{L^\infty} \| \rho_\epsilon^{\frac{\delta_p}{2}} \|_{L^\infty({\rm supp} \gamma_0^\epsilon)} \left( {\rm Area}({\rm supp} \gamma_0^\epsilon) \right)^{\frac{1}{p}}\\
\leq C \epsilon \cdot \epsilon \cdot \epsilon^{-\frac{\delta_p}{2}} \cdot \epsilon^{-\frac{1}{p}} = C \epsilon^{1 + \frac{1}{p}}.
\end{multline*}
Combining with \eqref{eqn523} one has that 
\beqn
\| \hat {\mc F}_2({\bm w}^*, {\bm v}_\epsilon) \|_{L_\epsilon^{p,\delta_p}(\Sigma_{0,\epsilon})} \leq C \epsilon.
\eeqn
By combining Step 2(a)---2(c) one can see that there exist $\tau>0$ and $C>0$ such that for $\epsilon$ sufficiently small, there holds
\beq\label{eqn524}
\| \hat {\mc F}_2({\bm w}^*, {\bm v}_\epsilon ) \|_{L_\epsilon^{p,\delta_p}({\mb H})} \leq C \epsilon^\tau.
\eeq

\vspace{0.2cm}

\noindent Proposition \ref{prop54} then follows from \eqref{eqn521} and \eqref{eqn524}.
\end{proof}

\subsection{Proof of Proposition \ref{prop55}}\label{subsection55}

The proof is by direct calculation. Since the variation of the deformation parameter has very simple effect on the augmented linearization, we assume that ${\bm w} = {\bm w}^*$. We also only prove \eqref{eqn511} for ${\bm w} = 0$. Abbreviate $J_{\epsilon, {\bm w}^*}$ by $J$ and $\delta_p$ by $\delta$. 

The proof is divided into the following two major steps.

\vspace{0.2cm}

\noindent {\bf Step 1.} Consider an intermediate gauged map
\beqn
\uds{\bm v}_\epsilon' = (u_\epsilon, \uds a_\epsilon)
\eeqn
with linearized operator 
\beqn
\uds{\hat D}_\epsilon':= \hat D_{\uds {\bm v}_\epsilon'}: \hat{\mc B}_\epsilon \to \hat{\mc E}_{\uds {\bm v}_\epsilon'} \cong \hat{\mc E}_\epsilon|_{{\bm v}_\epsilon}.
\eeqn
We first compare $\uds{\hat D}_\epsilon'$ and $\hat D_\epsilon$, whose domains and codomains are canonically identified without using parallel transports. For a gauged map ${\bm v}= (u, \phi, \psi)$, a domain-dependent almost complex structure $J$, and a section $\xi \in \Gammait(u^* TV)$, introduce
\beqn
I({\bm v})(\xi):= \nabla_s \xi + \nabla_\xi {\mc X}_\phi + J( \nabla_t \xi + \nabla_\xi {\mc X}_\psi) + (\nabla_\xi J) (\partial_t u +{\mc X}_\psi).
\eeqn
Then for an infinitesimal deformation ${\bm \xi} = (\xi, \eta, \zeta)$ of the gauged map and an infinitesimal deformation $t$ of the deformation parameter, one has
\beqn
\hat D_\epsilon(\xi, \eta, \zeta) = \left[ \begin{array}{c} I({\bm v}_\epsilon) (\xi) + {\mc X}_{\eta} + J {\mc X}_\eta \\
\nabla_s^{a_\epsilon} \eta  + \nabla_t^{a_\epsilon} \zeta + d\mu(u_\epsilon) (J_V \xi) \\
\nabla_s^{a_\epsilon} \zeta - \nabla_t^{a_\epsilon} \eta + d\mu(u_\epsilon) (\xi)     \end{array} \right],
\eeqn
and
\beqn
\uds{\hat D}_\epsilon' (t, \xi, \eta, \zeta) = \left[ \begin{array}{c} I(\uds{\bm v}_\epsilon' ) (\xi) +  {\mc X}_{\eta} + J  {\mc X}_\eta \\
\nabla_s^{\uds a_\epsilon} \eta  + \nabla_t^{\uds a_\epsilon} \zeta + d\mu(u_\epsilon) (J_V \xi) \\
\nabla_s^{\uds a_\epsilon} \zeta - \nabla_t^{\uds a_\epsilon} \eta + d\mu(u_\epsilon) (\xi)     \end{array} \right].
\eeqn
Then by the definition of norms and the uniform Sobolev embedding (Lemma \ref{Sobolev}),
\beqn
\| I( \uds {\bm v}_\epsilon')(\xi) - I( {\bm v}_\epsilon)(\xi) \|_{L_\epsilon^{p,\delta}}\leq C \| |\xi| |\upsilon | \|_{L_\epsilon^{p,\delta}} \leq C \| {\bm \xi} \|_{L^\infty} \| \upsilon \|_{L_\epsilon^{p,\delta}} \leq C \| {\bm \xi} \|_{L_{m,\epsilon}^{1,p,\delta}} \| {\bm \rho} \|_{L_{m,\epsilon}^{1,p,\delta}}.
\eeqn
The other terms of the difference $\hat D_\epsilon(\xi, \eta, \zeta) - \uds{{\hat D}}_\epsilon'(\xi, \eta, \zeta)$ can be estimated similarly. Therefore one has (resetting the value of $C$)
\beqn
\| \uds{\hat D}_\epsilon'( {\bm \xi}) - \hat D_\epsilon ( {\bm \xi}) \|_{L_\epsilon^{p,\delta}} \leq C \| {\bm \xi}\|_{L_{m,\epsilon}^{1,p,\delta}}  \| {\bm \rho} \|_{L_{m,\epsilon}^{1,p,\delta}}.\eeqn

\vspace{0.1cm}

\noindent {\bf Step 2.} Now we compare the operator $\uds{ \hat D}_\epsilon'$ (which is the augmented linearization at $\uds {\bm v}_\epsilon'=(u_\epsilon, \uds a_\epsilon)$) and the operator $\uds{\hat D}_\epsilon$ (which is the linearization at $\uds {\bm v}_\epsilon = (\uds u_\epsilon, \uds a_\epsilon)$). We first consider the case when ${\bm \xi} = (0, \eta, \zeta)$. In this case, one has
\beqn
(\uds{\hat D}_\epsilon - \uds{\hat D}_\epsilon')(0, \eta, \zeta) = \left[ \begin{array}{c} {\it pl}_\epsilon^{-1} \big( {\mc X}_\eta( \uds u_\epsilon)  + J {\mc X}_\zeta (\uds u_\epsilon) \big) - \big( {\mc X}_\eta(u_\epsilon) + J {\mc X}_\zeta(u_\epsilon) \big) \\
  0 \\ 0  \end{array} \right]
\eeqn
Hence by Lemma \ref{Sobolev}, one obtains
\beqn
\| ( \uds {\hat D}_\epsilon - \uds{\hat D}_\epsilon') (0, \eta, \zeta) \|_{L_\epsilon^{p,\delta}}  \leq C \| (0,\eta, \zeta)\|_{L_\epsilon^{p,\delta}} \| \nu  \|_{L^\infty} \leq C' \| {\bm \xi}\|_{L_{m,\epsilon}^{1,p,\delta}} \| {\bm \nu} \|_{L_{m,\epsilon}^{1,p,\delta}}.
\eeqn
Second, consider the situation when ${\bm \xi} = (\xi, 0, 0)$. Denote $\uds \xi = E(u_\epsilon, \uds u_\epsilon)(\xi)$. Then 
\beq\label{eqn525}
( \uds{\hat D}_\epsilon - \uds{\hat D}_\epsilon' )(\xi, 0, 0) = \left[ \begin{array}{c}  {\it pl}_\epsilon^{-1} \big( I( \uds{\bm v}_\epsilon ) (\uds \xi)  \big)  - I( \uds{\bm v}_\epsilon') (\xi) \\
0 \\
d \mu( \uds u_\epsilon ) ( \uds \xi ) - d\mu(u_\epsilon)  (\xi)
\end{array} \right].
\eeq
For the last entry, we claim that for some $C$ independent of $u_\epsilon$, $\uds u_\epsilon$, and $\xi$ such that
\beq\label{special}
|d\mu(\uds u_\epsilon)( \uds \xi) - d\mu(u_\epsilon)(\xi)|\leq C \left( |\mu(u_\epsilon)| |\rho | |\xi| + |d\mu(u_\epsilon)( \rho )||\xi| + |d\mu(u_\epsilon)(\xi) ||\rho | \right).
\eeq
Indeed, if $\rho = 0$, then $d\mu(\uds u_\epsilon)(\uds \xi) - d\mu(u_\epsilon) (\xi) = 0$; if $\mu(u_\epsilon) = 0$, then since $\mu^{-1}(0)$ is totally geodesic, one has
\beqn
\mu(u_\epsilon) = 0 \Longrightarrow |d\mu(\uds u_\epsilon)( \uds \xi) - d\mu(u_\epsilon)(\xi)|\leq C \left( |d\mu(u_\epsilon)(\rho )||\xi| + |d\mu(u_\epsilon)(\xi) ||\rho | \right).
\eeqn
Then \eqref{special} follows. Therefore, since $\mu(u_\epsilon) = 0$ over $\Sigma_{0,\epsilon}$, one has
\begin{multline}\label{eqn527}
\| d\mu(\uds u_\epsilon)(\uds \xi) - d\mu(u_\epsilon)(\xi) \|_{L_\epsilon^{p,\delta}(\Sigma_{0,\epsilon})}\\ \leq C \left( \| \xi \|_{L^\infty} \| d\mu(u_\epsilon) (\nu_\epsilon) \|_{L_\epsilon^{p,\delta} } + \|\nu_\epsilon \|_{L^\infty} \| d\mu(u_\epsilon) (\xi) \|_{L_\epsilon^{p,\delta} } \right)\leq C \| {\bm \xi} \|_{L_{m,\epsilon}^{1,p,\delta}} \| {\bm \nu} \|_{L_{m,\epsilon}^{1,p,\delta}}.
\end{multline}
On the other hand, over each $\Sigma_{i, 0}$, for all $\alpha>0$ and a certain $C(\alpha)>0$, one has $|\mu(u_\epsilon)| \leq C(\alpha) |z- z_{i,\epsilon}|^{-2 + \alpha}$, which has finite $L^{p,\delta}$-norm. Hence
\begin{multline}\label{eqn528}
\| d\mu(\uds u_\epsilon)(\uds \xi) - d\mu(u_\epsilon)(\xi) \|_{L_\epsilon^{p,\delta}(\Sigma_{i,\epsilon})}\\ \leq C \left( \| \xi \|_{L^\infty} \| \rho \|_{L^\infty} \| \mu(u_\epsilon) \|_{L_\epsilon^{p,\delta}(\Sigma_{i, \epsilon})} + \| \xi \|_{L^\infty} \| d\mu(u_\epsilon) (\rho ) \|_{L_\epsilon^{p,\delta} } + \|\rho \|_{L^\infty} \| d\mu(u_\epsilon) (\xi) \|_{L_\epsilon^{p,\delta} } \right)\\
\leq C \| {\bm \xi} \|_{L_{m,\epsilon}^{1,p,\delta}} \| {\bm \rho} \|_{L_{m,\epsilon}^{1,p,\delta}}.
\end{multline}
Then \eqref{eqn527} and \eqref{eqn528} give the estimate of the third entry of \eqref{eqn525}, namely
\beq\label{eqn529}
\| d\mu(\uds u_\epsilon)(\uds \xi) - d\mu(u_\epsilon)(\xi) \|_{L_\epsilon^{p,\delta}} \leq C \| {\bm \xi} \|_{L_{m,\epsilon}^{1,p,\delta}} \| {\bm \nu} \|_{L_{m,\epsilon}^{1,p,\delta}}.
\eeq

Now we estimate the first entry of \eqref{eqn525}. We remind the reader that in the case of pseudoholomorphic curves (i.e., when the gauge field is zero), one has the following pointwise estimate
\begin{multline*}
\big| {\it pl}_\rho^{-1} \big( \nabla_s \uds \xi + J \nabla_t \uds \xi + (\nabla_{\uds \xi} J)(\partial_t \uds u_\epsilon)\big) - \big( \nabla_s \xi + J \nabla_t \xi + (\nabla_\xi J)(\partial_t u_\epsilon)\big) \big| \\
 \leq C \left( | du_\epsilon| |\rho | |\xi|  + |\nabla \rho| |\xi| + |\rho | |\nabla \xi| \right).
 \end{multline*}
See the detailed proof of this fact in \cite[Proof of Proposition 3.5.3]{McDuff_Salamon_2004}). In our situation when the covariant derivative is replaced by the covariant derivative $\nabla^{\uds a_\epsilon}$ which has extra terms coming from the gauge field, we claim that similar pointwise estimate still holds. Indeed we have
\begin{multline}\label{eqn530}
\big| {\it pl}_\rho^{-1} I( \uds u_\epsilon, \uds \phi_\epsilon, \uds \phi_\epsilon )(\uds \xi) - I(u_\epsilon, \uds\phi_\epsilon, \uds \psi_\epsilon)(\xi) \big|\\
 \leq C \left( | \partial_s u_\epsilon + {\mc X}_{\uds \phi_\epsilon}| |\rho | |\xi|  + | \partial_t u_\epsilon + {\mc X}_{\uds \psi_\epsilon}||\rho||\xi| + |\nabla^{\uds a_\epsilon} \rho| |\xi| + |\rho | |\nabla^{\uds a_\epsilon} \xi| \right).
\end{multline}
The proof is a tedious reproduction of the proof of \cite[Proposition 3.5.3]{McDuff_Salamon_2004} incorporating the gauge fields; the detail is left to the reader. Then from \eqref{eqn530}, one has
\begin{multline*}
\big\| {\it pl}_\rho^{-1} I( \uds u_\epsilon, \uds \phi_\epsilon, \uds \phi_\epsilon )(\uds \xi) - I(u_\epsilon, \uds\phi_\epsilon, \uds \psi_\epsilon)(\xi) \big\|_{L_\epsilon^{p,\delta}} \\
\leq C \Big( \big( \big\| \partial_s u_\epsilon + {\mc X}_{\uds \phi_\epsilon} \|_{L_\epsilon^{p,\delta}} + \|\partial_t u_\epsilon + {\mc X}_{\uds \psi_\epsilon} \|_{L_\epsilon^{p,\delta}} \big) \| \rho \|_{L^\infty} \| \xi \|_{L^\infty} + \| \nabla^{\uds a_\epsilon} \rho\|_{L_\epsilon^{p,\delta}} \| \xi \|_{L^\infty} + \| \nabla^{\uds a_\epsilon} \xi \|_{L_\epsilon^{p,\delta}} \| \rho \|_{L^\infty} \Big) \\
\leq C \Big( \big( \| \partial_s u_\epsilon + {\mc X}_{\phi_\epsilon} \|_{L_\epsilon^{p,\delta}} + \| \partial_t u_\epsilon + {\mc X}_{\psi_\epsilon} \|_{L_\epsilon^{p,\delta}} +  \| \uds a_\epsilon - a_\epsilon \|_{L_\epsilon^{p,\delta}}\big) \| \rho \|_{L^\infty} \| \xi \|_{L^\infty} \\
+ \| \nabla^{a_\epsilon} \rho\|_{L_\epsilon^{p,\delta}} \| \xi \|_{L^\infty} + \| \nabla^{ a_\epsilon} \xi \|_{L_\epsilon^{p,\delta}} \| \rho \|_{L^\infty}\Big).
\end{multline*}
The first three terms in the last parentheses are uniformly bounded. Hence it follows that 
\beqn
\big\| {\it pl}_\rho^{-1} I( \uds u_\epsilon, \uds \phi_\epsilon, \uds \phi_\epsilon )(\uds \xi) - I(u_\epsilon, \uds\phi_\epsilon, \uds \psi_\epsilon)(\xi) \big\|_{L_\epsilon^{p,\delta}} \leq C \| {\bm \rho} \|_{L_{m,\epsilon}^{1,p,\delta}} \| {\bm \xi} \|_{L_{m,\epsilon}^{1,p,\delta}}.	
\eeqn
Together with \eqref{eqn529}, one obtains
\beqn
\| \hat {\uds D}_\epsilon'( {\bm \xi}) - \hat {\uds D}_\epsilon({\bm \xi}) \|_{L_{\epsilon}^{p,\delta}} \leq C \| {\bm \rho}\|_{L_{m,\epsilon}^{1,p,\delta}} \| {\bm \xi} \|_{L_{m,\epsilon}^{1,p,\delta}}.
\eeqn

\vspace{0.2cm}

\noindent This completes the proof of Proposition \ref{prop55}.

\subsection{Proof of Proposition \ref{prop56}}\label{subsection56}

Now we construct the approximate right inverse along the gauged map ${\bm v}_\epsilon = (u_\epsilon, a_\epsilon)$. In this subsection we abbreviate $\delta_p = 2- \frac{4}{p}$ by $\delta$. 

We make the following assumptions for the purpose of simplifying the notations. Namely, we assume that $m + m^+ = 1$ so that the singular configuration has only two components, one affine vortex ${\bm v}_1$ (over ${\mb H}$) and one holomorphic disk ${\bm v}_0$. The case with more affine vortex components has no more complexity other than notations.\footnote{This assumption violates the stability condition but the stability condition is not required for constructing the approximate solution and proving Proposition \ref{prop54}, \ref{prop55}, and \ref{prop56}.} In this simplified situation, we also assume the coordinate of the nodal point is
\beqn
z_1 = z_{1, \epsilon} = 0.
\eeqn

To construct an approximate right inverse we need another cut-off function. We take $e \in (1, b)$ where $b$ is the number used in choosing the cut-off function $\gamma_i^\epsilon$ (see \eqref{eqn51} and \eqref{eqn52}). Introduce cut-off functions $\chi_0^\epsilon, \chi_1^\epsilon: {\mb H} \to [0, 1]$ satisfying the following conditions:
\begin{align}\label{eqn531}
&\ {\rm supp} \chi_0^\epsilon \subset {\mb H} \setminus B\Big( 0, \frac{1}{e \sqrt\epsilon} \Big),\ \chi_0^\epsilon|_{\Sigma_{0, \epsilon}} \equiv 1;\ &\ {\rm supp} \chi_1^\epsilon \subset B \Big( 0, \frac{e}{\sqrt\epsilon} \Big),\ \chi_1^\epsilon |_{\Sigma_{1,\epsilon}} \equiv 1.
\end{align}
Moreover, require that  
\beq\label{eqn532}
|\nabla \chi_i^\epsilon (z)| \leq \frac{2}{ \log e} \frac{1}{|z|},\ \ i = 0, 1.
\eeq
Notice that when $\nabla \chi_0^\epsilon \neq 0$ or $\nabla \chi_1^\epsilon \neq 0$, the approximate solution agrees with one of the constant gauged map ${\bm c}_1  =  (x_1, 0)$ (see \eqref{eqn55}).

We use parallel transport to identify tangent vectors along nearby maps. By our construction, $u_\epsilon$ is close to $u_{0, \epsilon}$ over $\Sigma_{0, \epsilon}^\gg$ and is close to $u_1$ over $\Sigma_{1,\epsilon}^\gg$. Then one use the parallel transport associated to the connection $\nabla$ chosen by Lemma \ref{connection} to define
\begin{align*}
&\ {\it pl}_0: u_{0,\epsilon}^* TV|_{\Sigma_{0, \epsilon}^\gg} \to u_\epsilon^* TV|_{\Sigma_{0, \epsilon}^\gg},\ &\ {\it pl}_1: u_1^* TV|_{\Sigma_{1, \epsilon}^\gg} \to u_\epsilon^* TV|_{\Sigma_{1, \epsilon}^\gg}.
\end{align*}
Using ${\it pl}_0$, ${\it pl}_1$ and $\chi_0^\epsilon$, $\chi_1^\epsilon$, we define the maps
\begin{align*}
&\ {\it cut}: \hat{\mc E}_\epsilon \to \hat{\mc E}_{0, \epsilon}  \oplus \hat{\mc E}_1 = \hat{\mc E}_{\Gamma, \epsilon},\ &\ {\it paste}: \hat{\mc B}_{\Gamma, \epsilon} \to \hat{\mc B}_\epsilon.
\end{align*}
as follows. For ${\bm \nu} \in \hat{\mc E}_\epsilon$, define ${\it cut}({\bm \nu}) =  ( {\bm \nu}_0, {\bm \nu}_1)$ where
\beqn
{\bm \nu}_i (z) = \left\{ \begin{array}{cc} {\it pl}_i^{-1} [{\bm \nu}(z)],&\ z \in \Sigma_{i, \epsilon}; \\
                                                                      0, &\  z\notin \Sigma_{i,\epsilon}, \end{array}\right.\ \ i = 0, 1.
\eeqn
On the other hand, for $({\bm \xi}_0, {\bm \xi}_1) \in \hat{\mc B}_{\Gamma, \epsilon} \subset \hat{\mc B}_{0,\epsilon} \oplus \hat{\mc B}_1$, define
\beq\label{eqn533}
{\bm \xi}_\epsilon(z):= {\it paste}({\bm \xi}_0, {\bm \xi}_1 )(z) = \left\{ \begin{array}{cl} \displaystyle {\bm \xi}_1 (z), &\ z \in \Sigma_{1,\epsilon}^\ll;\\
\displaystyle \xi_\epsilon^H(z) + {\bm \xi}_\epsilon^G(z), &\ z \in \Sigma_{1,\epsilon}^\gg \setminus \Sigma_{i,\epsilon}^\ll
\end{array} \right.
\eeq
where
\beq\label{eqn534}
\xi_\epsilon^H(z): = \left\{ \begin{array}{cl} \displaystyle  \xi_1^H(z) + \chi_0^\epsilon(z) \Big( {\it pl}_0 ( \xi_0^H (z) - e^{- \lambda_1 \theta_1} \xi_0^H (x_1) \Big),\ &\ z \in \Sigma_{1,\epsilon} \setminus \Sigma_{1,\epsilon}^\ll;\\
                                     \displaystyle \xi_0^H(z) + \chi_1^\epsilon(z) \Big( {\it pl}_1 ( \xi_1^H( z )  - e^{- \lambda_1 \theta_1} \xi_1^H(x_1) \Big),\ &\ z \in \Sigma_{1, \epsilon}^\gg \setminus \Sigma_{1,\epsilon}
\end{array} \right.
\eeq
and
\beq\label{eqn535}
{\bm \xi}_\epsilon^G(z): = \chi_0^\epsilon(z) {\it pl}_0 ({\bm \xi}_0^G) + \chi_1^\epsilon(z) {\it pl}_1 ( {\bm \xi}_1^G ).
\eeq
By the matching condition for $\hat{\mc B}_{\Gamma, \epsilon}$, the value of $\xi_1^H$ an $\xi_2^H$ agrees at the nodal point, therefore ${\bm \xi}_\epsilon$ is continuous and indeed contained in $\hat{\mc B}_\epsilon$. By abuse of notation, we denote the map 
\beqn
({\bm w}, {\bm \xi}_0, {\bm \xi}_1) = ({\bm w}, {\it paste}({\bm \xi}_0, {\bm \xi}_1 ))
\eeqn
still by
\beqn
{\it paste}: W_{\rm def} \oplus \hat{\mc B}_{\Gamma, \epsilon} \to W_{\rm def} \oplus \hat{\mc B}_\epsilon.
\eeqn
Finally, define the ``approximate right inverse'' 
\beq\label{eqn536}
\hat {\mc Q}_\epsilon^{\rm app} = {\it paste} \circ \hat {\mc Q}_{\Gamma, \epsilon} \circ {\it cut}: \hat{\mc E}_\epsilon \to W_{\rm def} \oplus \hat{\mc B}_\epsilon.
\eeq
Here $\hat {\mc Q}_{\Gamma, \epsilon}$ is the operator given by Proposition \ref{inverse}.

The following proposition verifies that for suitable values of all relevant parameters the above operator is indeed an approximate right inverse. 

\begin{prop}\label{prop59}
Suppose $\log e \geq 20 c(p,\delta_p) c(Q)$. Then there exist $\epsilon_4$ and $C$ (which also depend on $b$) such that for $\epsilon \in (0, \epsilon_4)$, with respect to the norm $\| \cdot \|_{L_\epsilon^{p,\delta_p}}$ on $\hat{\mc E}_\epsilon$ and the norm $\| \cdot \|_{L_{m, \epsilon}^{1,p,\delta_p}}$ on $\hat{\mc B}_\epsilon$, there holds
\beq\label{eqn537}
\| \hat {\mc Q}_\epsilon^{\rm app}\|\leq C,
\eeq
and 
\beq\label{approx}
\| \hat {\mc D}_\epsilon \circ \hat {\mc Q}_\epsilon^{\rm app} - {\rm Id} \| \leq \frac{1}{2}.
\eeq
\end{prop}

\begin{proof}
We first show that \eqref{approx} implies \eqref{eqn537}. Indeed, if \eqref{approx} is true, then it follows that $\hat {\mc D}_\epsilon$ is surjective. Since the index of $\hat {\mc D}_\epsilon$ is zero, it follows that $\hat {\mc D}_\epsilon$ is invertible. Then any operator $\hat {\mc Q}_\epsilon$ satisfies \eqref{eqn537} have an upper bound on its norm. 

Now we prove \eqref{approx}. Given ${\bm \nu} \in \hat{\mc E}_\epsilon$, denote
\beqn
({\bm w}, {\bm \xi}_0, {\bm \xi}_1) = \hat {\mc Q}_{\Gamma, \epsilon} ({\it cut}({\bm \nu})) = \hat {\mc Q}_{\Gamma, \epsilon} ({\bm \nu}_0, {\bm \nu}_1).
\eeqn
We would like to show that for appropriate value of $e$ and sufficiently small $\epsilon$, there holds
\beqn
\| \hat{\mc D}_\epsilon ( {\it paste} ({\bm w}, {\bm \xi}_0, {\bm \xi}_1)) - {\bm \nu} \|_{L_\epsilon^{p,\delta}} \leq \frac{1}{2} \Big( \|{\bm \nu}_0 \|_{L_\epsilon^{p,\delta}} + \| {\bm \nu}_1 \|_{L^{p,\delta}} \Big).
\eeqn
By the relation between the weight function $\rho_\epsilon$ and the component-wise weight functions (see
Subsection \ref{subsection52}), this implies \eqref{approx}. 

By definition of $\hat {\mc Q}_{\Gamma, \epsilon}$ and ${\it cut}$, we have
\beqn
{\bm \nu} = {\it pl}_0  \big( \hat {\mc D}_{0, \epsilon} ( {\bm w}, {\bm \xi}_0 ) \big) + {\it pl}_1 \big( \hat {\mc D}_1 ( {\bm w}, {\bm \xi}_1 ) \big).
\eeqn
Therefore
\beqn
\begin{split}
\hat {\mc D}_\epsilon \circ \hat {\mc Q}_\epsilon^{\rm app} ({\bm \nu}) - {\bm \nu} = &\ \hat {\mc D}_\epsilon \Big( {\bm w}, {\it paste}({\bm \xi}_0, {\bm \xi}_1) \Big)  - \Big( {\it pl}_0 \big( \hat {\mc D}_{0, \epsilon}({\bm w}, {\bm \xi}_0) \big) + {\it pl}_1 \big( \hat {\mc D}_1 ({\bm w}, {\bm \xi}_1) \big) \Big) \\
= &\ \hat D_\epsilon \Big( {\it paste}({\bm \xi}_0, {\bm \xi}_1 ) \Big) - \Big( {\it pl}_0 \big( \hat D_{0, \epsilon}( {\bm \xi}_0 ) \big) + {\it pl}_1  \big( \hat D_1 ({\bm \xi}_1) \big) \Big).
\end{split}
\eeqn
Here the last equation follows because over the region where the deformation parameter ${\bm w}$ deforms the equation, the approximate solution agrees with the singular configuration. Now we estimate the last line in different regions as follows.

\vspace{0.2cm}

\noindent (a) Inside $\Sigma_{1,\epsilon} \setminus {\rm supp} \chi_0^\epsilon$ one has ${\bm v}_\epsilon = {\bm v}_{i,\epsilon}'$ (see \eqref{eqn56}) and one has
\beqn
\hat {\mc D}_\epsilon \Big( \hat {\mc Q}_\epsilon^{\rm app} ({\bm \nu}) \Big) - {\bm \nu} = \hat D_{{\bm v}_{1,\epsilon}'} \Big( {\it pl}_1 ({\bm \xi}_1) \Big) - {\it pl}_1 \big( \hat D_1 ({\bm \xi}_1) \big)
\eeqn
If we write ${\bm v}_{i, \epsilon}' = \exp_{{\bm v}_i} {\bm \xi}_{i, \epsilon}'$, then Lemma \ref{lemma51} says that
\beqn
\lim_{\epsilon \to 0} \| {\bm \xi}_{1, \epsilon}'\|_{L_m^{1,p,\delta}} = 0.
\eeqn
Then by using the same method as proving Proposition \ref{prop55}, we have
\beq\label{eqn539}
\Big\| \hat {\mc D}_\epsilon \Big( \hat {\mc Q}_\epsilon^{\rm app} ({\bm \nu}) \Big) - {\bm \nu} \Big\|_{L_\epsilon^{p,\delta}( \Sigma_{1;\epsilon} \setminus {\rm supp} \chi_0^\epsilon)}  \leq C \| {\bm \xi}_{1, \epsilon}'\|_{L_m^{1,p,\delta}} \| {\bm \xi}_1 \|_{L_m^{1,p,\delta}} \leq C(\epsilon) \| {\bm \xi}_1\|_{L_m^{1,p,\delta}}
\eeq
where $C(\epsilon)$ denotes a number which converges to zero as $\epsilon \to 0$.

\vspace{0.1cm}

\noindent (b) Similarly to the above situation, inside $\Sigma_{0,\epsilon} \setminus {\rm supp} \chi_1^\epsilon$, one has
\beqn
\hat {\mc D}_\epsilon \Big( \hat {\mc Q}_\epsilon^{\rm app} ({\bm \nu}) \Big) - {\bm \nu} = \hat D_\epsilon \Big( {\it pl}_0 ({\bm \xi}_0 ) \Big) - {\it pl}_0 \big( \hat D_0 ({\bm \xi}_0) \big)
\eeqn
and
\beq\label{eqn540}
\Big\| \hat {\mc D}_\epsilon (\hat {\mc Q}_\epsilon^{\rm app}({\bm \nu})) - {\bm \nu}   \Big\|_{L^{p, \delta}_\epsilon (\Sigma_{0, \epsilon} \setminus {\rm supp} \chi_1^\epsilon)}  \leq C(\epsilon)   \| {\bm \xi}_0 \|_{L_{m, \epsilon}^{1,p,\delta}}.
\eeq

\vspace{0.1cm}

\noindent (c) In the neck region $N_\epsilon:= {\rm supp} \chi_0^\epsilon \cap {\rm supp} \chi_1^\epsilon $, by our construction of the approximate solution, ${\bm v}_\epsilon = {\bm c}_1 =  ( x_1, 0)$. By definition of ${\it paste}$, one has
\beqn
{\it paste}({\bm \xi}_0, {\bm \xi}_1) = \left[ \begin{array}{c} \chi_0^\epsilon  {\it pl}_0 (\xi_0^H (z)) + \chi_1^\epsilon (z) {\it pl}_1 (\xi_1^H(z)) - \xi_1^H(x_1) \\
\chi_0^\epsilon (z) {\it pl}_0 ({\bm \xi}_0^G(z)) + \chi_1^\epsilon (z) {\it pl}_1 ( {\bm \xi}_1^G(z))
\end{array} \right].
\eeqn
The constant vector $\xi_1^H(x_1)$ is in the kernel of $\hat D_\epsilon$. Hence one has the following simple manipulations (all norms below are the $L_\epsilon^{p,\delta}$-norm over $N_\epsilon$).
\beqn
\begin{split}
 &\ \big\| \hat {\mc D}_\epsilon \big( \hat {\mc Q}_\epsilon^{\rm app}({\bm \nu}) \big) - {\bm \nu}\big\|\\
= &\ \big\|  \hat D_\epsilon \big( \chi_0^\epsilon  {\it pl}_0 ( {\bm \xi}_0 ) + \chi_1^\epsilon {\it pl}_1 ( {\bm \xi}_1) \big) - {\it pl}_0 \big( \hat D_{0, \epsilon} ({\bm \xi}_0) \big) - {\it pl}_1 \big( \hat D_1 ({\bm \xi}_1) \big) \big\|  \\
\leq &\ \big\| \hat D_\epsilon \big( \chi_0^\epsilon  {\it pl}_0 ( {\bm \xi}_0 ) \big) - {\it pl}_0 \big( \hat D_{0,\epsilon} ( {\bm \xi}_0 ) \big) \big\|  + \big\| \hat D_\epsilon \big( \chi_1^\epsilon  {\it pl}_1 ( {\bm \xi}_1) \big) - {\it pl}_1 \big(  \hat D_1 ({\bm \xi}_1 ) \big) \big\| \\
\leq &\ \big\| \hat D_\epsilon \big( \chi_0^\epsilon {\it pl}_0 ( {\bm \xi}_0) \big) - {\it pl}_0  \big( \hat D_{0,\epsilon} ( \chi_0^\epsilon {\bm \xi}_0 )\big) \big\|  + \big\|{\it pl}_0 \big( \hat D_{0,\epsilon} ( {\bm \xi}_0 - \chi_0^\epsilon {\bm \xi}_0 ) \big) \big\|  \\
&\  + \big\| \hat D_\epsilon \big( \chi_1^\epsilon {\it pl}_1 ( {\bm \xi}_1) \big) - {\it pl}_1 \big( \hat D_1 ( \chi_1^\epsilon  {\bm \xi}_1 ) \big) \big\|  + \big\| {\it pl}_1 \big( \hat D_1 ( {\bm \xi}_1 - \chi_1^\epsilon {\bm \xi}_1) \big) \big\|  \\
= &\ \big\| \chi_0^\epsilon \big( \hat D_\epsilon ( {\it pl}_0 ( {\bm \xi}_0) ) - {\it pl}_0 \big( \hat D_{0,\epsilon} ({\bm \xi}_0 )  \big) \big\| + \big\|{\it pl}_0 \big( \hat D_{0,\epsilon} ( {\bm \xi}_0 - \chi_0^\epsilon {\bm \xi}_0 )  \big) \big\|\\
&\  + \big\| \chi_1^\epsilon \big( \hat D_\epsilon ( {\it pl}_1 ( {\bm \xi}_1) ) - {\it pl}_1 \big( \hat D_1 ( {\bm \xi}_1 ) ] \big) \big\|  + \big\|  {\it pl}_1 \big( \hat D_1 ( {\bm \xi}_1 - \chi_1^\epsilon {\bm \xi}_1) \big) \big\|.
\end{split}
\eeqn
The first and the third term can be bounded by the same method of deriving \eqref{eqn539} and \eqref{eqn540}, which gives
\beq
\begin{split}\label{eqn541}
&\ \big\| \chi_0^\epsilon \big(  \hat D_\epsilon ( {\it pl}_0 ( {\bm \xi}_0 )) - {\it pl}_0 ( \hat D_{0, \epsilon} ( {\bm \xi}_0 ) ) \big) \big\|_{L_\epsilon^{p,\delta}(N_\epsilon)}  \leq C(\epsilon)  \| {\bm \xi}_0 \|_{L_{m, \epsilon}^{1,p,\delta}} ;\\
&\  \big\| \chi_1^\epsilon  \big(  \hat D_\epsilon ( {\it pl}_1 ( {\bm \xi}_1) ) - {\it pl}_1 ( \hat D_i ( {\bm \xi}_1 ) ) \big) \big\|_{L_\epsilon^{p,\delta}(N_\epsilon)} \leq C(\epsilon) \| {\bm \xi}_1\|_{L_m^{1,p,\delta}}.
\end{split}
\eeq
To estimate the second and the fourth term, notice that $\hat D_1({\bm \xi}_1) = 0$ outside $\Sigma_{1,\epsilon}$ (the variation of the deformation parameter does not affect this region). Therefore
\beqn
\hat D_1( {\bm \xi}_1 - \chi_1^\epsilon {\bm \xi}_1 ) = \left[ \begin{array}{c} - (\ov\partial \chi_1^\epsilon  ) (\xi_1 - \xi_1^H(x_1)) \\
                                                                                 (\partial_s \chi_1^\epsilon )  \phi_1 + (\partial_t \chi_i^\epsilon ) \psi_1 \\
																																								 (\partial_s \chi_1^\epsilon )  \psi_1 - (\partial_t \chi_1^\epsilon ) \phi_1
\end{array}  \right].
\eeqn
Then for the $(TV)_H$ component, using \eqref{eqn532} and one has the estimate
\begin{multline*}
\big\| \big( \hat D_1 ( {\bm \xi}_1) - \chi_1^\epsilon {\bm \xi}_1) \big)^H \big\|_{L_\epsilon^{p,\delta}}
 \leq \big\| (\nabla \chi_1^\epsilon) (\xi_1^H - \xi_1^H(x_1)) \|_{L_\epsilon^{p,\delta}}\\
   \leq \frac{2}{\log e} \big\| \frac{ \xi_1^H - \xi_1^H(\infty)}{|z|} \big\|_{L^{p,\delta}} \leq \frac{2c(p,\delta)}{\log e} \Big( \| \xi_1^H \|_{L^\infty}  + \| \nabla^{a_1} \xi_1^H \|_{L^{p,\delta}} \Big) \leq \frac{2c(p,\delta)}{\log e} \| {\bm \xi}_1 \|_{L_m^{1,p,\delta}}\\
 \leq  \frac{2c(p,\delta) \| \hat  {\mc Q}_{\Gamma, \epsilon}\| }{\log e} \| {\it cut}({\bm \nu}) \| \leq \frac{2c(p,\delta) c(Q)}{\log e} \| {\bm \nu}\|_{L_\epsilon^{p,\delta}}.
\end{multline*}
Here $c(p,\delta)$ is the constant of the Hardy-type inequality (see \eqref{hardy}) and $c(Q)$ is the upper bound of the operator norm of $\hat {\mc Q}_{\Gamma, \epsilon}$ (see Proposition \ref{inverse}). For the complementary component, one has 
\beqn
\big\| \big( \hat D_1( {\bm \xi}_1 - \chi_1^\epsilon {\bm \xi}_1) \big)^G \big\|_{L_\epsilon^{p,\delta}} \leq C \| \nabla \chi_1^\epsilon \|_{L^\infty} \| {\bm \xi}_1^G \|_{L_\epsilon^{p,\delta}} \leq \frac{C\sqrt\epsilon}{\log e} \| {\bm \xi}_1 \|_{L_m^{1,p,\delta}}.
\eeqn
Therefore when $\epsilon$ is small, one has
\beq
\| \hat D_1 ({\bm \xi}_1) - \chi_1^\epsilon {\bm \xi}_1) \|_{L_\epsilon^{p,\delta}} \leq  \frac{4c(p,\delta) c(Q) }{\log e} \| {\bm \nu} \|_{L_\epsilon^{p,\delta}}.
\eeq
Similarly, when $\epsilon$ is small enough, one has
\beq\label{eqn543}
\Big\| \hat D_{0, \epsilon} ({\bm \xi}_0) - \chi_0^\epsilon {\bm \xi}_0 \|_{L_\epsilon^{p,\delta}} \leq \frac{4c(p,\delta) c(Q) } {\log e} \| {\bm \nu} \|_{L_\epsilon^{p,\delta}}.
\eeq

\vspace{0.2cm}

\noindent Putting \eqref{eqn539}---\eqref{eqn543} together, one can see that when $\epsilon$ is small enough, there holds
\beqn
\| \hat {\mc D}_\epsilon ( \hat {\mc Q}_\epsilon^{\rm app}( {\bm \nu})) - {\bm \nu} \|_{L_\epsilon^{p,\delta}} \leq \frac{ 10 c(p,\delta) c(Q) }{\log e} \| {\bm \nu} \|_{L_\epsilon^{p,\delta}}.
\eeqn
Since we chose $b$ such that $\log b > 100 c(p,\delta) c(Q)$ (see \eqref{eqn51}) and required $e$ such that $e< b$, one can see that \eqref{eqn537} holds when $\log e > 20 c(p, \delta) c(Q)$.
\end{proof}

Proposition \ref{prop56} then follows from Proposition \ref{prop59}.

\section{Surjectivity of the Gluing Map}\label{section6}

In this section we prove that the gluing map is surjective onto an open set of the moduli space. More precisely, we prove the following fact. 

\begin{thm}[Surjectivity of the gluing map] For $\tau$ sufficiently small, the gluing map 
\beqn
{\it Glue}: [0, \tau) \to \ov{\mc M}_{l, l^+}(V, L_V)
\eeqn
defined by Definition \ref{gluingmap} is a homeomorphism onto a neighborhood of $p$. 
\end{thm}

The injectivity part of the above theorem is trivial because the domains of the glued marked affine vortices for different gluing parameters are not isomorphic. On the other hand, the surjectivity part of the above theorem will follow from the uniqueness property contained in the implicit function theorem, if one can verify the following proposition.

\begin{prop}\label{prop62}
Suppose $({\mc C}_k, {\bm v}_k)$ is a sequence of marked affine vortices over ${\mb H}$ with $l$ boundary  markings and $l^+$ interior markings converging to the representative $({\mc C}, {\mc V})$ of ${\bm p}\in {\mc M}_\Gamma(V, L_V)$ (the representative is fixed in Section \ref{section4}). Then after appropriate gauge transformations, ${\bm v}_k$ can be written as 
\beqn
{\bm v}_k = \exp_{{\bm v}_{\epsilon_k}} {\bm \xi}_k
\eeqn
where $\epsilon_k \to 0$, ${\bm v}_{\epsilon_k}$ is the approximate solution constructed in Section \ref{section5} and ${\bm \xi}_k$ are contained in the Banach space $\hat{\mc B}_{\epsilon_k}$ (see \eqref{eqn59} and \eqref{eqn510}). Moreover, there holds
\beqn
\lim_{k \to \infty} \| {\bm \xi}_k \|_{L_{m,\epsilon_k}^{1,p,\delta_p}} = 0.
\eeqn
\end{prop}

\subsection{Proof of Proposition \ref{prop62}}

For simplicity, we make the assumption that the limiting stable marked affine vortex has only one nonconstant affine vortex component, which is an affine vortex component over ${\mb H}$ 
\beqn
{\bm v}_1 = (u_1, a_1) = (u_1, \phi_1, \psi_1).
\eeqn
On the other hand, the disk component is represented by a gauged map
\beqn
{\bm v}_0 = (u_0, a_0) = (u_0, \phi_0, \psi_0).
\eeqn
To fulfill the stability condition we assume that there are other constant affine vortex components (with extra marked points), but they do not affect how we construct the approximate solution and how we construct the gluing map (see Remark \ref{rem52}). Moreover, we regard the domain of ${\bm v}_0$ as the upper half plane ${\mb H}$ and assume that the boundary node between ${\bm v}_1$ and ${\bm v}_0$ is the origin $0 \in {\mb H}$. The case with more nonconstant affine vortex components, or no nonconstant affine vortex components can be obtained with only adding notational complexities. 

The convergence implies that the sequence of domain curves ${\mc C}_k$ converges to ${\mc C}$. Then there is a unique sequence of gluing parameters $\epsilon_k$ and a unique sequence of deformation parameters ${\bm w}_k$ such that 
\beqn
{\mc C}_k = {\mc C}_{\epsilon_k, {\bm w}_k}.
\eeqn
So that the domain of ${\bm v}_k$ is isomorphic to the domain of ${\bm v}_{\epsilon_k}$. 

We assume that the two components ${\bm v}_1$ and ${\bm v}_0$ are in the gauge used in the construction of the approximate solution. In particular, the affine vortices ${\bm v}_1$ satisfies conditions of Lemma \ref{lemma43}. Moreover, the values of ${\bm v}_1$ and ${\bm v}_0$ at the node are equal to a point 
\beqn
x_1 = \lim_{z \to \infty} u_1 (z) = u_0(0) \in L_V.
\eeqn

\subsubsection{The neck region}

We first estimate the distance between ${\bm v}_k$ and ${\bm v}_{\epsilon_k}$ over neck regions. For $b > a \geq 0$ we denote 
\beqn
N_{a, b}:= [a, b]\times [0, \pi]
\eeqn
whose coordinates are $(\tau, \theta)$. Via the exponential map one can identify this strip with the half annulus
\beqn
\exp(N_{a,b})= \Big\{ z \in {\mb H}\ |\ e^a \leq |z| \leq e^b \Big\}.
\eeqn
which is viewed as a subset of the domain of the approximate solution ${\bm v}_{\epsilon_i}$. Then the vortex equation over the half annulus can be rewritten in the polar form as
\beq\label{eqn61}
\begin{split}
\partial_s u + {\mc X}_\phi + J_V(\partial_t u + {\mc X}_\psi) = &\ 0,\\
\partial_s \psi - \partial_t \phi + [\phi, \psi] + e^{2\tau} \mu(u) = &\ 0,\\
u(\partial N_{a, b}) \subset &\ L_V.
\end{split}
\eeq
The energy of ${\bm v}$ then has the expression 
\beqn
E({\bm v}) = \| {\bm v}_s \|_{L^2(N_{a, b})}^2 + \| e^{\tau} \mu(u) \|_{L^2(N_{a,b})}^2.
\eeqn
For any small $\epsilon>0$ and large $T>0$, denote 
\beqn
N_\epsilon(T):= N_{T, - \log \epsilon -T};
\eeqn
\begin{align*}
&\ N_\epsilon^-(T):= N_{T, -\log \sqrt\epsilon},\ &\ N_\epsilon^+(T):= N_{-\log \sqrt\epsilon, -\log \epsilon -T}.
\end{align*}
We prove the following proposition. 

\begin{prop}\label{prop63}
Fix $p \in (2, 4)$. Then for any $\nu>0$, there exist a real number $T_\nu>0$ and an integer $k_\nu>0$ such that for all $k \geq k_\nu$, after a gauge transformation, the image of $u_k |_{N_{T_\nu}^{\epsilon_k}}$ is contained in a neighborhood of $x_1$ so that we can write 
\beqn
{\bm v}_k |_{N_{\epsilon_k}(T_\nu)} = ( \exp_{x_1} \xi_k, \phi_k, \psi_k)
\eeqn
and for ${\bm \xi}_k = (\xi_k, \phi_k, \psi_k)$ there holds
\beqn
\| {\bm \xi}_k \|_{L_{m,\epsilon_k}^{1,p,\delta_p} ( \exp(N_{\epsilon_k}(T_\nu)))} \leq \nu.
\eeqn
\end{prop}

This proposition relies on the fact that the energy decays exponentially over the neck region. Indeed, when the energy of a vortex over a long neck is very small, a certain annulus lemma holds. This property is similar to the property of holomorphic strips: if the energy of a holomorphic strip is less than a threshold, then the energy decays exponentially.

\begin{rem}
We remind the reader again that we follow the convention that in most cases, the letter $C$ represents a constant whose value is allowed to vary from line to line. 
\end{rem}

\begin{lemma}\cite[Proposition 57]{Wang_Xu} \cite[Proposition A.4]{Venugopalan_Xu} \label{annulus} For all (small) $\alpha > 0$, there exist $\varepsilon = \varepsilon(\alpha)>0$ and $C = C(\alpha)>0$ satisfying the following conditions. Given $a \geq 0$, $b \geq a + 2$, suppose ${\bm v} = (u, \phi,\psi)$ is a solution to \eqref{eqn61} over $N_{a, b}$. Suppose $E({\bm v}) \leq \varepsilon$. Then for $s \in [1, \frac{1}{2}( b-a)]$, there holds
\beqn
E({\bm v}; N_{a + s, b-s}) \leq C e^{- (2-2\alpha) s}  E({\bm v}; N_{a,b})
\eeqn
and 
\beqn
{\rm diam}_K(u(N_{a+s, b-s}))  \leq C e^{-(1-\alpha) s} \sqrt{ E({\bm v}; N_{a,b})}.
\eeqn
Here for a subset $S \subset V$, 
\beqn
{\rm diam}_K(S):= \sup_{x, x' \in S} \inf_{g\in K} {\rm dist}( x, gx').
\eeqn
\end{lemma}

The proof of Proposition \ref{prop63} is based on estimates on strips with a fixed length. The convergence towards the limiting stable affine vortex (see the no energy lost condition of Definition \ref{convergence}) implies that
\beqn
\lim_{T \to \infty} \lim_{k \to \infty} E({\bm v}_k; N_{\epsilon_k}(T) )= 0.
\eeqn
Hence for any given $\alpha>0$, one can take $T_0$ sufficiently large such that for all sufficiently large $k$, one has
\beqn
E({\bm v}_k; N_{\epsilon_k}(T_0) ) \leq \varepsilon(\alpha).
\eeqn
Then by Lemma \ref{annulus}, for all $\tau \in [T_0, -\log \epsilon_k - T_0]$, there holds
\beq\label{eqn62}
\sqrt{ E({\bm v}_k; N_{\tau-1, \tau+1})} + {\rm diam}_K ( u_k (N_{\tau-1, \tau+1})) \leq C e^{-(1-\alpha) d_{\epsilon_k}^{T_0} (\tau)}.
\eeq
Here 
\beqn
d_{\epsilon_k}^{T_0} (\tau):= \min \{ \tau - T_0, - \log \epsilon_k - T_0 - \tau\}.
\eeqn

We would like to choose a special gauge over the neck region. Fix a small $r>0$, define
\beqn
H_{x_1}:=H_{x_1}(r):= \exp_{x_1} \left( \{ \xi \in T_{x_1} V,\ d\mu(x_1)(J_V \xi) = 0 \ |\ |\xi| < r \right).
\eeqn
This is a local slice of the $K$-action. The convergence shows that for $T_0$ sufficiently large, the image $u_k (N_{\epsilon_k }(T_0))$ is contained in $K H_{x_1}$. Hence there is a unique gauge transformation on $N_{\epsilon_k }(T_0)$ such that 
\beqn
u_k (N_{\epsilon_k }(T_0) ) \subset H_{x_1}.
\eeqn
In this gauge we write the gauged map ${\bm v}_k$ as $(u_k, a_k)$. We can also write 
\beqn
u_k = \exp_{x_1} \xi_k'. 
\eeqn
The equation \eqref{eqn62} implies that in this gauge there holds
\beq\label{eqn63}
{\rm diam}(u_k (N_{\tau-1, \tau+1}))  \leq C e^{-(1-\alpha) d_{\epsilon_k }^{T_0} (\tau)},\ \forall \tau \in [T_0, -\log \epsilon_k - T_0 ].
\eeq

\begin{lemma}\label{lemma66}
For any $\alpha>0$, there exist $k_\alpha>0$ and $C = C(\alpha)>0$ such that for all $k \geq k_\alpha$ and $\tau \in [T_0, -\log \epsilon_k - T_0]$, there holds 
\beqn
\| \nabla \xi_k'\|_{L^p(N_{\tau-1, \tau+1})} 
+ \| a_k \|_{L^p(N_{\tau-1, \tau+1})} + e^{-\tau} \| \nabla a_k \|_{L^p(N_{\tau-1, \tau+1})} \leq C e^{-(1-\alpha) d_{\epsilon_k }^{T_0} (\tau)}.
\eeqn
\end{lemma}

\begin{proof}
An important ingredient in the proof is that over the strip $N_{\tau-1, \tau+1}$, the function $\sigma$ in the vortex equation \eqref{eqnb1} is $e^{2\tau} \sigma_0$ where $\sigma_0$ is independent of $\tau$. Hence when $\tau$ is large, one can use the {\it a priori} estimates of Lemma \ref{lemmab3}. The rest of the proof is divided into the following two steps. Fix $\alpha>0$.

\vspace{0.2cm}

\noindent {\it Step 1. We prove that there exists $C>0$ such that for sufficiently large $k$ and any $\tau \in [T_0, -\log \epsilon_k - T_0]$, there holds
\beq\label{eqn64}
\| \nabla \xi_k' \|_{L^p(N_{\tau-1, \tau+1})} \leq C e^{-(1-\alpha) d_{\epsilon_k}^{T_0} (\tau)}.
\eeq
}

\noindent Using the map $E_2 = E_2(x_1, \xi_k')$, by Lemma \ref{lemmaa2} and the equation \eqref{eqn61}, one has
\beqn
0 = {\bm v}_{k,s} + J_V {\bm v}_{k,t} = {\mc X}_{\phi_k'} + E_2(\nabla_s \xi_k') + J_V( {\mc X}_{\psi_k'} + E_2(\nabla_t \xi_k')).
\eeqn
Projecting onto the $(TV)_H$ direction, one obtains
\begin{multline*}
|\nabla_s \xi_k' + J_V \nabla_t \xi_k'| \leq C \left( |d\mu(u_k)({\bm v}_{k,s})| + |d\mu(u_k )({\bm v}_{k ,t})| \right) \\
\leq C e^{-\frac{2}{p} \tau} \sqrt{ E({\bm v}_k; N_{\tau-2, \tau+2})} \leq C e^{-(1 + \frac{2}{p} - \alpha) d_{\epsilon_k}^{T_0} (\tau)}.
\end{multline*}
Here the last line follows from Lemma \ref{lemmab3}, \eqref{eqnb6}, and \eqref{eqn62}. Furthermore, \eqref{eqn63} implies that there exists $\exp_{x_1} (\xi_k'')\in H_{x_1}(r)$ such that 
\beqn
\| \xi_k' - \xi_k'' \|_{L^\infty(N_{\tau-2, \tau+2})} \leq C {\rm diam}(u_k (N_{\tau-2, \tau+2})) \leq C e^{-(1-\alpha ) d_{\epsilon_k}^{T_0} (\tau)}.
\eeqn
The fact that $L_V$ is totally geodesic implies that $\xi_k' - \xi_k''|_{\partial N_{\tau-2, \tau+2}} \in TL_V$. Then by the interior elliptic estimate for the $\ov\partial$-operator with totally real boundary condition, one has 
\beqn
\| \nabla \xi_k' \|_{L^p(N_{\tau-1, \tau+1})} = \| \nabla( \xi_k' - \xi_k'') \|_{L^p(N_{\tau-1, \tau+1})} \leq C e^{-(1-\alpha ) d_{\epsilon_k}^{T_0} (\tau)}.
\eeqn

\vspace{0.2 cm}

\noindent {\it Step 2. We prove that there exists $C>0$ such that for all sufficiently large $k$ and any $\tau \in [T_0, - \log \epsilon_k - T_0]$, there holds 
\beq\label{eqn65}
\| a_k \|_{L^p(N_{\tau-1, \tau+1})} + e^{-\tau} \| \nabla a_k \|_{L^p(N_{\tau-1, \tau+1})} \leq C e^{-(1-\alpha) d_{\epsilon_k}^{T_0} (\tau)}.
\eeq
}
\vspace{0.1cm}

\noindent It is similar to the proof of Lemma \ref{lemmab7}. Indeed, by Lemma \ref{lemmaa2}, one has
\beqn
{\mc X}_{\phi_k} = {\bm v}_{k,s} - E_2(x_1, \xi_k') \nabla_s \xi_k'.
\eeqn
Then by the annulus lemma and \eqref{eqn64}, one has
\beqn
\| \phi_k \|_{L^p(N_{\tau-1, \tau+1})} \leq C \left( \| {\bm v}_{k,s}\|_{L^p(N_{\tau-1, \tau+1})} + \| \nabla_s \xi_k' \|_{L^p(N_{\tau-1, \tau+1})} \right)\leq C  e^{-(1-\alpha) d_{\epsilon_k}^{T_0} (\tau)}.
\eeqn
To estimate the derivative of $\phi_k$, we use Lemma \ref{lemmaa3} and obtain a calculation similar to \eqref{eqnb28}, namely
\beqn
{\mc X}_{\nabla_t \phi_k} = \nabla_t {\bm v}_{k,s} - E_{22}( \nabla_s \xi_k', \nabla_t \xi_k') - E_2 \nabla_t \nabla_s \xi_k'.
\eeqn
Using the fact that $\xi_k'$ is in the kernel of $d\mu(x_1) \circ J_V$, one obtains 
\begin{multline*}
\| \nabla_t \phi_k  \|_{L^p(N_{\tau-1, \tau+1})} \leq C \| d\mu(x_1)(J_V  E_2^{-1} {\mc X}_{\nabla_s \xi_k'}) \|
_{L^p(N_{\tau-1, \tau+1})} \\
\leq   C \left( \| \nabla_t {\bm v}_{k,s} \|_{L^p(N_{\tau-1, \tau+1})} +  \| \nabla \xi_k' \|_{L^{2p}(N_{\tau-1, \tau+1})}^2 + \| d\mu(x_1) (J_V \nabla_t \nabla_s \xi_k') \|_{L^p(N_{\tau-1, \tau+1})} \right)\\
\leq C \left( \| \nabla_t^{a_k} {\bm v}_{k,s} \|_{L^p(N_{\tau-1, \tau+1})} + \| \beta_k' \|_{L^p(N_{\tau-1, \tau+1})} \| {\bm v}_{k, s} \|_{L^\infty(N_{\tau-1, \tau+1})} + \| \nabla \xi_k' \|_{L^{2p}(N_{\tau-1, \tau+1})}^2 \right)\\
\leq C e^{(1 - \frac{2}{p}) \tau} \sqrt{ E({\bm v}_k; N_{\tau-1, \tau+1})} \leq C e^\tau e^{-( 1 + \frac{2}{p} - \alpha) d_{\epsilon_k}^{T_0} (\tau)}.
\end{multline*}
Other components of derivatives of $a_k$ can be estimated similarly. Then \eqref{eqn65} is proved. 
\end{proof}

\begin{proof}[Proof of Proposition \ref{prop63}] The proof is a series of calculations. Choose $\alpha \in (0, \frac{2}{p})$. 

\vspace{0.2cm}

\noindent {\it Step 1. We estimate the difference $\beta_k = a_k - a_{\epsilon_k} = a_k - \gamma_1^{\epsilon_k} \check a_1$ of connections over the inner half $\exp(N_{\epsilon_k}^-(T))$ of the neck $\exp(N_{\epsilon_k}(T))$.}

\vspace{0.1cm}

\noindent Recall that the Sobolev norm is weighted by the weight function $\rho_{\epsilon_k}$ to the power $\frac{\delta_p}{2} = 1- \frac{2}{p}$. The value of $\rho_{\epsilon_k}$ over the strip $N_{\tau-1, \tau+1} \subset N_{\epsilon_k}^-(T)$ is roughly $e^{2\tau}$. Moreover, when changing the cylindrical metric to the Euclidean metric, the $L^p$-norms of 1-forms is multiplied by roughly $e^{-(1-\frac{2}{p})\tau}$. Then by Lemma \ref{lemma66}, one has  
\begin{multline}\label{eqn66}
\| a_k \|_{L^{p, \delta_p}(\exp(N_{\epsilon_k}^-(T)))} \leq C \left( \sum_{\tau \in {\mb Z} \cap [T, -\log \sqrt{\epsilon_k}] } e^{\delta_p \tau} e^{- (1-\frac{2}{p}) \tau} \| a_k \|_{L^p(N_{\tau-1, \tau+1})}\right)\\
\leq C\left( \sum_{\tau \in {\mb Z} \cap [T, +\infty) }  e^{(1- \frac{2}{p}) \tau} e^{-(1-\alpha) (\tau - T_0)} \right)\leq C e^{-(\frac{2}{p}-\alpha)T}.
\end{multline}
On the other hand, by Lemma \ref{lemma43}, one has 
\beqn
\| \gamma_1^{\epsilon_k} \check a_k \|_{L^{p,\delta_p}(\exp(N_{\epsilon_k}^-(T))} \leq \| \check a_1 \|_{L^{p,\delta_p}({\mb A} \setminus B_{e^T})} \leq C e^{-(\delta_0 - \delta_p) T}. 
\eeqn
Then for $\varepsilon = \min \{  \frac{2}{p} - \alpha, \delta_0 - \delta_p \}$ one has
\beq\label{eqn67}
\| \beta_k \|_{L^{p,\delta_p}(\exp( N_{\epsilon_k, T}^-) )} \leq \| a_k \|_{L^{p,\delta_p}( \exp( N_{\epsilon_k, T}^-) )} + \| \gamma_1^{\epsilon_k} \check a_1 \|_{L^{p,\delta_p}( \exp( N_{\epsilon_k, T}^-) )} \leq C e^{-\varepsilon T}.
\eeq
Using Lemma \ref{lemma66} and the fact that $| a_{\epsilon_k}|$ is bounded, calculation similar to \eqref{eqn66} yields
\beqn
\| \nabla^{a_{\epsilon_k}} a_k \|_{L^{p,\delta_p}(\exp(N_{\epsilon_k}^-(T)))} \leq C e^{- (\frac{2}{p} - \alpha) T}.
\eeqn
Moreover, 
\begin{multline*}
\| \nabla^{a_{\epsilon_k}} (\gamma_1^{\epsilon_k} \check a_1)\|_{L^{p,\delta_p}(\exp(N_{\epsilon_k}^-(T))}\\
 \leq \| \gamma_1^{\epsilon_k} \nabla^{a_{\epsilon_k}} \check a_1 \|_{L^{p,\delta_p}(\exp(N_{\epsilon_k}^-(T))} + \| (\nabla \gamma_1^{\epsilon_k}) \check a_1 \|_{L^{p,\delta_p}(\exp(N_{\epsilon_k}^-(T))}\\
\leq C \left( \| \nabla \check a_1 \|_{L^{p,\delta_p}({\mb A} \setminus B_{e^T})} + \| \check a_1 \|_{L^{p,\delta_p}({\mb A} \setminus B_{e^T})} \right) \leq C e^{-(\delta_0 - \delta_p) T}.
\end{multline*}
Therefore
\beq\label{eqn68}
\| \nabla^{a_{\epsilon_k}} \beta_k \|_{L^{p,\delta_p}(\exp(N_{\epsilon_k}^-(T))} \leq C \left( e^{-\varepsilon T} + (\epsilon_k)^\varepsilon \right).
\eeq

\vspace{0.2cm}

\noindent {\it Step 2. We estimate $\beta_k$ over the other half of the neck $\exp(N_{\epsilon_k}^+(T))$.}

\vspace{0.1cm}

\noindent Notice that in this region the weight function $\rho_{\epsilon_k}$ is a constant $\epsilon_k^{-1}$. Then one has 
\begin{multline*}
\| a_k \|_{L^p(\exp( N_{\epsilon_k}^+(T))} \leq C \left( \sum_{\tau \in {\mb Z} \cap [-\log \sqrt{\epsilon_k}, - \log \epsilon_k - T]} \epsilon_k^{-\frac{\delta_p}{2}} e^{-(1-\frac{2}{p}) \tau} \| a_k \|_{L^p(N_{\tau-1, \tau+1})}\right)\\
\leq C \left( \sum_{\tau \in {\mb Z} \cap [-\log \sqrt{\epsilon_k}, -\log \epsilon_k - T]} e^{ (1- \frac{2}{p})(\log \epsilon_k -  \tau)} e^{-(1-\alpha ) ( -\log \epsilon_k - T_0 - \tau)  } \right) \\
\leq C \sum_{\tau \in {\mb Z} \cap (-\infty, -\log \epsilon_k - T ]} e^{-(\frac{2}{p}-\alpha ) (\log \epsilon_k - \tau)} \leq C e^{-(\frac{2}{p}-\alpha ) T}.
\end{multline*}
By the exponential decay of the connection form $\check a_0$, one obtains 
\beq\label{eqn69}
\| \beta_k \|_{L^p(\exp(N_{\epsilon_k}^+(T))} \leq C \| a_k \|_{L^p(\exp(N_{\epsilon_k}^+(T))} + \| s_\epsilon^* \check a_0 \|_{L^p(\exp(N_{\epsilon_k}^+(T))} \leq C \left( e^{-(\frac{2}{p}-\alpha) T} + e^{-(1- \delta_p) T} \right).
\eeq
We omit the calculation for the derivative $\nabla^{a_{\epsilon_k}}  \beta_k$. Indeed we will obtain
\beq\label{eqn610}
\| \nabla^{a_{\epsilon_k}} \beta_k \|_{L^p(\exp(N_{\epsilon_k}^+(T))} \leq C \left( e^{- (\frac{2}{p}-\alpha) T} + e^{- (1- \delta_p) T}\right).
\eeq

\vspace{0.2cm}

\noindent {\it Step 3. We estimate $\nabla \xi_k$.}

\vspace{0.1cm}

\noindent By very similar calculations, one can prove that
\beqn
\| \nabla \xi_k' \|_{L^{p,\delta_p}(\exp(N_{\epsilon_k}(T))} \leq C e^{- (\frac{2}{p}-\alpha) T}
\eeqn
and 
\beqn
\| \nabla u_{\epsilon_k} \|_{L^{p,\delta_p}(\exp(N_{\epsilon_k}(T))} \leq C e^{-(\delta_0 - \delta_p) T}.
\eeqn
Since $u_k$ and $u_{\epsilon_k}$ are both close to $x_1$, $\xi_k$ is a function of $u_k$ and $u_{\epsilon_k}$. Then
\beqn
|\nabla \xi_k| \leq C (|\xi_k'||\nabla u_{\epsilon_k }| + |\nabla u_k ||\gamma_i^{\epsilon_k} \check \vartheta_1|) \leq C ( |\xi_k'| |\nabla \check \vartheta_1| + |\xi_k'||\check \vartheta_1| |\nabla \gamma_i^{\epsilon_k}| + |\nabla \xi_k'| |\check \vartheta_1|).
\eeqn
Therefore, one has
\begin{multline*}
\| \nabla^{a_{\epsilon_k}}  \xi_k \|_{L^{p,\delta_p}(\exp(N_{\epsilon_k}^-(T))} \leq C \left( \| \nabla \xi_k \|_{L^{p,\delta_p}(\exp(N_{\epsilon_k}^-(T))} + \|a_{\epsilon_k}\|_{L^{p,\delta_p}(\exp(N_{\epsilon_k}^-(T))}  \| \xi_k \|_{L^\infty (N_{\epsilon_k}^-(T))} \right)\\
\leq C \left( \| \nabla \check \vartheta_1\|_{L^{p,\delta_p}(\exp(N_{\epsilon_k, T}^-))} + \| \nabla \xi_k' \|_{L^{p,\delta_p}(\exp(N_{\epsilon_k}^-(T))} + e^{-(\delta_0 - \delta_p) T} + \| (\nabla \gamma_1^{\epsilon_k}) \check \vartheta_1 \|_{L^{p,\delta_p}(\exp(N_{\epsilon_k}^-(T))} \right)\\
\leq C \left( e^{-(\delta_0 - \delta_p) T} + e^{-(\frac{2}{p}-\alpha) T} + \| (\nabla \gamma_1^{\epsilon_k}) \check \vartheta_1 \|_{L^{p,\delta_p}(\exp(N_{\epsilon_k}^-(T))} \right).
\end{multline*}
To estimate the last term, using the fact that $|\nabla \gamma_1^{\epsilon_k}| \approx \sqrt{\epsilon_k}$, the area of the support of $\nabla \gamma_1^{\epsilon_k}$ is a multiple of $\epsilon_k^{-1}$, and \ref{exponential}, one has
\begin{multline*}
\| (\nabla \gamma_1^{\epsilon_k}) \check \vartheta_1 \|_{L^{p,\delta_p}(\exp(N_{\epsilon_k}^-(T))} \leq C \sqrt{\epsilon_k} \cdot ({\rm Area} ({\rm supp} \nabla \gamma_1^{\epsilon_k} ) )^{\frac{1}{p}} \cdot \sup_{{\rm supp} \nabla \gamma_1^{\epsilon_k}} \left( |\check \vartheta_1| \rho_{\epsilon_k}^{\frac{\delta_p}{2}}\right)\\
\leq C \epsilon_k^{\frac{1}{2}  - \frac{1}{p} + \frac{\delta_0}{2} - \frac{\delta_p}{4} - \frac{\delta_p}{2}} = C \epsilon_k^{\frac{\delta_0 - \delta_p}{2}}.
\end{multline*}
Therefore one has
\beqn
\| \nabla^{a_{\epsilon_k}}  \xi_k \|_{L^{p,\delta_p}(\exp(N_{\epsilon_i}^-(T))}  \leq C \left( e^{-(\frac{2}{p}-\alpha) T} + e^{-(\delta_0-\delta_p)T} + \epsilon_k^{\frac{\delta_0-\delta_p}{2}}\right).
\eeqn
A similar calculation yields an estimate for the outer half of the neck
\beqn
\| \nabla^{a_{\epsilon_k}}  \xi_k \|_{L^{p,\delta_p}(\exp(N_{\epsilon_k}^+(T))}  \leq C \left( e^{-(\frac{2}{p}-\alpha)T} + e^{-(1-\delta_p) T} + \epsilon_k^{\frac{1-\delta_p}{2}}\right).
\eeqn
Therefore, for a certain $\varepsilon >0$, there holds
\beq\label{eqn611}
\| \nabla^{a_{\epsilon_k}}  \xi_k \|_{L^{p,\delta_p}(\exp(N_{\epsilon_k} (T))}  \leq C \left( e^{- \varepsilon T} + (\epsilon_k)^\varepsilon \right).
\eeq

\vspace{0.2cm}

Lastly notice that 
\beqn
\lim_{T \to \infty} \lim_{k \to \infty} \| \xi_k \|_{L^\infty(\exp(N_{\epsilon_k}(T))} = 0.
\eeqn
Then following \eqref{eqn67}---\eqref{eqn611}, for given $\nu>0$, one can choose $k_\nu$ sufficiently large and $T_\nu$ sufficiently large such that for $k \geq k_\nu$ there holds
\beqn
\| {\bm \xi}_k \|_{L_{m,\epsilon_k}^{1,p,\delta_p}(\exp(N_{\epsilon_k} (T))} \leq \nu.
\eeqn
This finishes the proof of Proposition \ref{prop63}.
\end{proof}

\begin{rem}
One can prove an analogue of Proposition \ref{prop63} if there are nonconstant affine vortex components over ${\mb C}$. We mention that in that case one needs to use the annulus lemma for cylinders proved by Ziltener in \cite{Ziltener_Decay}.
\end{rem}

\subsubsection{The affine vortex component}

Fix the $T_\nu$ given by Proposition \ref{prop63}. Now we compare ${\bm v}_k$ and ${\bm v}_{\epsilon_k}$ in fixed compact region complementary to the necks, namely the half disk $B_{e^{T_\nu}}^+(0)\subset {\mb H}$. The convergence of the affine vortex component implies that ${\bm v}_k |_{B_{e^{T_\nu+1}}}$ converges modulo gauge transformation to ${\bm v}_1|_{B_{e^{T_\nu+1}}}$ with all derivatives. Suppose in this gauge the gauged map ${\bm v}_k$ is $(u_k^\sim, a_k^\sim)$. This may not agree with the gauge we specified in the proof of Proposition \ref{prop63}, but the two can be glued together in the following way without making the estimate essentially worse.

Over the half annulus $\exp(N_{T_\nu, T_{\nu} + 1})$, we can write 
\begin{align*}
&\ u_k = \exp_{u_1} \xi_k,\ &\ u_k^\sim = \exp_{u_1} \xi_k^\sim.
\end{align*}
Since $u_1(\exp(N_{T_\nu, T_\nu+1}))$ is contained in a neighborhood of $\mu^{-1}(0)$ where the $K$-action is free and $\xi_k, \xi_k^\sim$ are very small pointwise, there is a unique gauge transformation $g_k = \exp h_k$ where $h_k: \exp(N_{T_\nu, T_\nu +1}) \to K$ such that 
\beqn
u_k  = g_k u_k^\sim.
\eeqn
Moreover, one has the estimates
\beqn
\begin{split}
\| \xi_k \|_{L^p(\exp(N_{T_\nu, T_\nu +1}))} + \| \nabla^{a_{\epsilon_k}} \xi_k \|_{L^p(\exp(N_{T_\nu, T_\nu +1}))} \leq &\ C \nu,\\
\| \xi_k^\sim \|_{L^p(\exp(N_{T_\nu, T_\nu +1}))} + \| \nabla^{a_{\epsilon_k}} \xi_k^\sim \|_{L^p(\exp(N_{T_\nu, T_\nu +1}))} \leq &\ C \nu
\end{split}
\eeqn
where $C$ is independent of $i$ and $\nu$. Hence one can obtain the estimate 
\beqn
\| h_k \|_{L^p(\exp(N_{T_\nu, T_\nu +1}))} + \| \nabla^{a_{\epsilon_k} } h_k \|_{L^p(\exp(N_{T_\nu, T_\nu +1}))} \leq C \nu.
\eeqn
Then by using a cut-off function one can concatenate $u_k$ and $u_k^\sim$. Therefore one has the following proposition. 

\begin{prop}\label{prop68}
For any $\nu>0$, there exist $k_\nu>0$ and $T_\nu>0$ such that for all $k \geq k_\nu$, after appropriately gauge transforming ${\bm v}_k$, over the region 
\beqn
W_{k, T_\nu}:= \Big\{ z \in {\mb H}\ |\ |z| \leq \epsilon_k^{-1} e^{-T_\nu} \Big\} 
\eeqn
one can write 
\beqn
{\bm v}_k = \exp_{{\bm v}_{\epsilon_k}} {\bm \xi}_k
\eeqn
and there holds
\beqn
\| {\bm \xi}_k \|_{L_{m,\epsilon_k}^{1,p,\delta_p}(W_{k, T_\nu})} \leq \nu.
\eeqn
\end{prop}

\subsubsection{Near infinity}

Now we compare ${\bm v}_k$ and ${\bm v}_{\epsilon_k}$ near the infinity of the holomorphic disk component. The convergence towards the limiting stable affine vortex implies that 
\beqn
\lim_{T\to \infty} \lim_{k \to \infty} E({\bm v}_k; {\mb H} \setminus B_{\epsilon_k^{-1} e^{T}} ) = 0.
\eeqn
Hence for $T_0>0$ and $k_0$ sufficiently large, for all $k \geq k_0$ one has
\beqn
E({\bm v}_k; {\mb H} \setminus B_{\epsilon_k^{-1} e^{T_0}} ) \leq \varepsilon(\alpha)
\eeqn
where $\varepsilon(\alpha)$ is the constant in the annulus lemma (Lemma \ref{annulus}). One can also define the slice $H_{x_\infty}(r_0)$ of the $K$-action through $x_\infty$. Then one can put ${\bm v}_k$ in a gauge such that 
\beqn
u_k({\mb H} \setminus B_{\epsilon_k^{-1} e^{T_0}} ) \subset H_{x_\infty}(r_0).
\eeqn
In this gauge one can also write 
\beqn
{\bm v}_k = \exp_{{\bm v}_{\epsilon_k}} {\bm \xi}_k.
\eeqn
Then it is a similar procedure to prove the following fact.

\begin{prop}\label{prop69}
For all $\nu>0$, there exist $T_\nu>0$ and $k_\nu>0$ such that for all $k \geq k_\nu$ there holds
\beqn
\| {\bm \xi}_k \|_{L_{m,\epsilon_k}^{1,p,\delta_p} ({\mb H}\setminus B_{\epsilon_k^{-1} e^{T_\nu}})} \leq \nu.
\eeqn
\end{prop}

The details are left to the reader. 

\subsubsection{In the compact region}

It remains to estimate the difference of ${\bm v}_k$ and ${\bm v}_{\epsilon_k}$ over
\beqn
Z_{k, T_\nu}:= \Big\{ z \in {\mb H}\ |\ \epsilon_k^{-1} e^{-T_\nu} \leq |z| \leq \epsilon_k^{-1} e^{T_\nu}\Big\}.
\eeqn
Notice that in this region the approximate solution ${\bm v}_{\epsilon_k}$ is equal to the rescaling of the disk component ${\bm v}_\infty$ over the region 
\beqn
Z_{T_\nu}:= \Big\{ z \in {\mb H}\ |\ e^{-T_\nu} \leq |z| \leq e^{T_\nu}\Big\}.
\eeqn
Recall that $s_{\epsilon_k}: Z_{k, T_\nu} \to Z_{T_\nu}$ is the rescaling $z \mapsto \epsilon_k z$. Denote 
\beqn
\tilde {\bm v}_k := (s_{\epsilon_k}^{-1})^* {\bm v}_k = (\tilde u_k, \tilde a_k).
\eeqn
This is a sequence of solutions to the equation \eqref{eqnb1} over $Z_{T_\nu}$ for $\epsilon = \epsilon_k$ and $\sigma \equiv 1$. The convergence towards $u_\infty$ implies that modulo gauge transformation, $\tilde u_k$ converges uniformly to $u_\infty$. Then one can write 
\beqn
\tilde {\bm v}_k = \exp_{{\bm v}_\infty} \tilde {\bm \xi}_k,\ {\rm where}\ \tilde {\bm \xi}_k = (\tilde \xi_k, \tilde \beta_k) \in \Gammait( Z_{T_\nu}, u_\infty^* TV \oplus {\mf k} \oplus {\mf k}).
\eeqn
Theorem \ref{thmb2} implies that after appropriate gauge transformations, there holds 
\beqn
\lim_{k \to \infty} \left( \| \tilde \beta_k \|_{L^p(Z_{T_\nu})} + \epsilon_k \| \nabla^{a_0} \tilde \beta_k \|_{L^p(Z_{T_\nu})} + \| \tilde \xi_k^H \|_{L^\infty(Z_{T_\nu})} +  \epsilon_k^{-1} \| \tilde \xi_k^G\|_{L^p(Z_{T_\nu})} + \|\nabla^{a_0} \tilde \xi_k \|_{L^p(Z_{T_\nu})} \right)  = 0.
\eeqn
Since $Z_{T_\nu}$ is a compact region which is independent of $\epsilon_k$, by using the auxiliary norm \eqref{eqn425}, one can see that the above is equivalent to 
\beqn
\lim_{k \to \infty} \left( \| \tilde \xi_k^H \|_{L^p (Z_{T_\nu})} +  \| (\epsilon_k^{-1} \tilde \xi_k^G, \tilde \beta_k ) \|_{L_{1,p;\epsilon_k}^{\rm aux}(Z_{T_\nu})} \right) = 0.
\eeqn
Then by \eqref{eqn426}, the above implies the following result. 

\begin{prop}\label{prop610}
For any $\nu>0$, there exist $k_\nu>0$ such that for all $k \geq k_\nu$, after appropriately gauge transforming ${\bm v}_k$, one can write 
\beqn
{\bm v}_k = \exp_{{\bm v}_\infty} {\bm \xi}_k,\ {\rm where}\ {\bm \xi}_k \in \Gammait( Z_{k, T_\nu}, u_{\epsilon_k}^* TV \oplus {\mf k} \oplus {\mf k})
\eeqn
and there holds
\beqn
\| {\bm \xi}_k \|_{L_{m,\epsilon_k}^{1,p, \delta_p}(Z_{k, T_\nu})}  \leq \nu.
\eeqn
\end{prop}

Lastly, the gauge transformations in Proposition \ref{prop68}, \ref{prop69}, and \ref{prop610} may not agree. One can use the same argument as in the proof of Proposition \ref{prop68} to glue these gauge transformations to obtain a sequence of global gauge transformations over ${\mb H}$ such that after gauge transforming ${\bm v}_k$, the statement of Proposition \ref{prop62} holds.

\appendix







\section{Derivatives of the Exponential Map}

This appendix is nearly identical to part of \cite[Appendix C]{Gaio_Salamon_2005}, which we include here for convenience. Let $M$ be a complete Riemannian manifold. For any $x \in M$, $\xi \in T_x M$, and $i, j \in \{1, 2\}$ there exist linear maps
\begin{align*}
&\ E_i(x, \xi): T_x M \to T_{\exp_x \xi} M, &\ E_{ij}(x, \xi): T_x M \oplus T_x M \to T_{\exp_x \xi} M
\end{align*}
defined by the following identities
\beq\label{eqna1}
d \exp_x \xi = E_1(x, \xi) dx + E_2(x, \xi) \nabla \xi. 
\eeq
\beqn
\begin{split}
\nabla E_1(x, \xi) w = &\ E_{11}(x, \xi)(w, dx) + E_{12}(x, \xi) (w, \nabla \xi) + E_1(x, \xi) \nabla w,\\
\nabla E_2(x, \xi) w = &\ E_{21}(x, \xi) (w, dx) + E_{22}(x, \xi)(w, \nabla \xi) + E_2(x, \xi) \nabla w.
\end{split}
\eeqn

\begin{lemma}\cite[Lemma C.1]{Gaio_Salamon_2005}
For any $x \in X$, $\xi \in T_x X$, $\eta \in {\mf k}$, 
\beqn
{\mc X}_\eta (\exp_x \xi) = E_1(x, \xi) {\mc X}_\eta(x) + E_2(x, \xi) \nabla_\xi {\mc X}_\eta.
\eeqn
\end{lemma}

Denote the infinitesimal action at a point $x\in M$ by
\beqn
L_x: {\mf k} \to T_x M,\ L_x(\eta) = {\mc X}_\eta(x).
\eeqn
Let $M^* \subset M$ be the open subset where the infinitesimal $K$-action is injective. Namely
\beqn
M^*:= \{ x \in M\ |\ {\rm ker} L_x = \{0\} \}.
\eeqn
Let $\Omega \subset {\mb H}$ be an open subset and $\uds u: \Omega \to M^*$ be a $C^1$-map. Then there exists a unique connection form $\uds a = \uds \phi ds + \uds \psi dt$ on $\Omega$ defined via the relation 
\beqn
\partial_s \uds u + {\mc X}_{\uds \phi},\ \partial_t \uds u + {\mc X}_{\uds \psi} \in ({\rm Im} L_u)^\bot.
\eeqn
Denote 
\beqn
\uds {\bm v} = (\uds u, \uds a) =  (\uds u, \uds \phi, \uds \psi)
\eeqn
and call it the gauged map induced from $\uds u$. 

\begin{lemma} \label{lemmaa2} \cite[Lemma C.3]{Gaio_Salamon_2005}
Let $\uds u: \Omega \to M^*$ be a $C^2$-map inducing the gauged map $\uds {\bm v}$ as above. Let ${\bm v} = (u, \phi, \psi)$ be another gauged map with 
\beqn
\phi = \uds \phi + \eta,\ \psi = \uds \psi + \zeta,\ u = \exp_{\uds u} \xi,\ \xi \in \Gammait( \Omega, \uds u^* TM).
\eeqn
Then one has
\begin{align*}
&\ {\mc X}_\eta = {\bm v}_s - E_1(\uds u, \xi)( \uds {\bm v}_s)  - E_2( \uds u, \xi) (\nabla_s^{\uds a} \xi),\ &\ \ {\mc X}_\zeta = {\bm v}_t - E_1( \uds u, \xi) (\uds {\bm v}_t) - E_2( \uds u, \xi) (\nabla_t^{\uds a} \xi).
\end{align*}
\end{lemma}

\begin{lemma}\label{lemmaa3}\cite[Lemma C.5]{Gaio_Salamon_2005}
Under hypothesis of Lemma \ref{lemmaa2}, suppose in addition 
\beqn
\xi \in ({\rm Im} L_{\uds u})^\bot,
\eeqn
then (abbreviate $E_i(\uds u, \xi), E_{ij}(\uds u, \xi)$ by $E_i$ and $E_{ij}$)
\beqn
\begin{split}
L_u {\nabla_t^a \eta}  = &\ \nabla_t^a {\bm v}_s + \nabla_{{\mc X}_\eta} {\mc X}_\zeta - \nabla_{{\bm v}_t} {\mc X}_\eta - \nabla_{{\bm v}_s} {\mc X}_\zeta  -E_{11} ( \uds {\bm v}_s, \uds {\bm v}_t) - E_{12} ( \uds {\bm v}_s, \nabla_t^{\uds a} \xi) \\
                             &\ - E_{21}(\nabla_s^{\uds a } \xi, \uds {\bm v}_t ) - E_{22} (\nabla_s^{\uds a} \xi, \nabla_t^{\uds a} \xi) - E_1 \nabla_t^{\uds a} \uds {\bm v}_s - E_2 \nabla_t^{\uds a}\nabla_s^{\uds a} \xi;\\
L_u \nabla_s^a \eta =  &\ \nabla_s^a {\bm v}_s + \nabla_{{\mc X}_\eta} {\mc X}_\eta - 2 \nabla_{{\bm v}_s} {\mc X}_\eta - E_{11} ( \uds {\bm v}_s, \uds {\bm v}_s )- E_{12} ( \uds {\bm v}_s, \nabla_s^{\uds a} \xi) \\
 &\ - E_{21} (\nabla_s^{\uds a} \xi, \uds {\bm v}_s ) - E_{22} ( \nabla_s^{\uds a} \xi, \nabla_s^{\uds a} \xi) - E_1 \nabla_s^{\uds a} \uds {\bm v}_s - E_2 \nabla_s^{\uds a} \nabla_s^{\uds a} \xi.											
\end{split}
\eeqn
\end{lemma}

\section{Estimates about Adiabatic Limits}

In this appendix we derive certain estimates about the adiabatic limit of the affine vortex equation over a fixed compact domain. We first fix the notations and set up the problem. Let 
\beqn
\Omega \subset {\mb H}
\eeqn
be an open subset whose closure is compact. Let 
\beqn
\sigma: \Omega \to [1, +\infty)
\eeqn
be a smooth function. Let
\beqn
J_\epsilon \in {\mc J}_{\rm tame}(V)
\eeqn
be a family of $K$-invariant $\omega_V$-tamed almost complex structures which depends on $\epsilon \geq 0$ smoothly.\footnote{By the graph construction, one can assume that $J_\epsilon$ is independent of domain coordinates. See \cite[Appendix A]{Gaio_Salamon_2005} for more details.} When $\epsilon > 0$, for gauged maps ${\bm v} = (u, \phi, \psi)$ from $\Omega$ to $V$, the {\it $\epsilon$-vortex equation} reads 
\beq\label{eqnb1}
\begin{aligned}
\partial_s u + {\mc X}_\phi + J_\epsilon( \partial_t u + {\mc X}_\psi) = 0,\\
\partial_s \psi - \partial_t \phi + [\phi, \psi] + \epsilon^{-2} \sigma \mu(u) = 0,\\
u(\partial \Omega) \subset L_V.
\end{aligned}
\eeq
The limiting case of the above equation when $\epsilon = 0$ is regarded as 
\beq\label{eqnb2}
\begin{aligned}
\partial_s u + {\mc X}_\phi + J_0 ( \partial_t u + {\mc X}_\psi) = 0,\\
\mu(u) \equiv 0,\\
u(\partial \Omega) \subset L_V.
\end{aligned}
\eeq

\begin{notation}
To simplify notations, in this appendix, we abbreviate the base almost complex structure $J_V$ by $J$ and the Levi--Civita connection of the metric $h_V$ by $\nabla$. We still follow the notational convention that the letter $C$ represents a constant whose value can vary from line to line.
\end{notation}

\begin{thm}\label{thmb2}
Let $p \in (2, 4)$. Let $\epsilon_k$ be a sequence of positive numbers converging to zero. Let ${\bm v}_k = (u_k, \phi_k, \psi_k)$ be a sequence of solutions to \eqref{eqnb1} for $\epsilon = \epsilon_k$. Let ${\bm v}_\infty = (u_\infty, \phi_\infty, \psi_\infty)$ be another gauged map from $\Omega$ to $V$ solving \eqref{eqnb2}. Suppose there holds
\beq
\limsup_{k \to \infty} \Big(\| \partial_s u_k + {\mc X}_{\phi_k}(u_k)\|_{L^\infty(\Omega)} + \epsilon_k^{-1} \| \sqrt\sigma \mu(u_k) \|_{L^\infty(\Omega)} \Big) < \infty
\eeq
and $u_k$ converges to $u_\infty$ uniformly over all compact subsets of $\Omega$. Then for any precompact open subset $Q \subset \Omega$, after gauge transforming ${\bm v}_k$, one can write 
\beqn
u_k (z) = \exp_{u_\infty} \xi_k,\ \xi_k \in \Gammait( u_\infty^* TV).
\eeqn
Moreover, if we denote $\beta_k = (\phi_k - \phi_\infty) ds + (\psi_k - \psi_\infty) dt$, then there holds
\beqn
\lim_{k \to \infty} \Big( \| \beta_k \|_{L^p(Q)} + \epsilon_k \| \nabla^{a_\infty} \beta_k \|_{L^p(Q)} + \| \xi_k^H \|_{L^\infty (Q)} + \epsilon_k^{-1} \| \xi_k^G \|_{L^p(Q)} + \|\nabla^{a_\infty} \xi_k \|_{L^p(Q)} \Big) =0.
\eeqn
\end{thm}

\subsection{{\it a priori} estimates}

The first step of the proof is an {\it a priori} estimate for the adiabatic limit of the vortex equation. 

\begin{lemma}(cf. \cite[Lemma 9.3]{Gaio_Salamon_2005})\label{lemmab3}
Let $p\geq 2$ be a real number. For any $M>0$ and any compact subset $Q \subset \Omega$ there exist $C = C(M, \Omega, Q) > 0$ and $\epsilon(M)$ which satisfy the following condition. Suppose $\epsilon \in (0, \epsilon(M)]$ and ${\bm v}$ is a solution to \eqref{eqnb1} over $\Omega$ such that
\beqn
u(\Omega) \subset U_M
\eeqn
and 
\beq\label{eqnb4}
\sup_{z\in \Omega} \left( |{\bm v}_s(z)| +  \epsilon^{-1} \sqrt\sigma |\mu(u(z))| \right) \leq M.
\eeq
Then one has
\begin{align}
\label{eqnb5}\| \mu(u) \|_{L^p(Q)} \leq &\ C \epsilon^{1 + \frac{2}{p}} \left( \| {\bm v}_s \|_{L^2(\Omega)} + \epsilon^{-1} \| \sqrt\sigma \mu(u)\|_{L^2(\Omega)} \right),\\
\| d\mu(u)({\bm v}_s) \|_{L^p(Q)} + \| d\mu(u) ( {\bm v}_t )\|_{L^p(Q)} \leq &\ C \epsilon^{\frac{2}{p}} \left( \| {\bm v}_s \|_{L^2(\Omega)} + \epsilon^{-1} \| \sqrt\sigma \mu(u)\|_{L^2(\Omega)} \right),\label{eqnb6}\\
\| \nabla^a {\bm v}_s \|_{L^p(Q)} + \| \nabla_t^a {\bm v}_s \|_{L^p(Q)} \leq &\ C \epsilon^{-  1 + \frac{2}{p}}.\label{eqnb7}
\end{align}
\end{lemma}

\begin{proof}
In the case when $\Omega$ has empty boundary, one can prove a stronger estimate in which the right hand side of \eqref{eqnb7} has a factor proportional to the square root of the total energy (see \cite[Lemma 9.3]{Gaio_Salamon_2005}). Here we modify the proof of Gaio and Salamon to extend to the case when $\Omega$ has a nonempty boundary. First we prove the $p = 2$ case. Introduce $u_0, v_0: \Omega \to {\mb R}$ given by
\beqn
u_0:= \frac{1}{2} \left( |{\bm v}_s|^2 + \epsilon^{-2} \sigma |\mu(u)|^2 \right),
\eeqn
\beqn
v_0:= \frac{1}{2} \left( | \nabla_s^a {\bm v}_s|^2 + |  \nabla_t^a {\bm v}_s|^2 + \epsilon^{-4}  \sigma^2 |{\mc X}_{\mu(u)}(u)|^2 + \epsilon^{-2} \sigma | d\mu(u)({\bm v}_s)|^2 + \epsilon^{-2} \sigma | d\mu(u)({\bm v}_t) |^2 \right).
\eeqn
Via explicit computation, it was shown in \cite[p.119]{Gaio_Salamon_2005} that there exists $C>0$ such that
\beq\label{eqnb8}
\Delta u_0 \geq v_0 - C u_0.
\eeq
Choose a precompact open subset $Q' \subset \Omega$ containing $Q$ in the interior. 

We use a certain mean value estimate (\cite[Lemma 9.2]{Gaio_Salamon_2005}) to derive the desired estimate for $p =2$. When $\partial \Omega = \emptyset$, the differential inequality \eqref{eqnb8} implies that for some $C>0$, 
\begin{multline}\label{eqnb9}
\epsilon^{-1} \|\mu(u)\|_{L^2(Q)} + \| d\mu(u) ({\bm v}_s) \|_{L^2(Q)} + \| d\mu(u)({\bm v}_t) \|_{L^2(Q)} + \epsilon \|  \nabla_s^a {\bm v}_s \|_{L^2(Q)} + \epsilon \| \nabla_t^a {\bm v}_s \|_{L^2(Q)} \\
\leq C \epsilon \left( \| {\bm v}_s \|_{L^2(Q')} + \epsilon^{-1} \|\mu(u) \|_{L^2(Q')} \right) \leq C\epsilon.
\end{multline}
Here the last inequality follows from the \eqref{eqnb4} and the fact that $Q'$ has finite area. The mean value estimate also implies that 
\beq\label{eqnb10}
\| {\bm v}_s \|_{L^\infty(Q)} + \epsilon^{-1} \| \mu(u) \|_{L^\infty(Q)} \leq C \left(  \| {\bm v}_s \|_{L^2(Q')} + \epsilon^{-1} \| \mu(u) \|_{L^2(Q')} \right) \leq C.
\eeq
To use the mean value estimate for the case when $\partial \Omega \neq \emptyset$, one has to verify that 
\beq\label{eqnb11}
\frac{\partial}{\partial {\bm n}} u_0|_{\partial \Omega} \equiv 0.
\eeq
This is another place where we need the properties of the metric $h_V$. Indeed, 
\beqn
\partial_t |\mu(u)|^2 |_{\partial \Omega} = 2 \langle \partial_t \mu(u), \mu(u) \rangle|_{\partial \Omega} = 0
\eeqn
as $u(\partial \Omega) \subset L_V \subset \mu^{-1}(0)$. On the other hand, to compute $\partial_t |{\bm v}_s|^2$, we change the gauge locally so that $\psi \equiv 0$. So one has ${\bm v}_s + J_\epsilon \partial_t u = 0$. 
\begin{multline*}
\frac{1}{2} \partial_t | {\bm v}_s|^2 = \langle  \nabla_t (\partial_s u + {\mc X}_\phi), {\bm v}_s \rangle
= \langle  \nabla_s ( \partial_t u ) + {\mc X}_{\partial_t \phi} +  \nabla_{\partial_t u} {\mc X}_\phi, {\bm v}_s \rangle\\
= \langle  \nabla_s (J_\epsilon {\bm v}_s), {\bm v}_s \rangle  + \langle {\mc X}_{\partial_t \phi}, {\bm v}_s \rangle + \langle  \nabla_{\partial_t u} {\mc X}_\phi, {\bm v}_s \rangle:= {\rm I} + {\rm II} + {\rm III}.
\end{multline*}
Since $J_\epsilon = J$ over the boundary, $L_V$ is totally geodesic, and $J (TL_V)$ is orthogonal to $TL_V$ (see Lemma \ref{metric}), one has 
\beqn
{\rm I} |_{\partial \Omega} = \langle  \nabla_s (J_\epsilon {\bm v}_s), {\bm v}_s \rangle|_{\partial \Omega} = \langle  \nabla_s (J {\bm v}_s), {\bm v}_s \rangle|_{\partial\Omega} = 0.
\eeqn
By the equation \eqref{eqnb1}, one has
\beqn
{\rm II} |_{\partial \Omega} = \langle {\mc X}_{\partial_t \phi}, {\bm v}_s \rangle|_{\partial \Omega} = \langle \sigma \epsilon^{-2} {\mc X}_{\mu(u)}, {\bm v}_s \rangle|_{\partial \Omega} = 0.
\eeqn
Lastly, notice that ${\rm III}$ only depends on $\partial_t u$ pointwise. Hence at each boundary point $z_0 \in \partial \Omega$, one can extend $\partial_t u(z_0)$ to a $K$-invariant vector field over an open neighborhood of $u(z_0)$, which is still denoted by $\partial_t u$. Hence 
\beqn
{\rm III}|_{\partial \Omega} = \langle \nabla_{\partial_t u}{\mc X}_\phi, {\bm v}_s\rangle|_{\partial \Omega} = \langle  \nabla_{{\mc X}_\phi} (\partial_t u), {\bm v}_s \rangle|_{\partial\Omega} = 0
\eeqn
as ${\mc X}_\phi$ is tangent to $L_V$, $\partial_t u$ is orthogonal to $L_V$, and $L_V$ is totally geodesic with respect to $ \nabla$ (see Lemma \ref{metric}). Therefore \eqref{eqnb11} is true. So one can use the mean value estimate and see that \eqref{eqnb9} and \eqref{eqnb10} hold when $\partial \Omega \neq \emptyset$. Then the $p=2$ case of \eqref{eqnb5}---\eqref{eqnb7} follow from \eqref{eqnb9}.

We derive \eqref{eqnb5}---\eqref{eqnb7} for all $p$ by interpolation. Indeed, by the $p=2$ case of \eqref{eqnb9} and \eqref{eqnb10}, one has
\begin{multline*}
\|\mu(u)\|_{L^p(Q)} = \left( \int_Q |\mu(u)|^p \right)^{\frac{1}{p}} \leq \left( \sup_Q |\mu(u)|\right)^{\frac{p-2}{p}} \left( \int_Q |\mu(u)|^2 \right)^{\frac{1}{p}} \\
\leq C \epsilon^{1+\frac{2}{p}} \big( \| {\bm v}_s\|_{L^2(\Omega)} + \epsilon^{-1} \| \sqrt\sigma \mu(u) \|_{L^2(\Omega)} \big)
\end{multline*}
and
\begin{multline*}
\| d\mu(u) ({\bm v}_s) \|_{L^p(Q)} + \| d\mu(u)({\bm v}_t)\|_{L^p(Q)} = \left( \int_Q |d \mu(u)( {\bm v}_s ) |^p \right)^{\frac{1}{p}} + \left( \int_Q |d\mu(u)({\bm v}_t)|^p \right)^{\frac{1}{p}} \\
 \leq C \left( \sup_Q |{\bm v}_s| \right)^{\frac{p-2}{p}} \left( \| d\mu(u) ({\bm v}_s) \|_{L^2(Q)}^{\frac{2}{p}} + \| d\mu(u)( {\bm v}_t) \|_{L^2(Q)}^\frac{2}{p} \right) \\
 \leq C \epsilon^{\frac{2}{p}} \big( \|{\bm v}_s\|_{L^2(\Omega)} + \epsilon^{-1} \| \sqrt\sigma \mu(u) \|_{L^2(\Omega)} \big).
\end{multline*}
To derive \eqref{eqnb7}, we first prove that 
\beq\label{eqnb12}
\|  \nabla_s^a {\bm v}_s \|_{L^\infty(Q)} + \|  \nabla_t^a {\bm v}_s \|_{L^\infty(Q)} \leq C \epsilon^{-1}.
\eeq
Indeed, let $\varphi_\epsilon: {\mb H} \to {\mb H}$ be the map $z \mapsto \epsilon z$. Then $\varphi_\epsilon$ pulls back ${\bm v}$ to a solution ${\bm v}' = (u', \phi', \psi')$ to the equation 
\beqn
\partial_s u' + {\mc X}_{\phi'} + J_\epsilon (\partial_t u' + {\mc X}_{\psi'}) = 0,\ \partial_s \psi' - \partial_t \phi' + [\phi', \psi'] + (\sigma \circ \varphi_\epsilon) \mu(u') = 0
\eeqn
over $\varphi_\epsilon^{-1}(\Omega)$. The condition \eqref{eqnb4} implies that 
\beqn
| {\bm v}_s' | + \sqrt{\sigma \circ \varphi_\epsilon} |\mu(u')| \leq M \epsilon.
\eeqn
Denote $a' = \phi' ds + \psi' dt$. Then by the elliptic estimate of the vortex equation (notice that $J_\epsilon$ and $\sigma \circ \varphi_\epsilon$ satisfy uniform bounds on all derivatives independent of $\epsilon$), one has 
\beqn
\sup_{\varphi_\epsilon^{-1}(Q)} \Big( | \nabla_s^{a'} {\bm v}_s'| + | \nabla_t^{a'} {\bm v}_s'| \Big) \leq C\epsilon. 
\eeqn
This is equivalent to \eqref{eqnb12}. Therefore 
\begin{multline*}
\| \nabla_s^a {\bm v}_s \|_{L^p(Q)} + \|  \nabla_t^a {\bm v}_s \|_{L^p(Q)} \\
\leq C \left( \|  \nabla_s^a{\bm v}_s \|_{L^\infty(Q)}^{\frac{p-2}{p}} \|  \nabla_s^a {\bm v}_s \|_{L^2(Q)}^{\frac{2}{p}} + \|  \nabla_t^a {\bm v}_s \|_{L^\infty(Q)}^{\frac{p-2}{p}} \|  \nabla_t^a {\bm v}_s \|_{L^2(Q)}^{\frac{2}{p}} \right) \leq C \epsilon^{-1 + \frac{2}{p}}.
\end{multline*}
This finishes the proof. 
\end{proof}

\subsection{Projection to $\mu^{-1}(0)$}

The following theorem is an analogue of \cite[Step 2 of Proof of Theorem 10.1]{Gaio_Salamon_2005}. 

\begin{thm}\label{thmb4}
For any $M>0$ and any compact subset $Q \subset \Omega$, there exist constants $C(M, \Omega, Q)>0$ and $\epsilon(M)>0$ satisfying the following conditions. Let ${\bm v} = (u, \phi, \psi)$ be a solution to the $\epsilon$-vortex equation \eqref{eqnb1} for $\epsilon \in (0, \epsilon(M)]$. Suppose
\beqn
\sup_{z\in \Omega} \left( |{\bm v}_s(z)| + \epsilon^{-1} \sqrt{\sigma(z)} |\mu(u(z))| \right) \leq M.
\eeqn
Then there exists a map $h: \Omega \to {\mf k}$ satisfying the following conditions. 

\begin{enumerate}
\item If we denote $\dot u = \exp_u (J{\mc X}_h)$, then there holds 
\beqn
\mu( \dot u) \equiv 0.
\eeqn

\item Define $\dot \phi$ and $\dot \psi$ via the conditions
\begin{align*}
&\ \partial_s \dot u + {\mc X}_{\dot \phi}(\dot u) \in (TV)_H,\ &\ \partial_t \dot u + {\mc X}_{\dot \psi}(\dot u) \in (TV)_H
\end{align*}
and denote 
\beqn
\dot {\bm \xi}:= (\dot \xi, \dot \beta) = ( J {\mc X}_h, ( \phi - \dot \phi) ds + (\psi - \dot \psi) dt).
\eeqn
If we view $\dot {\bm \xi}$ as an infinitesimal deformation of ${\bm v}$, then  
\beq\label{eqnb13}
\| \dot \beta\|_{L^p(Q)} + \epsilon \| \nabla^a \dot \beta \|_{L^p(Q)} + \epsilon^{-1} \| \dot \xi \|_{L^p(Q)} + \| \nabla^a \dot \xi \|_{L^p(Q)} \leq C \epsilon^{\frac{2}{p}}.
\eeq

\item There holds
\beq\label{eqnb14}
\| \dot{\bm v}_s + J_\epsilon \dot{\bm v}_t \|_{L^p(Q)} \leq C \epsilon^{\frac{2}{p}}.
\eeq
and
\beq\label{eqnb15}
\|  \nabla^a \dot {\bm v}_s \|_{L^p(Q)} + \| \nabla^a \dot {\bm v}_t \|_{L^p(Q)} \leq C \epsilon^{- 1 + \frac{2}{p}}.
\eeq
\end{enumerate}
\end{thm}

\begin{proof}
For given $M$, when $\epsilon$ is sufficiently small, the image of $u$ is sufficiently close to the level set $\mu^{-1}(0)$. Then the pointwise existence and uniqueness of $h$ and $\dot u$ follows from the implicit function theorem. Hence $\dot {\bm v} = (\dot u, \dot \phi, \dot \psi)$ is defined over $\Omega$. Moreover, for a certain $C>0$ there holds throughout $\Omega$ that 
\beq\label{eqnb16}
|J{\mc X}_h| \leq C |\mu(u)|.
\eeq
In the rest of the proof, we derive the expected estimates in several steps. 

\vspace{0.2cm}

\noindent {\it Step 1. There is a constant $C>0$ such that throughout $\Omega$ there holds
\beq\label{eqnb17}
| \nabla^a (J {\mc X}_h)| \leq C \Big( | \mu(u)| + | d\mu(u) ({\bm v}_s)| + | d\mu(u)({\bm v}_t)|  \Big).
\eeq}

\noindent By Lemma \ref{lemmaa2}, one has
\beq\label{eqnb18}
{\mc X}_{\dot\eta} + \dot{\bm v}_s - E_1( {\bm v}_s) - E_2 ( \nabla_s^{a} (J {\mc X}_h) ) 
=  0.
\eeq
Here $E_1 = E_1( u, J {\mc X}_h)$, $E_2 = E_2( u, J{\mc X}_h)$. Since $\dot {u}$ is contained in $\mu^{-1}(0)$, applying $d\mu(\dot u)$ to the above identities gives  
\beqn
d\mu( \dot{u}) ( E_2  (  \nabla_s^{a} (J {\mc X}_h) )= - d\mu ( \dot{u}) (E_1 ({\bm v}_s)),
\eeqn
The left hand side roughly contains the term we would like to bound since $E_2$ is close to the identity.
By the smoothness of $d\mu$, one has
\beqn
| d\mu( \dot{u}) \cdot E_2 ( \nabla_s^a (J {\mc X}_h)) - d\mu(u) ( \nabla_s^a (J {\mc X}_h) ) | \leq C |J {\mc X}_h| | \nabla_s^a (J {\mc X}_h)| \leq C |\mu(u)| |\nabla_s^a J {\mc X}_h|.
\eeqn
\beqn
|d\mu( \dot{u}) ( E_1({\bm v}_s) ) - d\mu(u) ({\bm v}_s) | \leq C |J{\mc X}_h| |{\bm v}_s| \leq C |{\bm v}_s||\mu(u)| \leq C |\mu(u)|.
\eeqn
Moreover, since $J$ is $K$-invariant, one has
\beqn
\nabla_s^a J {\mc X}_h = J \nabla_s^a {\mc X}_h + (\nabla_{{\bm v}_s} J) {\mc X}_h= J {\mc X}_{\nabla_s^a h} + J \nabla_{{\bm v}_s} {\mc X}_h + (\nabla_{{\bm v}_s} J) {\mc X}_h.
\eeqn
Therefore, 
\beqn
|\nabla_s^a J {\mc X}_h - J {\mc X}_{\nabla_s^a h}| \leq C |{\bm v}_s| |{\mc X}_h| \leq C|\mu(u)|.
\eeqn
Moreover, since $|J{\mc X}_{\nabla_s^a h}| \leq C |d\mu(u)(J {\mc X}_{\nabla_s^a h})|$, one has
\begin{multline*}
| \nabla_s^a J {\mc X}_h| \leq C \left( | d\mu(u)( \nabla_s^a (J {\mc X}_h) )| +  |\mu(u)| \right)\\
\leq C \left(  |d\mu(u)( \nabla_s^a (J{\mc X}_h)) - d\mu( \dot u) ( E_2( \nabla_s^a (J {\mc X}_h)))| + |d\mu( \dot u)( E_1({\bm v}_s) - d\mu(u) ({\bm v}_s)| + |d\mu(u)({\bm v}_s)| + | \mu(u) | \right)\\
\leq C \left( |\mu(u)| | \nabla_s^a (J{\mc X}_h)|  + |\mu(u)| + |d\mu(u)({\bm v}_s)|\right).
\end{multline*}
Therefore, when $|\mu(u)|$ is sufficiently small, one has
\beqn
| \nabla_s^a J{\mc X}_h| \leq C \Big( |d\mu(u) ( {\bm v}_s )| + |\mu(u)|\Big).
\eeqn
Similarly one can derive
\beqn
| \nabla_t^a J{\mc X}_h| \leq C \Big( |d\mu(u) ({\bm v}_t)| + |\mu(u)|\Big).
\eeqn
Therefore \eqref{eqnb17} follows.

\vspace{0.2cm}

\noindent {\it Step 2. There is a constant $C>0$ such that throughout $\Omega$ there holds
\beq\label{eqnb19}
|\dot{\bm v}_s + J_\epsilon \dot{\bm v}_t |\leq C \Big( | \mu(u) | + | d\mu(u)({\bm v}_s)| + | d\mu(u) ({\bm v}_t) |\Big). 
\eeq
Then \eqref{eqnb14} follows from \eqref{eqnb19}, \eqref{eqnb5}, and \eqref{eqnb6}.
}
\vspace{0.1cm}

\noindent 
Since $\dot{\bm v}_s$ and $\dot{\bm v}_t$ are contained in $(TV)_H$, \eqref{eqnb18} implies that 
\beqn
\begin{split}
\dot{\bm v}_s + J_\epsilon \dot{\bm v}_t = &\ P_H \big( E_1({\bm v}_s) + E_2( \nabla_s^a J{\mc X}_h) + J_\epsilon E_1({\bm v}_t) + J_\epsilon E_2(  \nabla_t^a J{\mc X}_h)\big) \\
                           = &\ P_H  \big( (J_\epsilon E_1 - E_1 J_\epsilon )({\bm v}_t) + E_2( \nabla_s^a J{\mc X}_h) + J_\epsilon E_2( \nabla_t^a J{\mc X}_h) \big).
\end{split}
\eeqn
Then \eqref{eqnb5} and \eqref{eqnb17} one has
\beqn
| \dot{\bm v}_s + J_\epsilon \dot{\bm v}_t | \leq C \Big( |J{\mc X}_h| |{\bm v}_t| + | \nabla_s^a J{\mc X}_h| + |  \nabla_t^a J {\mc X}_h | \Big) \leq C \Big( |\mu(u)| + |d\mu(u)({\bm v}_s)| + |d\mu(u)({\bm v}_t)|\Big).
\eeqn
So \eqref{eqnb19} follows. 

\vspace{0.2cm}

\noindent {\it Step 3. There is a constant $C>0$ such that throughout $\Omega$ there holds
\beq\label{eqnb20}
|\dot \beta| \leq C \Big( |\mu(u)| + |d \mu(u)({\bm v}_s)| + |d\mu(u)({\bm v}_t)| \Big).
\eeq
}

\noindent 
Since $\dot{\bm v}_s$ is contained in $(TV)_H$, by applying $d\mu(\dot u) \circ J$ to \eqref{eqnb18}, one obtains
\beqn
d\mu(\dot u)(J {\mc X}_{\dot \eta}) = d\mu(\dot u)(J E_1({\bm v}_s)) + d\mu(\dot u)(J E_2 \nabla_s^a (J {\mc X}_h)).
\eeqn
It follows that
\begin{multline*}
| d\mu(\dot u)(J {\mc X}_{\dot \eta} ) | \leq C \left( | d\mu(\dot u)(J E_1({\bm v}_s)) | + | d\mu(\dot u)(J E_2 \nabla_s^a(J {\mc X}_h)) | \right)\\
\leq C \left( | d\mu(u) (J {\bm v}_s ) | + |J{\mc X}_h||{\bm v}_s| + | d\mu(u) (J  \nabla_s^a (J {\mc X}_h))| + |J{\mc X}_h||\nabla_s^a (J{\mc X}_h)| \right) \\
			 \leq  C \left( |  d\mu(u)( {\bm v}_t) | + |\mu(u)| + | \nabla_s^a J {\mc X}_h| \right).
			 \end{multline*}
One can derive the estimate for $\dot\zeta$. Notice that $|\dot \eta|$ is comparable to $|d\mu(\dot u)(J {\mc X}_{\dot \eta})|$. Then \eqref{eqnb20} follows from the above estimate and \eqref{eqnb17}.

\vspace{0.2cm}

\noindent {\it Step 4. There exists a constant $C>0$ such that throughout $\Omega$ there holds
\beq\label{eqnb21}
 |  \nabla_s^a \dot \beta | \leq C \Big( |{\bm v}_s| + |  \nabla^a {\bm v}_s| + |J {\mc X}_h| + | \nabla^a J {\mc X}_h| + |\dot \beta|\Big).
\eeq
Then \eqref{eqnb13} follows from \eqref{eqnb16}, \eqref{eqnb17}, \eqref{eqnb20}, \eqref{eqnb21}, and \eqref{eqnb5}---\eqref{eqnb7}.
}

\vspace{0.1cm}

\noindent We first derive an estimate for second order derivatives of $h$. Apply $\nabla_s^a \circ d\mu(\dot u)$ to \eqref{eqnb18}, one obtains
\beqn
\nabla_s^a d\mu(\dot u)(E_2 \nabla_s^a J {\mc X}_h) + \nabla_s^a d\mu(\dot u)(E_1({\bm v}_s)) = 0.
\eeqn
Since the maps $E_1$, $E_2$, the almost complex structure, and the connection $\nabla$ are all $K$-invariant, one can see that for a certain $C>0$ there holds
\beqn
|\nabla_s^a d\mu(\dot u) (E_2 \nabla_s^a J {\mc X}_h) - d\mu(\dot u) ( E_2 \nabla_s^a \nabla_s^a J {\mc X}_h ) | \leq C | \nabla_s^a J {\mc X}_h| \left( |{\bm v}_s | + |\nabla_s^a J {\mc X}_h| \right) \leq C |\nabla_s^a J {\mc X}_h|,
\eeqn
and 
\beqn
|\nabla_s^a d\mu(\dot u) (E_1({\bm v}_s)) - d\mu(\dot u) ( E_1 ( \nabla_s^a {\bm v}_s)) | \leq C |{\bm v}_s| \left( |{\bm v}_s| + |\nabla_s^a J {\mc X}_h| \right)\leq C |{\bm v}_s|.
\eeqn
Therefore, one obtains 
\beqn
| d\mu(\dot u) ( E_2 \nabla_s^a \nabla_s^a J {\mc X}_h  )| \leq C (|{\bm v}_s| + |\nabla^a J {\mc X}_h|).
\eeqn
Furthermore, 
\beqn
|d\mu(\dot u) E_2 \nabla_s^a \nabla_s^a J {\mc X}_h - d\mu(u) \nabla_s^a \nabla_s^a J {\mc X}_h | \leq C | J {\mc X}_h| |\nabla_s^a \nabla_s^a J {\mc X}_h|
\eeqn
and
\beqn
| \nabla_s^a \nabla_s^a J {\mc X}_h - J {\mc X}_{ \nabla_s^a \nabla_s^a h} | \leq C \left( |{\bm v}_s| |{\mc X}_h| + |\nabla_s^a {\bm v}_s| |{\mc X}_h| \right);
\eeqn
also $|J {\mc X}_{\nabla_s^a \nabla_s^a h}|$ is comparable to $|d\mu(u) (J {\mc X}_{\nabla_s^a \nabla_s^a h})|$. So
\begin{multline}\label{eqnb22}
|\nabla_s^a \nabla_s^a J {\mc X}_h| \leq C \left( |d\mu(u)( \nabla_s^a \nabla_s^a J {\mc X}_h)| + |{\bm v}_s||{\mc X}_h| + |\nabla_s^a {\bm v}_s | |{\mc X}_h| \right)\\
\leq C \left( |J {\mc X}_h| + |\nabla_s^a J {\mc X}_h| +  |{\bm v}_s| +  |\nabla_s^a {\bm v}_s|  \right).
\end{multline}

Now we estimate the derivative of $\dot \eta$. By applying $\nabla_s^a \circ d\mu(\dot u) \circ J$ to \eqref{eqnb18}, one obtains
\beqn
\nabla_s^a d\mu(\dot u)(J {\mc X}_{\dot\eta} ) = \nabla_s^a d\mu(\dot u)(J E_1({\bm v}_s) ) + \nabla_s^a d\mu(\dot u)(JE_2 \nabla_s^a(J {\mc X}_h)).
\eeqn
Then by the previous bounds, particularly the second order estimate \eqref{eqnb22}, one has
\beqn
\begin{split}
|\nabla_s^a \dot \eta|\leq &\ C | d\mu(\dot u)( J {\mc X}_{\nabla_s^a \dot \eta})|  \\
\leq &\ C \big( |d\mu(\dot u)(J \nabla_s^a {\mc X}_{\dot \eta})| + |d\mu(\dot u) (J \nabla_{{\bm v}s} {\mc X}_{\dot \eta})| \big) \\
\leq &\ C \big( |\nabla_s^a d\mu(\dot u)(J {\mc X}_{\dot \eta})| + |{\bm v}_s| |\dot \eta| \big) \\
\leq &\ C \big( |\nabla_s^a d\mu(\dot u) (J E_1({\bm v}_s))| + |\nabla_s^a d\mu(\dot u)( J E_2 \nabla_s^a J {\mc X}_h)| + |\dot \beta| \big)\\
\leq &\ C \big( |{\bm v}_s|^2 + |{\bm v}_s| |\nabla_s^a J {\mc X}_h| + |\nabla_s^a {\bm v}_s| + |{\bm v}_s| |\nabla_s^a  J {\mc X}_h| + |\nabla_s^a J {\mc X}_h|^2 + |\nabla_s^a \nabla_s^a J {\mc X}_h|+ |\dot \beta| \big)\\
\leq &\ C \big( |{\bm v}_s| + |\nabla^a {\bm v}_s| + |J{\mc X}_h| + |\nabla^a J {\mc X}_h| + |\dot \beta| \big).
\end{split}
\eeqn
The estimates of $\nabla_t^a \dot \eta$, $\nabla_s^a \dot \zeta$, and $\nabla_t^a \dot \zeta$ can be derived similarly.

\vspace{0.2cm}

\noindent {\it Step 5. One has
\beq\label{eqnb23}
|\nabla^a \dot{\bm v}_s| \leq C \left( |{\bm v}_s| + |\nabla^a {\bm v}_s| + |\dot \beta| +  |\nabla^a \dot \beta| + |\nabla^a J {\mc X}_h | \right).
\eeq
Then \eqref{eqnb15} follows from the \eqref{eqnb5}---\eqref{eqnb7}, \eqref{eqnb17}, and \eqref{eqnb21}. }

\vspace{0.1cm}

\noindent By applying $\nabla_s^a$ to \eqref{eqnb18}, one has
\beqn
\nabla_s^a \dot{\bm v}_s = - \nabla_s^a {\mc X}_{\dot\eta} +  \nabla_s^a E_1({\bm v}_s) + \nabla_s^a E_2( \nabla_s^a J {\mc X}_h).
\eeqn
Then \eqref{eqnb23} can be derived from estimates obtained in previous steps. The case of $\nabla_t^a \dot {\bm v}_s$ is similar.
\end{proof}

\subsection{Proof of Theorem \ref{thmb2}}

Now we prove main result of this appendix. The first step is to put each ${\bm v}_k$ in a suitable gauge. By using the implicit function theorem, one can prove that for large $k$, after an appropriate gauge transformation, one can write 
\beqn
\dot u_k = \exp_{u_\infty} \ddot \xi_k
\eeqn
with 
\beqn
\ddot \xi_k \in \Gammait(u_\infty^* (TV)_H).
\eeqn
Notice that by the property of the metric $h_V$ (see Lemma \ref{metric}), one has
\beqn
\ddot \xi_k \in u^* (TV)_H \Longleftrightarrow \ddot \xi_k \in ( {\rm Im} L_{u_\infty})^\bot \Longleftrightarrow d\mu(u_\infty)(J \ddot \xi_k ) = 0.
\eeqn

\begin{lemma}\label{lemmab9}
For any compact subset $Q \subset \Omega$, one has
\beqn
\lim_{k \to \infty} \Big( \| \ddot\xi_k \|_{L^p(Q)} + \| \nabla^{a_\infty} \ddot\xi_k \|_{L^p(Q)}\Big) = 0.
\eeqn
\end{lemma}

\begin{proof}
Using \eqref{eqnb14} and the elliptic estimate for $\ov\partial$ operator (with totally real boundary conditions). Notice that the boundary restriction of $\dot u_k$ agrees with $u_k$ and since $L_V$ is totally geodesic, $\ddot \xi_k|_{\partial \Omega} \subset u_\infty^* TL_V$. The detail is left to the reader.
\end{proof}

Denote 
\beqn
\ddot \beta_k = \ddot \eta_k ds + \ddot \zeta_k dt:= \dot a_k - a_\infty. 
\eeqn

\begin{lemma}\label{lemmab6}
For any compact subset $Q \subset \Omega$ there is a constant $C>0$ such that
\beqn
\| \ddot \beta_k \|_{L^\infty(Q)} + \| \nabla^{a_\infty} \ddot \xi_k \|_{L^\infty(Q)} \leq C.
\eeqn
\end{lemma}

\begin{proof}
Suppress the index $k$ and abbreviate $u_\infty$ by $\uds u$. By the relation $\dot u = \exp_{\uds u} \ddot \xi$, using Lemma \ref{lemmaa2}, one has
\beq\label{eqnb24}
{\mc X}_{\ddot \eta} = \dot{\bm v}_s - E_1(\uds{\bm v}_s) - E_2( \nabla_s^{\uds a} \ddot \xi).
\eeq
Apply $d\mu(\uds u) \circ J \circ E_2^{-1}$ to both sides, one obtains
\beqn
d\mu( \uds u ) (J E_2^{-1} {\mc X}_{\ddot \eta}) = d\mu( \uds u)(J E_2^{-1}(\dot {\bm v}_s - E_1(\uds {\bm v}_s))) - d\mu(\uds u)(J \nabla_s^{\uds a} \ddot \xi).
\eeqn
The first term on the right hand side is uniformly bounded. Moreover, since $d\mu(\uds u)(J \ddot \xi) = 0$, the last term above is bounded by
\beqn
| d\mu( \uds u )(J \nabla_s^{\uds a} \ddot \xi)| \leq C |\uds {\bm v}_s| |\ddot \xi|.
\eeqn
By Lemma \ref{lemmab9} and the Sobolev embedding $W^{1,p} \hookrightarrow C^0$, this term is uniformly bounded. Therefore one has
\beqn
|d\mu(\uds u)(J E_2^{-1} {\mc X}_{\ddot \eta})| \leq C.
\eeqn
Moreover, since
\beqn
|d\mu(\uds u)(J E_2^{-1} {\mc X}_{\ddot \eta}) - d\mu(\uds u)(J {\mc X}_{\ddot \eta}) | \leq C |\ddot \eta||\ddot \xi|\leq C
\eeqn
and $|\ddot \eta|$ is comparable to $|d\mu(\uds u)(J {\mc X}_{\ddot \eta})|$, one obtains a uniform bound on $\ddot \eta$. Then using \eqref{eqnb24} again, one obtains a uniform bound on $|\nabla_s^{\uds a} \ddot \xi|$. It is similar to bound $|\ddot \zeta|$ and $|\nabla_t^{\uds a} \ddot \xi|$.
\end{proof}

Now we estimate the distance of the connection parts of ${\bm v}_\infty$ and ${\bm v}_k$.

\begin{lemma}\label{lemmab7} 
For all compact subsets $Q \subset \Omega$ there holds
\beqn
\lim_{k \to \infty} \left( \| a_k - a_\infty \|_{L^p(Q)} + \epsilon_k \| \nabla^{a_\infty} (a_k - a_\infty ) \|_{L^p(Q)} \right) = 0.
\eeqn
\end{lemma}

\begin{proof}
Since $\dot a_k$ and $a_\infty$ are pulled back by $\dot u_k$ and $u_\infty$, and $\dot u_k = \exp_{u_\infty} \dot \xi_k$, by using Lemma \ref{lemmab9}, one has
\beq\label{eqnb25}
\lim_{k \to \infty} \| \dot a_k - a_\infty \|_{L^p(Q)} \leq C  \lim_{k \to \infty} \Big( \| \dot \xi_k \|_{L^p(Q)} + \| \nabla^{a_\infty} \dot \xi_k \|_{L^p(Q)}\Big) = 0.
\eeq
Together with \eqref{eqnb13} (notice that $p< 4$ so $-1 + \frac{4}{p}> 0$) one has
\beq\label{eqnb26}
\lim_{k \to \infty} \| a_k - a_\infty \|_{L^p(Q)} = 0.
\eeq
To estimate $\nabla^{a_\infty} (a_k - a_\infty)$, we separate
\beqn
\nabla^{a_\infty} (a_k - a_\infty) = \nabla^{a_\infty} (a_k - \dot a_k) + \nabla^{a_\infty} (\dot a_k - a_\infty) = \nabla^{a_\infty} \dot \beta_k + \nabla^{a_\infty} \ddot \beta_k.
\eeqn
Then by \eqref{eqnb13}, \eqref{eqnb20}, and \eqref{eqnb26}, one has
\beq\label{eqnb27}
\lim_{k \to \infty} \epsilon_k \| \nabla^{a_\infty} \dot \beta_k \|_{L^p(Q)} \leq C \lim_{k \to \infty} \epsilon_k \left( \| \nabla^{a_k} \dot \beta_k \|_{L^p(Q)} + \| \dot \beta_k \|_{L^\infty(Q)} \| a_k -a_\infty \|_{L^p(Q)} \right) = 0.
\eeq

To estimate the $L^p$ norm of $\nabla^{a_\infty} \ddot \beta_k$, we only do it for $\nabla_t^{a_\infty} \ddot \eta$. The estimates for other components can be proved in a similar fashion. By Lemma \ref{lemmaa3}, one has the following calculation. (We suppress the index $k$, denote $\uds {\bm v} = {\bm v}_\infty$, and abbreviate $E_i (u_\infty, \ddot \xi)$, $E_{ij}(u_\infty, \ddot \xi)$ by $E_i$ and $E_{ij}$ respectively.)
\begin{multline}\label{eqnb28}
L_{\dot u} \nabla_t^{\dot a} \ddot \eta = \nabla_t^{\dot a} \dot {\bm v}_s + \nabla_{{\mc X}_{\ddot\eta}} {\mc X}_{\ddot\zeta} - \nabla_{\dot {\bm v}_t} {\mc X}_{\ddot \eta} - \nabla_{\dot {\bm v}_s} {\mc X}_{\ddot \zeta}  -E_{11} (\uds {\bm v}_s, \uds {\bm v}_t) - E_{12} ( \uds {\bm v}_s, \nabla_t^{\uds a} \ddot \xi) \\
                             - E_{21}(\nabla_s^{\uds a } \ddot \xi, \uds {\bm v}_t ) - E_{22} (\nabla_s^{\uds a} \ddot \xi, \nabla_t^{\uds a} \ddot \xi) - E_1 \nabla_t^{\uds a} \uds {\bm v}_s - E_2 \nabla_t^{\uds a}\nabla_s^{\uds a} \ddot \xi.
                             \end{multline}
Since $|L_{\dot u} \nabla_t^{\dot a} \ddot \eta|$ is comparable to $|\nabla_t^{\dot a} \ddot \eta|$. It suffices to bound each term on the right hand side above. 

\begin{enumerate}

\item Since $E_{ij}$ are bounded linear maps, by Lemma \ref{lemmab6} and the boundedness of $|\uds {\bm v}_s|$, $|\nabla^{\uds a} \uds{\bm v}_s|$, $|\dot {\bm v}_s|$, all terms on the right hand side of \eqref{eqnb28} have uniformly bounded $L^p$-norms over a compact $Q \subset \Omega$ except for the terms 
\begin{align*}
&\ \nabla_t^{\uds a} \dot {\bm v}_s,\ &\  E_2 \nabla_t^{\uds a} \nabla_s^{\uds a} \ddot \xi.
\end{align*}

\item By \eqref{eqnb15}, one has
\beqn
\| \nabla_t^{\dot a} \dot{\bm v}_s \|_{L^p(Q)} \leq C \epsilon^{- 2 + \frac{4}{p}}.
\eeqn

\item Notice that $|L_{\dot u} \nabla_t^{\dot a} \ddot \eta|$ is comparable to $|d\mu(\dot u)(J L_{\dot u}\nabla_t^{\dot a} \ddot \eta)|$. Furthermore, 
\beqn
|d\mu(\dot u)(J L_{\dot u} \nabla_t^{\dot a} \ddot \eta) - d\mu(\uds u)(J E_2^{-1} L_{\dot u} \nabla_t^{\dot a} \ddot \eta)| \leq C|\ddot \xi||\nabla_t^{\dot a} \ddot \eta|
\eeqn
and $|\ddot \xi|$ is very small. Hence for the last term $E_2 \nabla_t^{\uds a} \nabla_s^{\uds a} \ddot \xi$ of \eqref{eqnb28}, it suffices to consider instead 
\beqn
d\mu(\uds u)(J E_2^{-1} E_2 \nabla_t^{\uds a} \nabla_s^{\uds a} \ddot \xi) = d\mu( \uds u) (J \nabla_t^{\uds a} \nabla_s^{\uds a} \ddot \xi)\eeqn
Then since $d\mu(\uds u)(J \ddot \xi) = 0$, one has
\begin{multline}
|d\mu( \uds u )(J \nabla_t^{\uds a} \nabla_s^{\uds a} \ddot \xi)| \leq C \left(  | \uds{\bm v}_t| |\nabla_s^{\uds a} \ddot \xi| + |\nabla_t^{\uds a}( d\mu(\uds u)(J \nabla_s^{\uds a} \ddot \xi))| \right) \\
\leq C + C | \nabla_t^{\uds a} [ \nabla_s^{\uds a}, d\mu(\uds u) J] \ddot \xi | \leq C + C |\nabla_t^{\uds a} \uds {\bm v}_s||\ddot \xi| + C |\uds{\bm v}_s| |\nabla_t^{\uds a} \ddot \xi| \leq C.
\end{multline}
\end{enumerate}

Then since $p < 4$, one has 
\beqn
\lim_{k \to \infty} \epsilon_k \| \nabla_t^{\dot a_k} \ddot \eta_k \|_{L^p(Q)} = 0.
\eeqn
Other components of the derivatives of $\ddot \beta_k$ can be obtained similarly. Hence 
\beqn
\lim_{k \to \infty} \epsilon_k \| \nabla^{\dot a_k} \ddot \beta_k \|_{L^p(Q)} = 0.
\eeqn
Since $\dot a_k  - a_\infty = \ddot \beta_k$ is uniformly bounded, it follows that 
\beqn
\lim_{k \to \infty} \epsilon_k \| \nabla^{a_\infty} \ddot \beta_k \|_{L^p(Q)} \leq C \lim_{k \to \infty} \epsilon_k \left( \| \nabla^{\dot a_k} \ddot \beta_k \|_{L^p(Q)} + \| \ddot \beta_k \|_{L^\infty(Q)} \| \ddot \beta_k \|_{L^p(Q)} \right) = 0.
\eeqn
The lemma follows from the above and \eqref{eqnb27}.
\end{proof}

Now we have to estimate the difference between the maps $u_k$ and $u_\infty$. The uniform convergence $u_k \to u_\infty$ allows one to write $u_k = \exp_{u_\infty} \xi_k$. It follows from previous construction that
\beqn
u_k = \exp_{u_\infty} \xi_k = - \exp_{\dot u_k} J {\mc X}_{h_k} = - \exp_{\exp_{u_\infty} \dot\xi_k} J {\mc X}_{h_k}.
\eeqn
Because of the nonlinearity of the exponential map, $\xi_k$ may not agree with $\dot \xi_k - J {\mc X}_{h_k}$ but the difference is a higher order term. One can define a family of maps
\beqn
\Phi_z: (TV)_H|_{u_\infty(z)} \oplus {\mf k} \to T_{u_\infty(z)} V
\eeqn
such that 
\beqn
\xi_k (z) - ( \dot \xi_k (z) - J {\mc X}_{h_k (z)}(u_\infty(z)) ) = \Phi_z ( \dot\xi_k (z), h_k (z)). 
\eeqn
Moreover, for a certain $C>0$ there holds for all large $k$ that
\beqn
|\Phi_z(\dot\xi_k, h_k) |\leq C |\dot\xi_k | |h_k |
\eeqn
and
\beq\label{eqnb30}
| \nabla^{a_\infty} \Phi_z (\dot\xi_k, h_k)| \leq C \Big( | {\bm v}_{\infty, s}| |\dot\xi_k ||h_k | + |  \nabla^{a_\infty} \dot\xi_k ||h_k | + |\dot\xi_k  | | \nabla^{a_\infty} h_k | \Big).
\eeq
Therefore, by Lemma \ref{lemmab3} and Lemma \ref{lemmab9}, one has
\begin{multline*}
\lim_{k \to \infty} \Big( \| \xi_k^H \|_{L^p(Q)}  + \epsilon_k^{-1} \| \xi_k^G \|_{L^p(Q)} \Big) \\
 \leq C \lim_{k \to \infty} \Big(  \| \dot \xi_k \|_{L^p(Q)} + \epsilon_k^{-1} \| J {\mc X}_{h_k}\|_{L^p(Q)} + \epsilon_k^{-1} \| \Phi_z( \dot \xi_k, h_k ) \|_{L^p(Q)} \Big) \\
 \leq C \lim_{k \to \infty} \Big(  \| \dot \xi_k \|_{L^p(Q)} + \epsilon_k^{-1} \| J {\mc X}_{h_k}\|_{L^p(Q)} \Big) = 0.
\end{multline*}
By the fact that $a_\infty - a_k$ is uniformly bounded over a compact $Q$ and the estimate \eqref{eqnb30}, one has
\begin{multline*}
\lim_{k \to \infty} \|  \nabla^{a_\infty} \xi_k \|_{L^p(Q)} \leq \lim_{i \to \infty} \Big( \|  \nabla^{a_\infty} \dot \xi_k \|_{L^p(Q)} + \|  \nabla^{a_\infty} J{\mc X}_{h_k} \|_{L^p(Q)} + \|  \nabla^{a_\infty} \Phi_z (\dot \xi_k, h_k) \|_{L^p(Q)} \Big)\\
\leq \lim_{k \to \infty} \Big( \|  \nabla^{a_k} J {\mc X}_{h_k} \|_{L^p(Q)} + C \| a_\infty - a_k \|_{L^\infty(Q)} \| h_k \|_{L^p(Q)} + \|  \nabla^{a_\infty} \Phi_z (\dot \xi_k, h_k) \|_{L^p(Q)} \Big) = 0.
\end{multline*}
This finishes the proof of Theorem \ref{thmb2}.

\bibliography{mathref}

\bibliographystyle{amsplain}

\end{document}